\documentclass[10pt]{article}
\usepackage{latexsym}
\usepackage{hyperref}
\usepackage{amssymb}
\usepackage{amsmath}
\usepackage{mathrsfs}
\usepackage{graphicx}
\usepackage{verbatim}
\usepackage[all]{xy}

\newtheorem{theo}{Theorem}[section]

\newtheorem{lemma}[theo]{Lemma}

\newtheorem{corollary}[theo]{Corollary}
\newtheorem{prop}[theo]{Proposition}

\newtheorem{conjecture}[theo]{Conjecture}
\newenvironment{proof}{ \textbf{Proof.}}{$\Box$\medskip}
\newtheorem{defi}{Definition}[section]

\newcommand {\ZZ} {\mathbb {Z}}

\newcommand {\CC} {\mathbb {C}}

\newcommand {\QQ} {\mathbb {Q}}

\newcommand{\id}{\mathrm{id}}

\newcommand{\binomial}[2]{\left(\begin{array}{c}#1 \\ #2\end{array} \right)}

\newcommand{\doublebrace}[4]{\left\{\begin{array}{ll}#1 & #2 \\ #3 & #4\end{array}\right.}
\newcommand{\refth}[1]{Theorem \ref{#1}}
\newcommand{\refle}[1]{Lemma \ref{#1}}
\newcommand{\refsec}[1]{Section \ref{#1}}
\newcommand{\refcor}[1]{Corollary \ref{#1}}
\newcommand{\refprop}[1]{Proposition \ref{#1}}
\newcommand{\refcon}[1]{Conjecture \ref{#1}}

\newcommand{\refeq}[1]{(\ref{#1})}
\newcommand{\refdef}[1]{Definition \ref{#1}}

\newcommand{\gl}{\mathrm{gl}}
\newcommand{\sL}{\mathrm{sl}}

\newcommand{\eps}{\varepsilon}

\newcommand{\QL}{\QQ[x_1,\dots,x_n,\frac{1}{x_1},\dots,\frac{1}{x_n}][\frac{1}{1-x^\alpha}]_{\alpha \in \ZZ^n}}
\newcommand{\QQLaurent}{\QQ[x_1,\dots,x_n,\frac{1}{x_1},\dots,\frac{1}{x_n}]}
\def\cplus{\hbox{$\subset${\raise0.3ex\hbox{\kern -0.55em ${\scriptscriptstyle +}$}}\ }}

\def\clplus{\hbox{$\subset${\raise0.3ex\hbox{\kern -0.55em ${\scriptscriptstyle +}$}}\ }}
\def\crplus{\hbox{$\supset${\raise1.05pt\hbox{\kern -0.55em ${\scriptscriptstyle +}$}}\ }}
\def\nperp{\hbox{$\perp${\raise1.05pt\hbox{\kern -0.55em $/$}}\ }}
\def\sperp{\hbox{$\perp${\raise -1.5pt\hbox{\kern -.77em $-$}}\ }}

\author{Todor Milev}
\title{
Partial fraction decompositions and an algorithm for computing the vector partition function
}
\begin{document}
\maketitle
\maketitle
\abstract{ This paper gives an exposition of well known results on vector partition functions. The exposition is based on works of M. Brion, A. Szenes and M. Vergne and is geared toward explicit computer realizations. In particular, the paper presents two algorithms for computing the vector partition function with respect to a finite set of vectors $I$ as a quasipolynomial over a finite set of pointed polyhedral cones. We use the developed techniques to relate a result of P. Tumarkin and A. Felikson (and present an independent proof in the particular case of finite-dimensional root systems) to give bounds for the periods of the Kostant partition functions of $E_6$, $E_7$, $E_8$, $F_4$, $G_2$ (the periods are divisors of respectively 6, 12, 60, 12, 6). 

The first of the described algorithms has been realized and is publicly available under the Library General Public License v3.0 at \\ \url{http://vectorpartition.sourceforge.net/}. We include (non-unique) partial fraction decompositions for the generating functions of the Kostant partition function for $A_2$, $A_3$, $A_4$, $B_2$, $B_3$, $C_2$, $C_3$, $G_2$ in the appendix. }
\section{Introduction and notation}

Given a finite set $I$ of non-zero integral vectors with nonnegative coordinates (with respect to some basis), and a vector $\gamma$ with coordinates $(\gamma^1,\dots,\gamma^n)$, the vector partition function $P_I(\gamma)$ is by definition the number of ways we can split $\gamma$ as an integral sum with non-negative coefficients of the vectors in $I$. It is well known  that there exist finitely many pointed (i.e. with walls passing through the origin) polyhedral cones, with walls parallel to hyperplanes spanned by subsets of $I$ of rank $n-1$, such that $P_I$ is a quasipolynomial over each as a function of the coordinates $(\gamma^1,\dots,\gamma^n)$. We will reprove these results (\refcor{corChambers}); for original proofs, history, and further information on the subject we direct the reader to \cite{ConciniProcesiVergne}, \cite{Sturmfels} and the references therein. 

Our definition of quasipolynomial (\refsec{secRealization}) is done via quasinumbers (=linear combinations of indicator functions of affine translates of lattices) and is slightly different from the commonly accepted one
. The commonly accepted definition of quasipolynomials uses holomorphic functions of the form $e^{2\pi i \sum a_i\gamma^i}$ and we relate it to ours in the remarks after \refdef{defQNring} and in the appendix (indicator functions are not holomorphic).

In \refsec{secWallDist} we present an elementary algorithm for computing the vector partition function as a quasipolynomial. The algorithm allows us to prove that the partition function is a quasipolynomial over a finite set of closed pointed polyhedral cones (\refth{thMain}), but doesn't prove that the walls of the closed pointed polyhedral cones are parallel to hyperplanes spanned by subset of rank $n-1$ of the starting set $I$. 

In Sections \ref{secPFdecomposition} and \ref{secBVdecomposition} we present another algorithm, based on the Szenes-Vergne formula (\cite[Lemma 1.8]{SzenesVergne:ResidueFlas}). The algorithm is implied in \cite{SzenesVergne:ResidueFlas} but is not formulated in the exact same fashion as here; the current exposition is geared towards explicit computer realization. A similar approach is taken in  \cite{DeConciniProcesi:GeometryToricArrangements}. The walls of the chambers of quasipolinomiality produced by this algorithm are parallel to the hyperplanes spanned by subsets of $I$ of rank $n-1$, but a priori may fail to pass through the origin. However, two quasipolynomials coincide if they take the same value over an affine translate of a pointed polyhedral cone of non-zero $n$-volume (\refle{leTheQuasipolynomialLemma}) and so the combinatorial chambers produced by the algorithm in Sections \ref{secPFdecomposition} and \ref{secBVdecomposition} are ``glued'' by those produced by the  \refsec{secWallDist} algorithm and the other way round. This proves \refcor{corChambers} (see \cite{BBCV}). In practice, the algorithm in Sections \ref{secPFdecomposition} and \ref{secBVdecomposition} has turned out much faster than the one in \refsec{secWallDist}.

The algorithm in Sections \ref{secPFdecomposition} and \ref{secBVdecomposition} is based on the notion of a partial fraction decomposition in the ring $\QQ[x_1,\dots,x_n,$ $\frac{1}{x_1}, \dots, \frac{1} {x_n}]$ $[\frac{1}{1-x^\alpha}]_{\alpha \in \ZZ^n}$ (\refdef{defPartFracDecomposition}). The terminology ``partial fraction decomposition'' appears in the same context in \cite{SzenesVergne:ResidueFlas} and  \cite{DeConciniProcesi:GeometryToricArrangements}, however our definition is different to avoid computations with complex roots of unity (which eases computer realization of the algorithm). In \refsec{secTheTheoreticalAlgorithm} we show how to use the elongated Szenes-Vergne formula (\refle{leFlas}(c)) to produce a partial fraction decomposition. In \refsec{secBVdecomposition} we show how to obtain the power series expansion of the decomposed partial fractions by computing explicitly the form specified by \cite[Theorem 1]{BrionVergne}. 

The algorithm in Sections \ref{secPFdecomposition} and \ref{secBVdecomposition} allows us to give upper bounds on the periods of the Kostant partition functions of $E_6$, $E_7$, $E_8$, $F_4$, $G_2$ by proving that they are divisors of respectively 6, 12, 60, 12 and 6 (\refprop{propPeriods}). The periods for the classical root systems were first computed in \cite{BBCV}, and are 1 for type $A$ and 2 for the other series. \refprop{propPeriods} is based on the key \refle{leCanChooseCoeffOne}, which was first proved using different techniques in \cite{Tumarkin:SimpleVertex}. 

\refsec{secTheActualAlgorithm} presents a modification of \refsec{secTheTheoreticalAlgorithm} applicable to partial fraction decompositions of the generating function of the Kostant partition function of the classical root systems of type A, B, C and D (see \refsec{secRootSystems}). The modified algorithm \refsec{secTheActualAlgorithm} has turned out in practice to work considerably faster than the straightforward application of \refsec{secPFdecomposition}. 

We have realized the algorithms presented in Sections \ref{secPFdecomposition}, \ref{secBVdecomposition} and \ref{secImplementation} under the name ``vector partition function'' program. The program and its source code are publicly available at \url{http://vectorpartition.sourceforge.net/} under the Library General Public License v3.0 (free use, modification, and redistribution as long as the same rights are granted to recipients in all consequent redistributions). \refsec{secProgramFunctionality} gives a very brief description of some of the features the program and the appendix presents printouts of (non-unique) partial fraction decompositions (as formatted by the program) of the generating function of the Kostant partition function for $A_2$, $A_3$, $A_4$, $B_2$, $B_3$, $C_2$, $C_3$, and $G_2$.

The current version of the ``vector partition function'' program (as of September 2009) can compute algebraic expressions and combinatorial chamber subdivisions for the Kostant partition function of $A_5$, $A_4$, $A_3$, $A_2$, $B_4$, $B_3$, $B_2$, $C_4$, $C_3$, $C_2$, $D_4$ and $G_2$. The software package LattE can also compute algebraic expressions for the vector partition function. LattE was released in 2003 by a team directed by J. A. De Loera at University of California, Davis, and is available at \url{http://www.math.ucdavis.edu/~latte/software.htm}. Another project that computes vector partition functions is S. Verdoolaege's program ``Barvinok'', \url{http://www.kotnet.org/~skimo/barvinok/}\nocite{SvenVerdoolaege:Barvinok}. For a detailed description of the algorithm of the ``Barvinok'' program see \cite{Kevin&others:BarvinokAlgorithm}\footnote{Up to the author's knowledge, a detailed explanation of the exact relation between the partial fractions approach and the Barvinok algorithm approach for computing vector partition functions is missing from the literature.}. Software for computing values of the Kostant partition function for the classical root systems has also been written by M. Vergne and C. Cochet in MAPLE.



\textbf{Acknowledgements.} The author would like to thank Michele Vergne, Thomas Bliem and Alexey Petukhov for the valuable advice and discussions on vector partition functions, as well as Pavel Tumarkin for the valuable discussion on root systems and Coxeter groups.

\subsection{Notation}\label{secNotation}
The notation is mostly standard with the exception of the vector space/dual space notation. An ambient vector space of dimension $n$ will be fixed and denoted as $\CC^n$. Greek letters will denote vectors, except for the number $\pi$ and the letter $\tau$ which will always stand for quasinumber. For a preferred choice of basis $\eps_1,\dots\eps_n$, we will denote by $\eps_i^*$ the corresponding elements of the dual vector space defined by $\eps_i^*(\eps_j)=\delta_{ij}$. Then the expression $x^{\gamma}$ refers to $x_1^{\eps_1^*(\gamma)}\dots x_n^{\eps_n^*(\gamma)}$. We will use the notation $\gamma^i:= \eps_i^*( \gamma)$ for the $i^{th}$ coordinate with respect to the preferred basis $\eps_i$. Whenever we use the preceding three notations, a preferred choice of basis will be either explicitly stated or implied. 

$\QQ[x_1,\dots,x_n,\frac{1}{x_1},\dots,\frac{1}{x_n}][\frac{1}{1-x^\alpha}]_{\alpha \in \ZZ^n}$ will denote the subring of $\QQ(x_1,$ $\dots,$ $x_n)$ generated over $\QQ$ by $x_i,x_i^{-1}$ and $\frac{1}{1-x^\alpha}$ for all $\alpha\in\ZZ^n$.

\subsubsection{Root systems}\label{secRootSystems}
The Kostant partition function is the vector partition function with respect to the vectors with non-negative coordinates of sets of integral vectors called root systems, (coordinates are chosen with respect to a special basis called ``simple''). Root systems play a central role in the representation theory of Lie algebras and classify the simple finite-dimensional Lie algebras. There are 4 classical series of root systems, $A_n, B_n,C_n,D_n$, (indexed by a parameter $n$) and five exceptional root systems, $E_6,E_7,E_8,F_4,G_2$. In the literature positive root systems of the root systems $A_n$, $B_n$, $C_n$ and $D_n$ are often realized in non-simple basis coordinates as
\begin{eqnarray*}
A_n, n\geq 2&:& \Delta^+=\{\eta_i-\eta_j|i< j\in\{1,\dots,n+1\}\};\\
B_n, n\geq 2&:& \Delta^+=\{\eta_i\pm\eta_j|i< j\in\{1,\dots,n\}\}\cup\{\eta_i|i\in\{1,\dots n\}\};\\
C_n, n\geq 2&:& \Delta^+=\{\eta_i\pm\eta_j|i\leq j\in\{1,\dots,n\}\}\backslash\{0\};\\
D_n, n\geq 4&:& \Delta^+=\{\eta_i\pm\eta_j|i< j\in\{1,\dots,n\}\}.
\end{eqnarray*}
and in simple basis coordinates as
\begin{eqnarray*}
A_n&:&\eps_1:=\eta_1-\eta_2,\dots,\eps_{n}:=\eta_n-\eta_{n+1}; \\
B_n&:& \eps_1:=\eta_1-\eta_2,\dots,\eps_{n-1}:=\eta_{n-1}-\eta_{n},\eps_{n}:=\eta_{n};\\
C_n&:& \eps_1:=\eta_1-\eta_2,\dots,\eps_{n-1}:=\eta_{n-1}-\eta_{n},\eps_{n}:=2\eta_{n};\\
D_n&:& \eps_1:=\eta_1-\eta_2,\dots,\eps_{n-1}:=\eta_{n-1}-\eta_{n},\eps_{n}:=\eta_{n-1}+\eta_{n}.
\end{eqnarray*}
List of the positive root systems of exceptional Lie algebras $E_6$, $E_7$, $E_8$, $F_4$ and $G_2$ can be found at \\ \url{http://www.liegroups.org/dissemination/spherical/explorer/rootSystem.cgi} (the Atlas of Lie Groups and Representations team) and are automatically generated by the ``vector partition function'' program.
\subsubsection{List of specific notation}
\begin{itemize}
\item $\QQ[x_1,\dots,x_n,\frac{1}{x_1},\dots,\frac{1}{x_n}][\frac{1}{1-x^\alpha}]_{\alpha \in \ZZ^n}$ see \refsec{secNotation}.
\item $\tau_{M,c,d}(\gamma)=\doublebrace{1}{\mathrm{if~}M\gamma-c=0 ~(\mod d)}{0}{\mathrm{otherwise}}$, see \refdef{defQNring}. $\tau_{M,c ,d }$ is called quasinumber, the $M$ is a $m\times n$ integer matrix, $c$ is an $m\times 1$ integer column vector.
\item $B_{k}(x)=\sum_{t=0}^x t^k, k\geq 1$, $B_{0}(x)=x+1$ - Bernoulli sum polynomial in $x$. We consider $B_k$ as an element of $\CC[x]$.
\item ${B}_{k}^{l,m,d}(x):=\sum_{t=0}^x \tau_{l,m,d}(t) t^k$ - $\tau$-Bernoulli sum, $m,l\in \ZZ_{\geq 0}$ (considered as 1 by 1 matrices), $d\in \ZZ_{ >0 }$, $x\in\ZZ$. 
\item ${\tilde B}_{k}^{l,m,d}(x):=\sum_{t=0}^x \tau_{l,m,d}(t) t^k$ - the extension of the $\tau$-Bernoulli sum over $x\in\QQ$ via the right hand side of \refeq{eqTauBernoulli}, where $m,l\in \ZZ_{\geq 0}$, $d\in \ZZ_{ >0 }$. 
\item $x^\alpha:= x_1^{\alpha^1}\dots x_n^{\alpha^n}$. See \refsec{secNotation}.
\item $\partial_j:=\frac{\partial}{\partial x_j}$.
\item $\Lambda_{R}(\{\alpha_i\}_{i=1}^N)= \{a_1\alpha_1+\dots+a_N\alpha_N| a_i\in R\}$ - cone (respectively lattice) spanned by the $\alpha_i$'s over $R$, where $R$ is an arbitrary subset (resp. a subring) of $\CC$.
\item $P_{\alpha_1,\dots,\alpha_n}:\CC^n\to\ZZ$ - vector partition function with respect to the vectors $\alpha_i$.
\item $\ZZ^n:= \Lambda_{\ZZ}(\{\eps_i\})$ - short notation for the lattice generated by the vectors $\eps_i$, where $\eps_i$ is the implied preferred choice of basis. 
\item $\langle\bullet,\bullet\rangle$ - bilinear form defined in \refsec{secBVdecomposition}.
\item $\langle\bullet,\bullet\rangle_K$ - symmetric bilinear positive definite form on a root system induced from the Killing form and normalized as in the proof of \refle{leCanChooseCoeffOne}.
\end{itemize}

\section{Partial fraction decomposition algorithm}\label{secPFdecomposition}
\subsection{Partial fraction decomposition definition}
\begin{defi}\label{defGeometricDefinitions}
A \emph{closed pointed polyhedral cone} $C$ is a subset of $\QQ^n$ of non-zero $n$-volume given by homogeneous non-strict rational linear inequalities. The hyperplanes given by the defining inequalities are called \emph{supporting hyperplanes} of $C$ on the condition that they intersect $C$ in a set of non-zero $n-1$-volume. The intersection of a supporting hyperplane with $C$ will be called a \emph{wall of $C$}. 
\end{defi}
\textbf{Remark.} A wall of a closed pointed polyhedral cone is in turn a closed pointed polyhedral cone in its corresponding supporting hyperplane.

\begin{defi}\label{defFractionIndexingSet}
A \emph{fraction indexing set} is defined to be a finite collection of triples of $s$-tuples $\{(\alpha_1,\dots,\alpha_s,$ $l_1,\dots, l_s,$ $m_1,\dots,m_s)\}$, where $s\in\ZZ_{>0}$, all $\alpha_i \in\ZZ^n,$  all $l_i\in\ZZ_{>0},$ all $m_i\in\ZZ_{>0} $. 
\end{defi}
The following definition is similar to \cite[Proposition 3.2]{DeConciniProcesi:GeometryToricArrangements}. 
\begin{defi}\label{defPartFracDecomposition}
Let $q\in\QQ[x_1,\dots,x_n,\frac{1}{x_1},\dots,\frac{1}{x_n}][\frac{1}{1-x^\alpha}]_{\alpha \in \ZZ^n}$. A fraction indexing set for which all participating $s$-tuples $\alpha_1,\dots,\alpha_s$ of vectors are \textbf{linearly independent} is called a \emph{partial fraction indexing set}. \emph{A partial fraction decomposition of $q$} indexed by a partial fraction indexing set $I$ is defined to be a finite set of polynomials $\mathcal{P}_I = \{p_{ \alpha_{1}, \dots, \alpha_{s},l_1,\dots l_s,m_1, \dots,m_s}$ $~|~(\alpha_1,\dots,\alpha_s,l_1,\dots, l_s,m_1,\dots,m_s)\in I\} \subset\QQ[x_1 ,\dots, x_n, \frac{1} {x_1}, \dots, \frac{1} {x_n}]$  such that
\begin{equation}\label{eqPartFracDecomposition}
q= \sum_{\substack{(\alpha_1,\dots,\alpha_s,\\l_1,\dots, l_s,\\m_1,\dots,m_s)\in I}} p_{\substack{(\alpha_1,\dots,\alpha_s,\\l_1,\dots, l_s,\\m_1,\dots,m_s)}}(x) \frac{1} {(1-x^{l_1\alpha_{1}})^{m_1}}\dots \frac{1}{(1-x^{l_s\alpha_{s}})^{m_s}}\quad .
\end{equation}
The set of pointed polyhedral cones $\{\Lambda_{\QQ_{\geq 0}} ( \{ \alpha_{1}, \dots, \alpha_{s}\})$ $~|~ (\alpha_1, \dots, \alpha_s,$ $l_1, \dots, l_s,$ $m_1, \dots, m_s)$ $\in I\mathrm{~for~some~}m_i,l_i\}$ will be called \emph{cone support} of the partial fraction decomposition $\mathcal{P}_I$. The finite set of all vectors participating as generators of the elements of the cone support of $\mathcal{P}_I$ will be called \emph{support of $\mathcal{P}_I$}.
\end{defi}
\textbf{Remark.} We note that the ``elongations'' $l_i$ in the above definition could be substituted for a single ``greatest elongation'' by the geometric series sum formula (``the elongation formula''). However, this would surely yield great computer memory loss in realizations of the algorithm in \refsec{secTheTheoreticalAlgorithm} below. 
\subsection{Formulas used}
\begin{lemma}~\label{leFlas} 
\begin{itemize}
\item[(a)](Szenes-Vergne formula) Let $\alpha_1,\dots \alpha_k\in\ZZ^n$ such that $\sum_{i}\alpha_i\neq 0$ Then
\begin{eqnarray*}
&&\frac{1}{1-x^{\alpha_1}}\dots \frac{1}{1-x^{\alpha_k}} =\\ &&\frac{1}{1-x^{\sum\alpha_i}}\left(\sum_{j=1}^{j=k}\frac{x^{\alpha_1}\dots x^{\alpha_{j-1}}}{(1-x^{\alpha_1})\dots (1-x^{\alpha_{j-1}})}\frac{1}{(1-x^{\alpha_{j+1}})\dots(1-x^{\alpha_{k}})} \right).
\end{eqnarray*} 
\item[(b)] (``elongation formula'' (geometric series sum formula))
\begin{eqnarray*}
\frac{1}{1-x^\alpha}&=&\frac{1+x^\alpha+\dots+x^{(n-1)\alpha}}{1-x^{n\alpha}},
\end{eqnarray*} 
for $n\in\ZZ_{>0}$, and  
\begin{eqnarray*}
\frac{1}{1-x^\alpha}&=&\frac{-x^{-\alpha}-\dots-x^{n\alpha}}{1-x^{n\alpha}},
\end{eqnarray*} 
for $n\in\ZZ_{<0}$.
\item[(c)](the elongated Szenes-Vergne formula) Let $k\geq 1$ and $\alpha_1,\dots \alpha_k\in\ZZ^n, a_i\in\ZZ_{\neq 0}$ such that $\sum_{i}a_i\alpha_i\neq 0$ Then
\begin{eqnarray*}
&&\frac{1}{1-x^{\alpha_1}}\dots \frac{1}{1-x^{\alpha_k}}=\\ 
&&\frac{1}{1-x^{\sum a_i\alpha_i}}\left(\sum_{j=1}^{j=k}\frac{x^{a_1\alpha_1}\dots x^{a_{j-1}\alpha_{j-1}}}{(1-x^{\alpha_1})\dots (1-x^{\alpha_{j-1}})} \frac{p_{a_j}(x^{\alpha_j})} {(1-x^{\alpha_{j+1}}) \dots(1- x^{\alpha_{k}})} \right),
\end{eqnarray*} 
where $p_n(x^\alpha)$ are the numerators in the elongation formula (b). \textbf{Remark.} When $k=1$ (c) degenerates to (b).
\item[(d)] 
\begin{eqnarray*}
\frac{1}{1-x^\alpha} \frac{1}{1-x^\beta}&=& \frac{1}{1-x^{\alpha-\beta}}\left(\frac{1}{1-x^{\beta}}-\frac{x^{\alpha-\beta}}{1-x^{\alpha}} \right).
\end{eqnarray*} 
\begin{eqnarray*}
\frac{1}{1-x^{\alpha}} \frac{1}{1-x^\beta}&=& \frac{1}{1-x^{\alpha-2\beta}}\left(\frac{1}{1-x^{\beta}}-\frac{x^{\alpha-\beta}+x^{\alpha-2\beta}}{1-x^{\alpha}} \right).
\end{eqnarray*} 
\item[(e)] Let $p:=-x^{\alpha}(x^{-\beta}+\dots+x^{-n\beta})$. Then 
\begin{eqnarray*}
&&\frac{1}{(1-x^\alpha)^l} \frac{1}{(1-x^\beta)^m}=\\&& ~\sum_{t=1}^{l} \frac{1}{(1-x^\alpha)^t} \frac{p^{m}}{(1-x^{\alpha-n\beta})^{l+m-t}} \binomial{l+m-t-1}{m-1}\\&& +\sum_{t=1}^{m}\frac{1}{(1-x^\beta)^t} \frac{p^{m-t}}{(1-x^{\alpha-n\beta})^{l+m-t}} \binomial{l +m-t-1}{l-1}.
\end{eqnarray*}
\end{itemize}
\end{lemma}
\begin{proof}
(a) is a partial case of \cite[Lemma 1.8]{SzenesVergne:ResidueFlas}. (c) is a combination of (a) and (b). Formulas (d) and (e) are used in \refsec{secImplementation} and have been implemented together with formula (c) in the ``vector partition function'' program. The coefficients in (e) come from counting the number of ways one can get from point $(t,1)$ to point $(l,m)$ in a square grid by moving up and right with unit steps.
\end{proof}

\textbf{Remark.} One could realize a formula similar to (e) for more than two variables.

\subsection{An algorithm for computing partial fraction decompositions} \label{secTheAlg} \label{secTheTheoreticalAlgorithm}
We are given a finite set of integral non-zero vectors with non-negative coordinates $\bar\Delta= \{\alpha_1, \dots, \alpha_N\}\subset\ZZ^n$ and a finite set $\bar{\bar\Delta}$ of multiples of vectors in $\bar\Delta$.  
Let an element $q\in\QL$ be given in the form 
\begin{equation}\label{eqTheoreticalAlg}
q=\sum_{\substack{l_j\alpha_{i_j}\in\bar{\bar\Delta}, \forall j \\ m_i\in\ZZ_{>0}}} p_{\alpha_{i_1},\dots,\alpha_{i_k},m_1,\dots,m_k}(x)\frac{1}{(1-x^{\alpha_{i_1}})^{m_1}} \dots\frac{ 1}{( 1- x^{ \alpha_{i_k} })^{m_k}},
\end{equation}
where $p_{\alpha_{i_1},\dots,\alpha_{i_k},m_1,\dots,m_k}(x)\in\QQLaurent$. The algorithm presented here gives us a partial fraction decomposition $\mathcal{P}_I$ of $q$ with support $\bar\Delta$ (see \refdef{defPartFracDecomposition}). The denominators of the so produced partial fraction decomposition will involve in general multiples of vectors in $\Delta$ lying outside of $\bar{\bar\Delta}$.

The element $q$ will be stored in the computer's memory as a list of elements called \emph{fractions}, corresponding to each summand in \refeq{eqTheoreticalAlg}. For a fixed fraction the number $m_l$ will be called the multiplicity of the vector $\alpha_{i_l}$, and $\{\alpha_{i_1}, \dots, \alpha_{i_k}\}$ will be called the \emph{support} of the fraction. The fractions' supports are partially ordered by inclusion. Besides the data implied by \refeq{eqTheoreticalAlg}, for each fraction we store a boolean flag indicating whether our fraction is reduced and a boolean flag indicating whether we have a ``preferred linear dependence''. The ``reduced'' flag is set only if the vectors in the fraction's support are linearly independent. If the ``preferred linear dependence'' flag is set, we store an additional data structure representing a non-zero $\ZZ$-linear dependence $a_1\alpha_{i_1} +\dots + a_k\alpha_{i_k} =a_0\alpha_{i_0}$ with a distinguished index $i_0$ for which $a_0> 0$. The vector $\alpha_{i_0}$ we define to be the \emph{gaining multiplicity vector}.

\begin{itemize}
\item[Step 0] Initialize the starting data by setting all fractions to ``non-reduced'' and not having a ``preferred linear dependence''.
\item[Step 1] Select an arbitrary fraction from \refeq{eqTheoreticalAlg} that has a ``preferred linear dependence''. In case no fraction with a ``preferred linear dependence'' exists, select an arbitrary fraction from \refeq{eqTheoreticalAlg} that is not labeled as ``reduced''. Choose, if possible, in an arbitrary fashion a linear dependence $a_0\alpha_{i_0}= a_1\alpha_{i_1} +\dots +a_k \alpha_{ i_k}$, $a_0> 0$, $a_{i_j}\neq 0, \forall j$, and store it as a ``preferred linear dependence''. 

If no choice of ``preferred linear dependence'' is possible (i.e. the support of the fraction consists of linearly independent vectors) label the fraction as ``reduced'' and go to the beginning of Step 1. 

If all fractions from \refeq{eqTheoreticalAlg} are labeled as ``reduced'' terminate the program. 
\item[Step 2] Apply the ``elongated Szenes-Vergne formula'' (\refle{leFlas})(c) and the ``elongation formula''(\refle{leFlas}(b)) as indicated below 
\begin{eqnarray*}
&&\underbrace{\left(\prod_{j=1}^k \frac{1}{1-x^{\alpha_{{i_j}}}} \right)}_{ \mathrm{ apply ~\refle{leFlas} (c)~here}}\left(\prod_{j=1}^k \frac{1}{(1-x^{\alpha_{{i_j}}})^{m_j-1}} \right)\\&& \left(\underbrace{\frac{1}{(1-x^{\alpha_{i_0}})}}_{ \mathrm{ apply ~\refle{leFlas} (b)~with~elongation~}a_0\mathrm{~here}}\right)^{m_0} (\mathrm{ ~the~rest~of~the~multiples }).
\end{eqnarray*}
 This procedure will erase the selected fraction from the computer's memory and substitute it with new list of fractions. All newly produced fractions that have the same support as the originally selected fraction are set to inherit  the ``preferred linear dependence''. The newly produced fractions with a strictly smaller support are labeled as ``not having a preferred linear dependence''.

\textbf{Remark} In the ``vector partition function'' program the application of \refle{leFlas} (b) is made after all applications of the Szenes-Vergne formula are completed. This increases computational speed significantly, but requires additional accounting.
\item[Step 3] Add the vector $a_0\alpha_{i_0}$ to the set $\bar{\bar\Delta}$. Go to Step 1 with the so modified list of fractions $q$ and set $\bar{\bar\Delta}$.
\end{itemize}
\begin{lemma}
The above algorithm will come to a halt by labeling all fractions as reduced. The so produced sum will represent a partial fraction decomposition for the starting element $q$, whose support is a subset of $\bar\Delta$.
\end{lemma}
\begin{proof}
If a fraction participating in \refeq{eqTheoreticalAlg} has a ``preferred linear dependence'', Step 5 increases the multiplicity of the ``gaining multiplicity vector'' while decreasing the multiplicity of one of the other vectors in the support of the fraction. Thus, after finitely many steps, the support of all newly produced fractions will strictly decrease. Therefore all fractions in the final expression will end up being reduced.
\end{proof}

\textbf{Remark (by an idea of Thomas Bliem). } In case there is no ``preferred linear dependence'' chosen, there is freedom in choosing the index $a_0$\footnote{and the order of the vectors in the preferred linear dependence} in Step 1 and one can choose any strategy. One such choice leads to partial fraction decompositions that involve non-broken circuit bases in the sense of \cite{ConciniProcesi}. Suppose one wanted to compute a partial fraction decomposition of $\frac{1}{1-x^{\alpha_{i_1}}}, \dots,\frac{1} {1-x^{ \alpha_{ i_N}}}$. Suppose one knew by induction that a subset $S:=\alpha_{i_1},\dots, \alpha_{i_k}$ of the vectors $T:=\alpha_{i_1},\dots, \alpha_{i_N}$, $i_1< \dots < i_N$, was linearly independent. Then one would add vector $\alpha_{i_{k+1}}$ to $S$ and check whether the span of $ S\cup \{\alpha_{i_{k+1}}\}$ is maximal; if not, one would label vector $\alpha_{i_{k+1}}$ in Step 1 of our algorithm as ``gaining multiplicity'' and proceed to Step 5. If the set $ S\cup \{\alpha_{i_{k+1}}\}$ were of full rank, then one would substitute $S$ by it and continue with the next index. This ``simple strategy'' of choosing indices in fact yields a partial fraction decomposition with respect to non-broken circuit bases (\cite{ConciniProcesi}). By definition, an $n$-tuple of vectors $\alpha_{i_1}, \dots , \alpha_{ i_n}$, $i_1<\dots<i_n$ is a non-broken circuit base if $\alpha_{i_1}, \dots, \alpha_{i_k}$ and $\alpha_m$ are not linearly dependent for any index $m$ with $i_k\leq m\leq i_{k+1}$. We leave that to the reader to prove that our ``simple strategy'' of choosing indices in Step 1 will produce a partial fraction decomposition involving denominators whose corresponding vectors form non-broken circuit bases. It is important to note however, that the results in \refsec{secPeriods} are based on a different strategy on choosing the ``gaining multiplicity'' vector (namely, the vector with a smallest possible in absolute value coefficient in a non-trivial integral linear relation). This strategy has been realized in the latest version of the ``vector partition function'' program. 

For the generating function of the Kostant partition function of the classical root systems of type A, B, C or D (see \refsec{secRootSystems}) one can considerably simplify the algorithm, see \refsec{secTheActualAlgorithm}. 

\section{Brion-Vergne decomposition for partial fraction decompositions} \label{secBVdecomposition}

For any vector $\gamma$ define the differential operator $\xi_{\gamma}:=\sum \gamma^kx_k\partial_k$ in the ring $\QL$. 
Let $\alpha_1,\dots,\alpha_n$ be linearly independent vectors where $n$ is the dimension of the ambient vector space. Define the scalar product $\langle\bullet,\bullet\rangle$ by $\langle\alpha_i, \alpha_j \rangle= \sum \alpha_i^k\alpha_j^k$. Let $\beta_i$ be defined as the vector whose scalar product with $\alpha_k$ equals $\delta_{ik}$. We have that $\xi_{\beta_i} (x^{ \alpha_k })= \doublebrace{ x^{\alpha_i}} {\mathrm{if~} i=k}{0}{ \mathrm{ otherwise} }$. Then 
\begin{eqnarray*}
\left(\frac{\xi_{\beta_i}}{m_i}+\id\right)\left(\frac{1}{(1-x^{\alpha_i})^{m_i}}\right)= \frac{1}{(1-x^{\alpha_i})^{m_i+1}}
\end{eqnarray*}
and so we obtain a Brion-Vergne decomposition, \cite{BrionVergne}: 
\begin{eqnarray}\label{eqBVdecomposition}\notag
&&\left( \prod_{i=1}^n\left(\left(\frac{\xi_{\beta_i}}{1}+\id\right)\dots\left(\frac{\xi_{\beta_i}}{m_i-1}+\id\right) \right)\right) \left(\frac{1}{1-x^{\alpha_1}}\dots\frac{1}{1-x^{\alpha_n}}\right)= \\&& \frac{1}{(1-x^{\alpha_1})^{m_1}} \dots\frac{1}{(1-x^{\alpha_n})^{m_n}}\quad.
\end{eqnarray}
Let $\Lambda:=\Lambda_{\ZZ}(\alpha_1,\dots,\alpha_n)$. Now let us explore what happens when we apply $\xi_{\beta}$ to a power series $A(x):=\sum_{\gamma\in\Lambda} P(\gamma)x^{\gamma}$: 
\[
\xi_{\beta}\left(\sum_{\gamma\in\Lambda} P(\gamma)x^{\gamma}\right)= \sum_{\gamma\in\Lambda} P(\gamma)\langle\beta,\gamma\rangle x^{\gamma}\quad.
\]
Now let us explore what happens when we apply multiplication by $x^{\beta}$:
\[
x^{\beta}\sum_{\gamma\in\Lambda} P(\gamma)x^\gamma= \sum_{\gamma\in\Lambda} P(\gamma-\beta)x^{\gamma}.
\] 
Let $\sum_{\gamma\in\Lambda}\tau(\gamma)x^{\gamma}$ be the power series expansion of $\frac{1}{(1-x^{\alpha_1})\dots(1-x^{\alpha_n})}$ given by uncovering the brackets in the formal product $(1+x^{\alpha_1}+x^{2\alpha_1}+\dots)\dots(1+x^{\alpha_n}+x^{2\alpha_n}+\dots)$. Then 
\begin{equation}\label{eqQuasiNumberNE}
\tau(\gamma) = \doublebrace{1}{\mathrm{if~}\gamma\mathrm{~lies~in~}\Lambda_{\ZZ\geq 0}(\{\alpha_i\})}{0}{\mathrm{otherwise}}.
\end{equation}
Summarizing the above formulas, we get
\begin{eqnarray} \notag 
\sum_{\gamma\in\Lambda}P(\gamma)x^\gamma&:=&x^{\delta}\prod_{i=1}^n\frac{1}{(1-x^{\alpha_i})^{m_i}}\\&=&\sum_{\gamma\in\Lambda} \prod_{i=1}^n\prod_{l=1}^{m_i-1}\left(\frac{\langle\beta_i,\gamma-\delta\rangle}{l}+1 \right)\tau(\gamma-\delta)x^{\gamma} \notag \\
\label{eqtheVPfunction} P(\gamma) &=& \prod_{i=1}^n\prod_{l=1}^{m_i-1} \left(\frac{\langle\beta_i, \gamma-\delta \rangle} {l} + 1 \right)\tau(\gamma-\delta)  \notag \\
&=&\prod_{i=1}^n\binomial{\langle\beta_i,\gamma-\delta\rangle+m_i-1}{m_i-1}\tau(\gamma-\delta).
\end{eqnarray}
\subsection{Realizing \refeq{eqtheVPfunction} with quasipolynomials}\label{secRealization}

Formula \refeq{eqtheVPfunction} clearly raises the question of how to realize the function $\tau(\gamma-\delta)$. For the linear inequalities implied by ``$\gamma-\delta$ $\mathrm{~ lies ~ in ~}$ $\Lambda_{ \ZZ \geq 0}( \alpha_1,\dots,\alpha_n )$'', we can deal as in \cite{BBCV}. First, we fix an ``indicator'' vector in general position. Then we take into consideration only the fractions whose support generates cones that contain our ``indicator'' vector. In order for this approach to work we would need to start with an ``indicator'' vector    ``far away we are from the walls'', since there is a $-\delta$ shift in the argument of $\tau$. This naturally raises the question, how does the quasipolynomial expression for the vector partition function change when its argument is close to a wall of a cone from the cone support of a partial fraction decomposition.   

The answer to the this question is well known (see \cite{ConciniProcesiVergne}, \cite{Sturmfels}, \cite{SzenesVergne:ResidueFlas}) - the quasipolynomial expression for the vector partition function does not depend on the distance from the walls. We give our own proof (\refth{thMain}) in \refsec{secWallDist}. Of course, this immediately implies that the algebraic expressions for the vector partition function, which are different for different pointed polyhedral cones, take the same values over a common wall of neighboring pointed polyhedral cones. 

In order to represent in the computer's memory the function $\tau$, we store the linear inequalities in a separate structure and introduce the helping function $\bar \tau$ which ``forgets'' the linear inequalities:  
\[
\bar \tau(\gamma):= \doublebrace{1}{\mathrm{if~}\gamma\mathrm{~lies~in~}\Lambda_{\ZZ}(\alpha_1,\dots,\alpha_n)}{0}{\mathrm{otherwise}}.
\]
Let $B$ be the matrix 
$B:=\left(\begin{array}{ccc} b_{11}&\dots& b_{1n}\\ \vdots&&\vdots\\ b_{n1} &\dots &b_{nn} \end{array} \right) :=\left(\begin{array}{ccc} \alpha_1^1&\dots& \alpha_n^1\\ \vdots&&\vdots\\ \alpha_1^n &\dots &\alpha_n^n \end{array} \right)^{-1}$. Then $\bar \tau(\gamma-\delta)$ is non-zero if and only if 
\begin{equation}\label{eqDens1}
B\left(\begin{array}{c}\gamma^1-\delta^1\\ \vdots\\ \gamma^n-\delta^n\end{array}\right)\in\ZZ^{n}.
\end{equation} 
Set $c_i:=\sum b_{in}\delta^n$ and define $d_i$ to be the smallest positive integer for which $d_ic_i\in\ZZ,d_ib_{ik}\in\ZZ,~\forall{k}$ (in other words, $d_i$ is the least common multiple of the denominators involved in the $i^{th}$-line of \refeq{eqDens1}). Define $d$ to be the least common multiple of the $d_i$-s ($d$ will later be called ``period''). 

\noindent\textbf{Remarks.} $d$ is the least common multiple of the orders of the elements of the group $\ZZ^n/\Lambda_{\ZZ} (\{\alpha_i\})$. Thus $d$ is a (possibly proper) divisor of the index of the subgroup $\Lambda_{\ZZ} (\{\alpha_i\})$ in the group  $\ZZ^n$ (which in turn equals $\det(\alpha_i^j)$). The period $d$ does not change when we replace $\alpha_i, \alpha_j$ by $\pm\alpha_i, \pm\alpha_i\pm\alpha_j$ for any two indices $i$ and $j$ and any choice of signs (since the lattice doesn't change).

$\bar \tau(\gamma-\delta)$ is non-zero if and only if
\begin{equation}\label{eqTheSystem}
d B
\left(\begin{array}{c}
\gamma^1 \\
\vdots\\
\gamma^n
\end{array}\right)= 
d\left(\begin{array}{c}
c_1\\
\vdots\\
c_n
\end{array}\right) (\mod d).
\end{equation}
Let $M= dB$ be integer matrix of the above system, $c$ be the integer vector column of the system. Define $\tau_{\Lambda(\alpha_1,\dots,\alpha_n),\delta}(\gamma^{1},\dots,\gamma^{n}):=\tau_{M,c,d}(\gamma^{1},\dots,\gamma^{n}):=$ \\ $\doublebrace{1}{\mathrm{if~} \gamma_i \mathrm{ ~are~ integers}\\&\mathrm{~and ~the~ system~ \refeq{eqTheSystem}}\\&\mathrm{ ~has ~integer~ solution}} { 0} { \mathrm{otherwise}}$, where $\tau_{M, c,d} :\CC^{n}\to\QQ$. The notation $\tau_{\Lambda(\alpha_1,\dots,\alpha_n),\delta}$ is more suggestive mathematically, but for the computer realization it was more convenient to use the $\tau_{M,c,d}$-notation, and that is why we choose the latter for the rest of this the paper. We extend the definition of $\tau_{M,c,d}$ to $M$ being a $m\times n$ integer matrix and $c$ a $m\times 1$ vector column. This makes sense because by definition the function $\tau_{M,c,d}$ is supported on the a priori defined lattice $\ZZ^n$. 
\begin{defi}\label{defQNring}
The subring of the ring of functions $f:\CC^{n}\to\QQ$ generated by $\tau_{M,c,d}$ with coefficients in $\QQ$ we will call the ring of \emph{quasinumbers}. The elements $\tau_{M,c,d}$ will be called \emph{basic quasinumbers}.

The subring of the functions $f:\CC^n\to\QQ$ generated by the quasinumbers and $\eps_i^*$ we will call the ring of quasipolynomials.
\end{defi}

We note that relating expressions written with complex exponents and quasinumbers requires extra caution. One precise way of doing it is by saying that to every quasinumber $q$ over the lattice $\ZZ^n$ and for any lattice refinement $\frac{1}{d}\ZZ^n$ one can assign non-uniquely a holomorphic function $f$ of the form $\displaystyle f=\sum_{\substack{\{ a_i\in\QQ\} \\ p_{\{a_i\}}\in\QQ[x_1,\dots,x_n]}} p_{\{a_i\}} (\eps_1^*,\dots,\eps_n^*) e^{2\pi i (\sum a_i \eps_i^*)}$, where the $p_{\{a_i\}}$'s are polynomials with rational coefficients, such that the functions $f$ and $q$ have the same restrictions on the lattice $\frac{1}{d}\ZZ^n\subset\CC^n$. One can do that by combining the one variable formula in the appendix with integral linear substitutions. The non-uniqueness is clear: the zero function and the function $e^{2\pi i\eps_i^*}-1$ both correspond to the zero quasinumber for the lattice $\ZZ^n$, but not for the lattice $\frac{1}{2}\ZZ^n$.

Making Gaussian elimination transformations $\mod d$ to the system \refeq{eqTheSystem} does not change $\tau_{M,c,d}$. Such transformations have been realized in our ``vector partition function'' program and are used to determine when $\tau_{M,c,d}$ is zero and to multiply elements $\tau_{M,c,d}$ by  $\tau_{M',c',d'}$  to yield an expression of the same type $\tau_{M'',c'',d''}$. However, the program does not realize all additive relations in the ring of quasi-numbers. We finish this subsection with a definition and an important lemma.
\begin{defi}\label{defPeriod}
A least common period for a finite set $\{\tau_{M_i,c_i,d_i}\}$ we call the least common multiple of the $d_i$-s. Let the quasipolynomial $q$ be written via the generators $\eps_i^*$ and the elements of the finite set $\{ \tau_{ M_i, c_i, d_i}\}$. Then the least common period of the set $\{ \tau_{ M_i, c_i, d_i}\}$ is called a possible common period for $q$. The minimum over all possible common periods of $q$ is called the least common period for $q$ or simply the period of $q$.  
\end{defi}
\begin{lemma}\label{leTheQuasipolynomialLemma}
We are using the notation from \refdef{defGeometricDefinitions}. If two quasipolynomials $q_1$ and $q_2$ coincide on an affine translate $C'$ of a closed pointed polyhedral cone $C$ of non-zero measure, they coincide everywhere. 
\end{lemma}
\begin{proof}
Let $\Lambda$ be the intersection of the lattices naturally determined by the basic quasinumbers participating in the expressions for $q_1$ and $q_2$. Then both $q_1$ and $q_2$ are polynomials when restricted over each lattice translate $\delta+\Lambda$ as $\delta$ varies over the representatives of elements of the group $\ZZ^n/\Lambda$. Now the statement of the lemma follows from the corresponding statement for polynomials.   
\end{proof}

\subsection{An algorithm for computing the vector partition function with respect to a set $I$ of non-zero integral vectors with nonnegative entries}\label{secTheVPalgNonRealized}

\begin{itemize}
\item Step 1. Apply the algorithm from \refsec{secTheTheoreticalAlgorithm} to get a partial fraction decomposition of the generating function of the vector partition function, $\prod_{\alpha\in I}\frac{1}{1-x^\alpha}$.
\item Step 2. Fix a combinatorial chamber of interest and apply formula \refeq{eqtheVPfunction} to the partial fraction decomposition obtained in Step 1. One needs to fix an vector in general position which indicates the combinatorial chamber. 
\end{itemize}

\section{Partition function algebraic expressions do not depend on the distance from the walls}\label{secWallDist}

This section presents an independent, elementary, but computationally inefficient algorithm for calculating the vector partition function, whose aim is to prove that the vector partition function does not depend on the distance from the walls of the combinatorial chambers. \footnote{The author has not fully understood yet the interplay between the current algorithm and the preceding one.}  

Suppose we want to compute the vector partition function $P_{\alpha_1, \dots, \alpha_m}$ with respect to the vectors $\alpha_1,\dots, \alpha_{m-1}, \alpha_m$ and have already computed the partition function with respect to the first $m-1$ vectors (call it $P_{\alpha_1,\dots,\alpha_{m-1}}$). Our starting point is the following completely elementary observation: if one chooses exactly $t$ times the vector $\alpha_n$, one can decompose a vector $\gamma$ in exactly $P_{\alpha_1,\dots,\alpha_{m-1}}(\gamma-t\alpha_m)(\gamma)$ ways. Therefore 
\begin{equation}\label{eqAlg2Step1}
P_{\alpha_1,\dots,\alpha_m}(\gamma)= \sum_{t=0}^{\infty}P_{\alpha_1,\dots,\alpha_{m-1}}(\gamma-t\alpha_m) .
\end{equation}
We will gradually refine \refeq{eqAlg2Step1} to a complete algorithm how to compute the vector partition function. This elementary idea is the essence of this section\footnote{The author hopes that the resulting rather cumbersome formulas will not obscure it.}.   

Suppose as an induction hypothesis there exist finitely many closed pointed polyhedral cones $C_0,\dots,C_M$ such that $P_{\alpha_1, \dots, \alpha_{m-1}}$, restricted to each $C_i$, is a quasipolynomial; extend this quasipolynomial naturally to the entire space and call the corresponding element $P_{C_i}$. Suppose now we want to compute $P_{\alpha_1, \dots, \alpha_m}$ at a point $\gamma$ in general position, i.e. such that $\gamma$ lies in no support hyperplane of any of the $C_i$. Let $H$ be a wall of an arbitrary $C_j$ such that the ray $R:=\{ \gamma + t \alpha_m |t\in[0,\infty)\}$ intersects $H$. By the requirement that $\gamma$ be in general position, we get that $R$ intersects $H$ transversally, i.e. if $\nu_H$ is a normal vector for $H$ then $\langle\alpha_m,\nu_H\rangle\neq 0$. Assume in addition $\nu_H$ is chosen to point towards the interior of $C_j$ (i.e. $\langle \nu_H, \delta\rangle\geq 0,\forall \delta\in C_j$). 
\begin{defi}\label{defExitWall}
If $\langle \nu_H,\alpha_m,\rangle<0$ (respectively $\langle \nu_H,\alpha_m,\rangle>0$) we call such a wall $H$ \emph{an exit (}respectively \emph{entering) wall for the point $\gamma$ with respect to the direction $\alpha_m$},
We will use the shortened term ``exit (respectively entering) wall'' whenever the rest of the data is implied by the context. 
\end{defi}
$\nu_H$ must have rational coordinates, so assume in addition that it is rescaled to have integer coordinates. Then the point through which $R$ exits $C_j$ is $\gamma- \frac{\langle \gamma, \nu_{H} \rangle} {\langle \nu_{H}, \alpha \rangle}\alpha$; in coordinates this is the point $(\gamma^1- \alpha^1( \sum \gamma^i \nu_H^i),\dots ,\gamma^n -\alpha^m ( \sum\gamma^i\nu_H^i))$. Note that the coordinates of this exit point are a linear rational expression of the coordinates $\gamma^i$. Let $H_0,H_1,\dots,H_N$ be exit walls and $H_1,\dots,H_{N-1}$, be entering walls for the point $\gamma$ with respect to direction $\alpha_n$ and $C_0,\dots,C_N$, $N\leq M$ be pointed polyhedral cones such that $H_i$, $N\geq i\geq 1$ is a shared wall of $C_{i-1}$ and $C_{i}$ and $H_N$ is a wall of $C_N$. Assume in addition that the exit wall $H_N$ is a subset of a wall of $\Lambda_{\QQ}(\alpha_1,\dots,\alpha_N)$. Let the ray $R$ intersect $H_l$ in point $A_l$. Then $A_l=\gamma-\frac{ \langle \gamma , \nu_{H_l} \rangle} {\langle \alpha_m, \nu_{H_l}\rangle}\alpha_m$. Assume finally  that $C_0\cup\dots\cup C_N$ are all pointed polyhedral cones that have a common point with the segment $\{ \gamma + t \alpha_m |t\in[0,\gamma-\frac{ \langle \gamma , \nu_{H_l} \rangle} {\langle \alpha_m, \nu_{H_l}\rangle}\alpha_m]\}$. Set 
\begin{equation}\label{eqIntemediatePoints}
t_l:=\frac{ \langle \gamma , \nu_{H_l} \rangle} {\langle \alpha_m, \nu_{H_l}\rangle} , 1\leq l\leq N-1;\quad t_0:=0.
\end{equation}  The contribution of the first chamber $C_0$ to the sum \refeq{eqAlg2Step1} is 
\begin{equation}\label{eqContributionC0}
\sum_{t=0}^{\lfloor t_1\rfloor} P_{C_0}(\gamma -t \alpha_m ) = \sum_{ t=0}^{\lfloor t_{N-1}\rfloor}P_{C_0}(\gamma-t\alpha_m)-\sum_{ t=\lfloor t_1\rfloor+1}^{\lfloor t_{N-1}\rfloor} P_{C_0} (\gamma - t\alpha_m).
\end{equation}
The contribution any other chamber $C_j$, $N>j> 0$ is
\begin{equation}\label{eqContributionCj}
\sum_{t=\lfloor t_j \rfloor +1}^{\lfloor t_{j+1}\rfloor} P_{C_j}(\gamma -t \alpha_m ) = \sum_{ t=\lfloor t_j\rfloor+1 }^{\lfloor t_{N-1}\rfloor} P_{C_j} (\gamma- t \alpha_m)-\sum_{ t=\lfloor t_{j+1} \rfloor+1}^{\lfloor t_{N-1}\rfloor} P_{C_j} (\gamma - t\alpha_m).
\end{equation}

Note that if $\lfloor t_{j+1} \rfloor<\lfloor t_{j} \rfloor+1$
the latter expression is zero. Summing over all chambers $C_i$ we get the first refinement of \refeq{eqAlg2Step1}:
\begin{eqnarray}\label{eqAlg2Step2}\notag
P_{\alpha_1,\dots,\alpha_m}(\gamma)&=& P_{C_0}(\gamma)+ \sum_{j=0}^{N-2} \left(\sum_{t=\lfloor t_j\rfloor+1}^{t=\lfloor t_{N-1}\rfloor} P_{C_j}(\gamma-t\alpha_m)\right)\\&& - \sum_{j=1}^{N-1}\left(\sum_{t=\lfloor t_{j+1}\rfloor+1}^{t=\lfloor t_{N-1}\rfloor} P_{C_j}(\gamma-t\alpha_m)\right) .
\end{eqnarray}
Note that if the lattice spanned by the $\alpha_i$'s is unimodular, (for definition see \cite{Sturmfels}), \refeq{eqAlg2Step2} is immediately computable via Bernoulli polynomials. To compute \refeq{eqAlg2Step2} in closed form one needs to know in advance which are all possible exit/entering walls for the point $\gamma$ with respect to the direction $\alpha_m$ which in turn will account for new subdivisions of the chambers $C_i$. 

\begin{lemma}\label{leMainLemmaAlg2}
Suppose $C_0$ is a closed pointed polyhedral cone with walls with rational normal vectors such that there is a unique choice of exit walls $H_0,\dots,H_N$ for a point $\gamma$ in $C_0$ with respect to the direction $\alpha_m$. Let $C_0,\dots,C_N$ and $P_{C_0},\dots,P_{C_N}$ be as above. Suppose in addition each $P_{C_i}$ is a quasipolynomial. Then  $P_{\alpha_1,\dots,\alpha_{m}}$ is a quasipolynomial over $C_0$.
\end{lemma}  
\begin{proof}
(i) Assume $\gamma$ is in general position with respect to the supporting walls of $C_0,\dots,C_N$. Due to the uniqueness of the choice of exit walls such points $\gamma$ lie in the interior of $C_0$.

To prove the lemma holds for (i) all we need to do is show how to express \refeq{eqAlg2Step2} with quasinumbers. Assume the $t_i$'s are written as linear expressions in the $\gamma^i$'s with rational coefficients represented by irreducible fractions and let $d$ be the least common multiple of the denominators of the $t_i$'s. Assume that $d>1$.


Let $f:\frac{1}{d}\ZZ\to\CC$ be an arbitrary function. Then
\begin{equation}\label{eqTauInterpolation}
f\left(\lfloor x\rfloor\right)= \sum_{l=0}^{d-1}\tau_{1,l,d}(x)f\left(x-\frac{l}{d}\right),
\end{equation} 
where 1 stands for the one by one identity matrix. The above is a simple interpolation formula: multiplying by $\tau_{1,l,d}(x)$ we guarantee that ``we get a contribution'' only if $x=l (\mod d)$. 

We can carry out explicitly the substitution $P_{\alpha_1,\dots,\alpha_{n-1}|C_1}(\gamma-t\alpha)$ in the right hand side of \refeq{eqAlg2Step2} using
\begin{eqnarray}\label{eqInterpolationMultivariate}
\tau_{M,c,a}(\gamma-t\alpha)= \sum_{l=0}^{a-1} \tau_{1,l,a}(t)\tau_{M,c+l\alpha,a}(\gamma),
\end{eqnarray}
where 1 stands for the one by one identity matrix. Regrouping with respect to $t$ (reminder: the coordinates of $\gamma$ are free variables, $t$ is not) we get for each summand of \refeq{eqAlg2Step2}
\begin{eqnarray}\label{eqAlg2Step3}
&&\sum_{k,l} \left(\sum_{t=\lfloor t_j\rfloor+1}^{t=\lfloor t_{N-1}\rfloor} \tau_{l,m,d} (t) t^k\right) (\mathrm{some~ quasipolynomial~ in~ \gamma^i }).
\end{eqnarray}
The $\tau$-Bernoulli sum $B_k^{ l, m, d } (x) :=\sum_{t=0}^x \tau_{l,m ,d} (t)t^k, x\in\ZZ_{\geq 0} $ can be computed in closed form using quasinumbers - the $d$ different cases based on the remainder of $x~ \mod d$ we can express algebraically with quasinumbers as in \refeq{eqTauInterpolation}. We include the formula for completeness. 
\begin{eqnarray}\displaystyle
B_k^{l,m,d}(x)&=& \sum_{t=0}^x \tau_{l,m ,d} (t) t^k= \sum_{\substack{t=rd+s\\t=0}}^{t=x} \underbrace{\tau_{l,m ,d} (rd+s)}_{=\tau_{l,m ,d} (s) }~~ (rd+s)^k \notag\\&=& \sum_{s=0}^{d-1}\tau_{l,m ,d} (s)\sum_{r=0}^{\lfloor\frac{x}{d}\rfloor}( rd+s )^k\notag\\&=&\notag\sum_{s=0}^{d-1}\tau_{l,m ,d} (s) \sum_{p=0}^{d-1}\tau_{1,p,d}(x) \sum_{r=0}^{\frac{x}{d}-\frac{p}{d}}( rd+s )^k\\&=& \sum_{s=0}^{d-1}\tau_{l,m ,d} (s)   \sum_{p=0}^{d-1} \tau_{1,p,d}(x)\notag\\&& \sum_{q=0}^k \binomial{k}{q}d^qs^{k-q}B_{q} \left(\frac{x}{d}- \frac{p}{d} \right) \label{eqTauBernoulli},
\end{eqnarray} 
where $x\in\ZZ_{\geq 0}$, and $B_{q}(x)$ are the Bernoulli polynomials with rational coefficients defined by $B_{q}(x)=\sum_{t=0}^xt^q, \forall x\in\ZZ_{\geq 0} $. It is to be noted that the numbers $\tau_{l,m ,d} (s)$ in \refeq{eqTauBernoulli} are integer constants (zero or one) that can be directly substituted. Now extend the left hand side of \refeq{eqTauBernoulli} via the the expression on the right hand side over $x\in\QQ$. Call the resulting extension $\tilde B_k^{l,m,d}$. Then we can use our interpolation formula \refeq{eqTauInterpolation} and the properties of quasinumbers to compute \begin{eqnarray}\label{eqTauInterpolationFinal}
\sum_{t=\lfloor t_j\rfloor+1}^{t=\lfloor t_{N-1}\rfloor} \tau_{l,m,d} (t) t^k&=& B_k^{ l, m, d } (\lfloor t_{N-1}\rfloor)-B_k^{ l, m, d } (\lfloor t_{j}\rfloor)\notag\\&=& \sum_{s=0}^{d-1}\tau_{1,s,d}(t_{N-1})  \tilde B_k^{l,m,d}\left(t_{N-1}\right)-\notag\\&&\sum_{s=0}^{d-1}\tau_{1,s,d}(t_{j})\tilde B_k^{l,m,d}\left(t_{j}\right). 
\end{eqnarray} 
To carry out the substitutions $\tau_{1,s,d}(t_{N-1})$, $\tau_{1,s,d}(t_{j})$, we use the explicit definition \refeq{eqIntemediatePoints} and \refeq{eqInterpolationMultivariate} which completes the description of how to compute \refeq{eqAlg2Step2} in closed form. 

(ii) Suppose now $\gamma$ lies on a supporting hyperplane of $C_i$ for some $i$. Note that \refeq{eqContributionC0} and \refeq{eqContributionCj} continue to hold true, since if a point $t_j=\lfloor t_j\rfloor$ belongs to more than one chamber, it will be counted only for the first with respect to direction $\alpha_m$. Therefore the algebraic expressions, obtained by the above procedure, will be valid for $\gamma$ lying in a supporting hyperplane for a chamber $C_i$, which completes the proof.
\end{proof}

\refle{leMainLemmaAlg2} leads to an elementary plan for computing the vector partition function which we summarize  below. 
\begin{itemize}
\item[Step 0] We start with a list of one polyhedral cone determined by the first $\alpha_1,\dots,\alpha_n$ vectors (where $n$ is rank ambient vector space). We set the counting index $m$ to be $n+1$.
\item[Step 1] We find the first polyhedral cone $C_i$ over which all we have not computed $P_{\alpha_1, \dots, \alpha_m}$, such that we have computed $P_{\alpha_1, \dots, \alpha_m}$ over any chamber $C'$ that is adjacent to an entering wall for any point $\gamma\in C_i$ with respect to $\alpha_m$. If we can't find such a chamber we increment the counter $m$. If $m$ is greater than the largest index $N$ of the $\alpha_i$'s we terminate the program, else we go to the beginning of Step 1. 
\item[Step 2] We subdivide $C_i$ into pointed polyhedral sub-cones $C_1',\dots,C_l'$ such that for each $C_j'$ there there is a unique choice of exit walls as required by \refle{leMainLemmaAlg2}. We replace the chamber $C_i$ with the newly produced pointed polyhedral sub-cones.
\item[Step 3] For each newly obtained chamber $C_j'$ we compute $P_{C_j'}={P_{\alpha_1, \dots, \alpha_m}}_{|C_j'}$  in closed form via \refeq{eqAlg2Step2} and \refeq{leMainLemmaAlg2}.
\item[Step 4] We go back to Step 1. 
\end{itemize}
The fact that we can always carry out Step 2 and that the resulting number of chambers is finite is geometrically obvious from the fact that any combinatorial chamber has only finitely many faces; the detailed proof is left to the reader. This proves the following theorem.
\begin{theo}\label{thMain}
Let $\alpha_1,\dots,\alpha_N\subset\ZZ^n$ be non-zero integral vectors with nonnegative coordinates.
The above algorithm will produce a finite list of closed pointed polyhedral cones $C_1,\dots, C_M$ of full measure with $\bigcup_{ i=1}^M$ $C_i=$ $\Lambda_{\QQ\geq 0}$ $(\alpha_1,$ $\dots,\alpha_N)$ and compute the vector partition function ${P_{\alpha_1, \dots, \alpha_m}}_{|C_j}$ as a quasipolynomial over each $C_j$. The quasipolynomial  expressions are valid over the closures of the combinatorial chambers $C_j$. In particular, the expressions ${P_{\alpha_1, \dots, \alpha_N}}_{|C_i}$ and ${P_{\alpha_1, \dots, \alpha_N}}_{|C_j}$ are equal when restricted to  the intersection of $C_i$ and $C_j$.
\end{theo}
In fact the ``vector partition function'' program realizes this subdivision in arbitrary, but small enough, dimension with arbitrary, but small enough, input data\footnote{Our algorithm is of complexity $O(\binomial{N}{n})$ for a single chamber where $N$ is the number of faces and $n$ is the dimension of the vector space. This is practically applicable (and realized in the ``vector partition function'' program) on a 32-bit machine (low-class modern PC) for more than a thousand chambers only if each has up to about 30 faces in 5$^{th}$ dimension, provided we run the ``vector partition function'' program for hours. } . An explanation of the computational difficulties of solving this problem for higher dimensions and higher number of faces can be found in \cite{BaldoniLoeraVergne}.
\begin{corollary}\label{corChambers} Let $\alpha_1,\dots,\alpha_N\subset\ZZ^n$ be non-zero integral vectors with nonnegative coordinates of maximal rank.
\begin{itemize}
\item[(a)] For a fixed pointed polyhedral cone $C_i$ coming from the algorithm in the current section, the algorithm presented in \refsec{secTheVPalgNonRealized} will produce algebraic expressions valid over the entire chamber $C_i$.
\item[(b)] \cite{BBCV} Let $D_1,\dots, D_M$ be the cones spanned by $n-1$ element subsets of $\alpha_1,\dots,\alpha_N$ of rank $n-1$. Let $C_i$ be the closed polyhedral cones in which the boundaries of $D_1,\dots D_M$ split $\Lambda_{\QQ\geq 0}$ $( \alpha_{1}, \dots \alpha_{N})$. Then $P_{\alpha_1,\dots,\alpha_N}$ is a quasipolynomial over each $C_i$.
\item[(c)] (b) holds if we replace the cones $D_i$ by the cones in the cone support (\refdef{defPartFracDecomposition}) of any partial fraction decomposition of $\prod_{i=1}^N \frac{1}{1-x^\alpha_i}$.
\end{itemize}
\end{corollary}
\begin{proof}
(c) Suppose the contrary and let $p_{j},j\in J$ be polynomials giving a partial fraction decomposition for the generating function of the vector partition function (see \refdef{defPartFracDecomposition}) given by \refsec{secTheTheoreticalAlgorithm}. Let $C$ be a fixed pointed polyhedral cone defined by the walls participating in the partial fraction decomposition. Let $C_{i_1},\dots C_{i_k}$ be all pointed polyhedral cones given by the algorithm which have a full-measure intersection with $C$. \refth{thMain} gives us then quasipolynomials measuring the vector partition function over the entire $C_{i_j}$. On the other hand, we know from \refsec{secBVdecomposition} that the quasipolynomials given by the algorithm in \refsec{secTheVPalgNonRealized} measures the vector partition function over the intersection of all affine translates  $\delta + C$, where $\delta$ varies over the exponents of all monomials $x^\delta$ participating in the polynomials in the numerators of the partial fraction decomposition. Therefore, the quasipolynomial given by the first algorithm has the same value as the quasipolynomial given by the second algorithm over $\left(\bigcap_\delta( \delta+C)\right)\cap C_{i_j}$ for each $j$. By \refle{leTheQuasipolynomialLemma} the two quasipolynomials must be equal. 
\end{proof}

\section{On the periods of the Kostant partition function of $E_6$, $E_7$, $E_8$, $F_4$ and $G_2$}\label{secPeriods}
The periods for the Kostant partition function for the classical types root systems have been computed in \cite{BBCV}. For type A the period is one (i.e. ``no period''), for types B, C and D - 2. 
\begin{prop}~\label{propPeriods}
\begin{itemize}
\item[(a)] We are using \refdef{defPeriod}. The period of the Kostant partition function of a simple root system $\Delta$ is a divisor of the least common multiple of all possible orders of the elements of the group $\Lambda_{ \ZZ}(\Delta)/\Lambda_{\ZZ}(\Delta_1)$, where $\Delta_1$ runs over all root subsystems of $\Delta$ that have the same rank as $\Delta$.
\item[(b)]  The possible orders and the  upper bounds for the Kostant partition function periods are given in the below table.

\noindent\begin{tabular}{ccccccccc}
\multicolumn {9}{l}{Root system} \\ $A_n$&$B_n$&$C_n$&$D_n$& $E_6$&$E_7$&$E_8$&$F_4$& $G_2$\\ \\
\multicolumn {9}{l}{Possible orders of elements in $\Lambda_{\ZZ}(\Delta)/\Lambda_\ZZ(\Delta_1)$} \\1&1,2&1,2&1,2&1,2,3&1,2,3,4&1,2,3,4,5,6&1,2,3,4&1,2,3\\ \\
\multicolumn {9}{l}{Kostant function period is a divisor of} \\1&2&2&2&6&12&60&12&6\\
\end{tabular}
\end{itemize}
\end{prop}
\begin{proof}

(a)A priori, the algorithm in \refsec{secTheTheoreticalAlgorithm} gives us a partial fraction decomposition for which the fractions in \refeq{eqPartFracDecomposition} have fraction support (see \refsec{secTheTheoreticalAlgorithm} for definition of fraction support) consisting not necessarily of roots, but of root multiples instead. However for a root system, according to \refle{leCanChooseCoeffOne} below, in Step 1, \refsec{secTheTheoreticalAlgorithm} of our algorithm, we can always choose our ``preferred linear combination'' such that, after eventual re-indexing,  $a_0=1$. Therefore we can always get a partial fraction decomposition with denominators corresponding to roots. 

Let $\alpha_1,\dots \alpha_n$ be $n$ linearly independent roots, where $n$ is the rank of the root system, and let $\Delta_1$ be the root system generated by the $\alpha_i$'s. Let $\prod_i\left(\frac{1}{1-x^{\alpha_i}}\right) = \sum_{\gamma \in \Lambda_{\ZZ}(\Delta)} \tau_{dB^{ -1},0,d} (\gamma) x^\gamma $, where $\Lambda_\ZZ(\Delta)$ is the lattice generated by our root system, and where in the expression $\tau_{dB^{-1},0,d}$ we have that $B$ is the matrix $(\alpha_i^j)$, 0 stands for the zero vector and $d$ is the corresponding period. By the remarks after \refeq{eqDens1}, $d$ is the least common multiple of the orders of the elements of $\Lambda_\ZZ(\Delta)/\Lambda_\ZZ \{\Delta_1\}$. By \refeq{eqtheVPfunction} we see that the period corresponding to the partial fraction decomposition summands involving the roots $\alpha_i$ is $d$.    

(b) Let $\Delta_1$ be a root subsystem of $\Delta$ of the same rank as $\Delta$. Let $\alpha\in\Lambda_\ZZ(\Delta) \backslash \Lambda_\ZZ( \Delta_1)$ be a representative of an element of $\Lambda_\ZZ(\Delta)/\Lambda_\ZZ(\Delta_1)$. Modifying $\alpha$ by elements of $\Delta_1$ does not change the represented element, therefore we can assume that $\alpha$ is $\Delta_1$-lowest. All possibilities for $\Delta_1$-lowest vectors for the exceptional Lie algebras are given in the tables of \refle{leCanChooseCoeffOne} below, and by the linear relations there, we see that their orders are as in the given table. If $\Delta$ is of type $A$ then $\Delta_1=\Delta$ by the $\mathrm{rank~}\Delta_1 = \mathrm{rank~} \Delta $ requirement. We leave cases $B$, $C$ and $D$ to the reader with the note that by the rank requirement and \cite[Table 9]{Dynkin:Semisimple} $\Delta_1$ cannot have components of type $A$.
\end{proof}

The following lemma is key to \refprop{propPeriods}. It can be derived as a consequence of \cite[Lemma 4]{Tumarkin:SimpleVertex}. The present proof does not rely on \cite{Tumarkin:SimpleVertex} and was discovered independently by the author.
\begin{lemma}\label{leCanChooseCoeffOne}
Let $\Delta^+$ be the positive root system of a simple root system of arbitrary type (A,B,C,D,E,F or G) of rank $n$. Let $\alpha_1,\dots,\alpha_k\in \Delta^+$ be $k$ linearly dependent elements of $\Delta^+$. Then there exists a linear relation between the $\alpha_i$ for which one of the roots is with coefficient +1 or -1.
\end{lemma}
\begin{proof}
Throughout this proof, whenever we are considering exceptional Lie algebras (of types $E_6$, $E_7$, $E_8$, $F_4$ and $G_2$), the coordinates are picked in a simple root basis so that the corresponding Cartan matrices are as in \cite[Table 1,page 59]{Humpreys}, see \\ \url{http://www.liegroups.org/dissemination/spherical/explorer/rootSystem.cgi} for actual printouts of the root systems.

Take $s$ linearly dependent vectors, say $\alpha_1,\dots,\alpha_s$ such that any $s-1$ of them are linearly independent. Take $a_i\in\ZZ$ such that 
\begin{equation}\label{eqleCanChoose1_first}
a_1\alpha_1+\dots+a_s\alpha_s=0
\end{equation}
and so that $a_s$ is of minimal possible absolute value (i.e. choose a minimal linear relation). Let $\Delta_1$ be the root system generated by $\alpha_1,\dots,\alpha_{s-1}$, and let $\beta_1,\dots,\beta_{s-1}$ be a simple basis of $\Delta_1$ with respect to a chosen regular height function $h$, which we will specify only for type A and the exceptional Lie algebras. Then the roots $\beta_i$ are integral linear combinations of the roots $\alpha_1,\dots, \alpha_s$. \refeq{eqleCanChoose1_first} becomes 
\begin{equation}\label{eqleCanChoose1_second}
b_1\beta_1+\dots+ b_{s-1} \beta_{s-1} +a_s \alpha_s=0 
\end{equation}
for some integers $b_i$. Since $a_s\alpha_s=a_s(\alpha_s-\beta_i)+a_s\beta_i$ we see that the possible values of $a_s$ do not change if we change $\alpha_s$ by the $\beta_i$, so we can assume that $\alpha_s$ is the lowest weight with respect to $\Delta_1$ (i.e. that $\alpha_s-\beta_i$ is not a root $\forall i$). We note that the collection of roots $\Delta_1\cup \{\alpha_s\}$ cannot split since otherwise \refeq{eqleCanChoose1_first} would not be minimal and let us call $\Delta$ the simple root system generated by $\Delta_1\cup \{\alpha_s\}$. \cite[Table 9]{Dynkin:Semisimple} gives all possible maximal (but not necessarily of maximal rank) root subsystems up to conjugation for $\Delta_1$ when $\Delta$ is of type A, B, C and D (note that there are two cases for type D) and \cite[Table 10]{Dynkin:Semisimple}- when $\Delta$ is of exceptional type. We will reproduce the tables within our case by case study. We will treat the the case for $\Delta$ of type A separately from types B, C and D, and we will treat each exceptional case in its own right. 
\begin{itemize}
\item $\Delta$ is of type $A_{s-1}$. $\Delta_1$ is then of type $A_{k_1}+\dots+A_{k_l}$ with $\sum_r (k_r+1)=s$, and the requirement that $\Delta_1$ be of rank $s-1$ forces $\Delta_1$ to be of type $A_{s-1}$. Choose a regular height function $h$ so that $\alpha_s$ is the $h$-lowest weight. Then $\alpha_s+\beta_1+\dots+\beta_{s-1}=0$ and we have the desired linear combination with $a_s=\pm 1$.
\item $\Delta$ is of type $B_{s-1}, C_{s-1}$ or $D_{s-1}$. Then $\Delta_1$ is of type $A_{k_1}+\dots+A_{k_l}+$ algebras of types B, C and D with $\sum_r (k_r+1)+$ the rest of the ranks $=s-1$. The requirement that $\Delta_1$ be of rank $s-1$ forces all $k_i=0$. In order to avoid the tedious exploration of all possible cases (and there are quite a few), we will make a Cartan matrix determinant trick. Let $\langle\bullet,\bullet\rangle_K$ denote the scalar product on the roots induced by the Killing form and normalized so that the long roots have length $\sqrt{2}$ for types B and D and length 2 for type C (this way the most predominant root length is $\sqrt{2}$). For the time being, assume that \refeq{eqleCanChoose1_second} is a linear combination over $\QQ$ (instead of over $\ZZ$). To avoid confusion rename temporarily the coefficients in  \refeq{eqleCanChoose1_second} to $\tilde a_s$, $\tilde b_i$. Take the scalar products $\langle \beta_i, \bullet\rangle_K$ to obtain
\[
\left(\begin{array}{cccc} \langle\beta_1,\beta_1\rangle_K&\dots&\langle\beta_{s-1},\beta_{s-1}\rangle_K&\langle\beta_1,\alpha_{s}\rangle_K \\ \vdots&&\vdots&\vdots\\ \langle\beta_{s-1},\beta_1\rangle_K&\dots&\langle\beta_{s-1},\beta_{s-1}\rangle_K&\langle\beta_{s-1},\alpha_{s}\rangle_K \end{array} \right)\left(\begin{array}{c} \tilde b_1\\\vdots\\ \tilde b_{s-1}\\\tilde a_s\end{array} \right) =0.
\]
Now set $\tilde a_s=1$ and solve for the rest of the $\tilde b_i$'s. Then 
\begin{equation}\label{eqleCanChoose1third}
\begin{array}{lc}\left(\begin{array}{c}\tilde b_1\\ \vdots\\\tilde b_{s-1} \end{array} \right)&=-
{\underbrace{\left(\begin{array}{ccc} \langle\beta_1, \beta_1\rangle_K&\dots&\langle\beta_1,\beta_{s-1}\rangle_K \\ \vdots & &\vdots\\ \langle\beta_{s-1}, \beta_1\rangle_K& \dots&\langle\beta_{s-1}, \beta_{s-1} \rangle_K \end{array} \right)}_{=:M}~~~}^{-1} \\& \left(\begin{array}{c} \langle\beta_1,\alpha_s\rangle_K\\\vdots\\ \langle\beta_{s-1},\alpha_s\rangle_K\end{array} \right).\end{array}
\end{equation}
Now $M$ is a matrix corresponding to the simple basis $\beta_1,\dots,\beta_{s-1}$, which is block-diagonal (up to permuting the simple basis vectors $\beta_i$) with blocks of types B, C and/or D. The determinants of those matrices are proportional to the determinants of the Cartan matrices depending on the choice (normalization) of $\langle\bullet,\bullet\rangle_K$. For our particular choice of normalization, the determinants in question are 1 for type B and 4 for types C and D. Therefore, the denominators in the rational numbers $\tilde b_i$ are 1, 2 or 4. Thus we can rescale \refeq{eqleCanChoose1third} to an integral linear combination for which $a_s$ is a divisor of 4. So far we have used nothing specific about our choice of the index $s$; therefore we can apply the same reasoning to obtain that $a_i$ is a divisor of 4 for all $i$. Now at least one $a_i$ must be $\pm 1$ because otherwise all coefficients in  our linear relation \refeq{eqleCanChoose1_first} would be even, in contradiction with its minimality.  
\item $\Delta$ is of type $E_6$. All possibilities for $\Delta_1\neq \Delta$ come in two conjugacy classes (\cite[Table 10]{Dynkin:Semisimple}). Since all roots of $E_6$ are of equal length, we can choose the regular height function $h$ on the roots so that $\alpha_7$ is a $\Delta_1$-lowest weight (i.e. $\alpha_7-\beta$ is not a root for any $h$-positive root $\beta$). We present the Dynkin diagrams of the two conjugacy classes for $\Delta_1\neq \Delta$, an explicit simple coordinate realization of both, a list of all possibilities for $\alpha_7$, and an explicit linear relation giving us the possible values for $a_7$ and the $b_i$.
\setlength{\unitlength}{1.7 pt}
\begin{center}
\begin{tabular}{|cc|}\hline
\multicolumn{2}{|c|}{Case $3A_2$}\\
\begin{picture}(0,41)
\put(1,5){\circle{1}}
\put(3,5){\tiny{6}}
\put(1,5.5){\line(0,1){6}}
\put(1,12){\circle{1}}
\put(3,12){\tiny{5}}
\put(1,19){\circle{1}}
\put(3,19){\tiny{4}}
\put(1,19.5){\line(0,1){6}}
\put(1,26){\circle{1}}
\put(3,26){\tiny{3}}
\put(1,33){\circle{1}}
\put(3,33){\tiny{2}}
\put(1,33.5){\line(0,1){6}}
\put(1,40){\circle{1}}
\put(3,40){\tiny{1}}
\end{picture}
&
\raisebox{35pt}{$\begin{array}{ccccccc}
\beta_1=(0,& 0,& 0,& 0,& 0,& 1)\\
\beta_2=(1,& 1,& 1,& 1,& 1,& 0)\\
\beta_3=(1,& 0,& 1,& 1,& 0,& 0)\\
\beta_4=(0,& 1,& 1,& 2,& 2,& 1)\\
\beta_5=(0,& 0,& 1,& 0,& 0,& 0)\\
\beta_6=(0,& 1,& 0,& 1,& 0,& 0)\\
\end{array}$}\\
$\Delta_1$-lowest $\alpha_7$&Linear combination\\
(-1,-1,-2,-2,-1,-1)&$3\alpha_7+2\beta_1+\beta_2+2\beta_3+\beta_4+2\beta_5+\beta_6=0$\\
(-1,-2,-2,-3,-2,-1)&$3\alpha_7+\beta_1+2\beta_2+\beta_3+2\beta_4+\beta_5+2\beta_6=0$\\\hline
\end{tabular}

\begin{tabular}{|cc|}\hline
\multicolumn{2}{|c|}{Case $A_5+A_1$}\\
\begin{picture}(0,41)
\put(1,40){\circle{1}}
\put(3,40){\tiny{1}}
\put(1,39.5){\line(0,-1){6}}
\put(1,33){\circle{1}}
\put(3,33){\tiny{2}}
\put(1,32.5){\line(0,-1){6}}
\put(1,26){\circle{1}}
\put(3,26){\tiny{3}}
\put(1,25.5){\line(0,-1){6}}
\put(1,19){\circle{1}}
\put(3,19){\tiny{4}}
\put(1,18.5){\line(0,-1){6}}
\put(1,12){\circle{1}}
\put(3,12){\tiny{5}}
\put(1,5){\circle{1}}
\put(3,5){\tiny{6}}
\end{picture}& \raisebox{35pt}{$\begin{array}{ccccccc}
\beta_1=&(1,&0,&0,&0,&0,&0)\\
\beta_2=&(0,&0,&1,&0,&0,&0)\\
\beta_3=&(0,&0,&0,&1,&0,&0)\\
\beta_4=&(0,&0,&0,&0,&1,&0)\\
\beta_5=&(0,&0,&0,&0,&0,&1)\\
\beta_6=&(1,&2,&2,&3,&2,&1)\end{array}$}
\\
$\Delta_1$-lowest $\alpha_7$& Linear combination\\
(-1,-1,-2,-3,-2, -1) & $2\alpha_7 +\beta_1+2\beta_2+3\beta_3+2\beta_4+\beta_5+\beta_6=0$\\
\hline
\end{tabular}
\end{center}

Suppose $\Delta_1$ is of type $3A_2$. Then $\alpha_1,\dots,\alpha_6$ must each belong to exactly one copy of $A_2$. By the following mini-lemma we get the desired linear combination with a coefficient $\pm 1$ for some $\alpha_i$.

\textbf{Mini-Lemma.} Let $\beta_1$, $\beta_2$ be a simple basis of A$_2$. Suppose $\alpha_1,\alpha_2$ are roots of A$_2$ such that $\beta_1$ and $\beta_2$ are linear combinations of $\alpha_1$ and $\alpha_2$. (a) Let $\beta_1+ 2\beta_2 =a_1\alpha_1+a_2\alpha_2$. (b) Let $\pm\beta_1\pm \beta_2 =a_1\alpha_1+a_2\alpha_2$. Then either $a_1$ or $a_2$ is $\pm 1$. 

\textbf{Proof.} (a) Depending on the angle between $\alpha_1$ and $\alpha_{2}$ we get either $\beta_1= \alpha_{1 \mathrm{~or~}2}$ and $\beta_{2}= \alpha_{2 \mathrm{~resp.~} 1}$ or $\beta_{1\mathrm{~or~}2} = \pm(\alpha_{ 1} - \alpha_{2})$ and $\beta_{2 \mathrm{~resp.~} 1} = \alpha_ {2\mathrm{~or~}1}$. In all possible cases, we get either $a_1$ or $a_2$ to be $\pm 1$. We leave (b) to the reader. $\Box$  

If for some $\alpha_i$ the root system $\Delta_{8-i}$ generated by $\alpha_j, j\neq i$, is of type $3A_2$ we are done by repeating the preceding considerations for the index $i$. Similarly we are done if $\Delta_{8-i}$ is of type $E_6$ - then $\alpha_i$ is an integral linear combination of the rest of the roots. Suppose finally $\Delta_{8-i}$ are all of type $A_5+A_1$. Then choosing different height functions we see from the second table that all $a_i$ must be equal to 2, in contradiction with the minimality of \refeq{eqleCanChoose1_first}.  
 
\item  $\Delta$ is of type $E_7$. All possibilities for $\Delta_1\neq\Delta$ come in six conjugacy classes (\cite[Table 10]{Dynkin:Semisimple}). Since all roots of $E_7$ are of equal length, we can choose the regular height function $h$ on the roots so that $\alpha_8$ is a $\Delta_1$-lowest weight (i.e. $\alpha_8-\beta$ is not a root for any $h$-positive root $\beta$). 
\begin{center}
\begin{tabular}{|cc|}\hline
\multicolumn{2}{|c|}{Case $A_2+A_5$}\\
\begin{picture}(0,50)
\put(1,5){\circle{1}}
\put(3,5){\tiny{7}}
\put(1,5.5){\line(0,1){6}}
\put(1,12){\circle{1}}
\put(3,12){\tiny{6}}
\put(1,12.5){\line(0,1){6}}
\put(1,19){\circle{1}}
\put(3,19){\tiny{5}}
\put(1,19.5){\line(0,1){6}}
\put(1,26){\circle{1}}
\put(3,26){\tiny{4}}
\put(1,26.5){\line(0,1){6}}
\put(1,33){\circle{1}}
\put(3,33){\tiny{3}}
\put(1,40){\circle{1}}
\put(3,40){\tiny{2}}
\put(1,40.5){\line(0,1){6}}
\put(1,47){\circle{1}}
\put(3,47){\tiny{1}}
\end{picture}
&
\raisebox{39pt}{$\begin{array}{ccccccc}
\beta_1=(1,& 0,& 0,& 0,& 0,& 0,& 0)\\
\beta_2=(0,& 0,& 1,& 0,& 0,& 0,& 0)\\
\beta_3=(0,& 0,& 0,& 0,& 1,& 0,& 0)\\
\beta_4=(0,& 0,& 0,& 0,& 0,& 1,& 0)\\
\beta_5=(0,& 0,& 0,& 0,& 0,& 0,& 1)\\
\beta_6=(1,& 1,& 2,& 3,& 2,& 1,& 0)\\
\beta_7=(0,& 1,& 0,& 0,& 0,& 0,& 0)\\
\end{array}$}
\\
$\Delta_1$-lowest $\alpha_8$&Linear combination\\
(-1,-1,-2,-2,-2,-2,-1)&$3\alpha_8+\beta_1+2\beta_2+2\beta_3+4\beta_4+3\beta_5+2\beta_6+\beta_7=0$\\
(-2,-2,-3,-4,-3,-2,-1)&$3\alpha_8+2\beta_1+\beta_2+\beta_3+2\beta_4+3\beta_5+4\beta_6+2\beta_7=0$\\\hline
\end{tabular}
\end{center}
By direct observation our mini-lemma from the E$_6$ case applies to the vectors $\beta_1$ and $\beta_2$ and since two of the vectors $\alpha_i$ must be in the A$_2$ component we are done.

\begin{center}
\setlength{\unitlength}{1.7 pt}
\begin{tabular}{|cc|}\hline
\multicolumn{2}{|c|}{Case $D_6+A_1$}\\
\begin{picture}(0,41)
\put(8,40){\circle{1}}
\put(10,40){\tiny{2}}
\put(-6,40){\circle{1}}
\put(-4,40){\tiny{1}}
\put(1.35,33.35){\line(1,1){6}}
\put(0.65,33.35){\line(-1,1){6}}
\put(1,33){\circle{1}}
\put(3,33){\tiny{3}}
\put(1,32.5){\line(0,-1){6}}
\put(1,26){\circle{1}}
\put(3,26){\tiny{4}}
\put(1,25.5){\line(0,-1){6}}
\put(1,19){\circle{1}}
\put(3,19){\tiny{5}}
\put(1,18.5){\line(0,-1){6}}
\put(1,12){\circle{1}}
\put(3,12){\tiny{6}}
\put(1,5){\circle{1}}
\put(3,5){\tiny{7}}
\end{picture}& \raisebox{35pt}{$\begin{array}{ccccccc}
\beta_1=(1,& 0,& 0,& 0,& 0,& 0,& 0)\\
\beta_2=(0,& 1,& 0,& 1,& 1,& 1,& 1)\\
\beta_3=(0,& 0,& 1,& 0,& 0,& 0,& 0)\\
\beta_4=(0,& 0,& 0,& 1,& 0,& 0,& 0)\\
\beta_5=(0,& 0,& 0,& 0,& 1,& 0,& 0)\\
\beta_6=(0,& 0,& 0,& 0,& 0,& 1,& 0)\\
\beta_7=(1,& 2,& 2,& 3,& 2,& 1,& 0)\\
\end{array}$}
\\
$\Delta_1$-lowest $\alpha_8$&Linear combination\\
(-2,-2,-3,-4,-3,-2,-1)&$2\alpha_8+3\beta_1+2\beta_2+4\beta_3+3\beta_4+2\beta_5+\beta_6+\beta_7=0$\\\hline
\end{tabular}
\end{center}

\begin{center}
\begin{tabular}{|cc|}\hline
\multicolumn{2}{|c|}{Case $A_3+A_1+A_3$}\\
\begin{picture}(0,50)
\put(1,5){\circle{1}}
\put(3,5){\tiny{7}}
\put(1,5.5){\line(0,1){6}}
\put(1,12){\circle{1}}
\put(3,12){\tiny{6}}
\put(1,12.5){\line(0,1){6}}
\put(1,19){\circle{1}}
\put(3,19){\tiny{5}}
\put(1,26){\circle{1}}
\put(3,26){\tiny{4}}
\put(1,33){\circle{1}}
\put(3,33){\tiny{3}}
\put(1,33.5){\line(0,1){6}}
\put(1,40){\circle{1}}
\put(3,40){\tiny{2}}
\put(1,40.5){\line(0,1){6}}
\put(1,47){\circle{1}}
\put(3,47){\tiny{1}}
\end{picture}
&
\raisebox{39pt}{$\begin{array}{ccccccc}
\beta_1=(1,& 0,& 0,& 0,& 0,& 0,& 0)\\
\beta_2=(0,& 0,& 1,& 0,& 0,& 0,& 0)\\
\beta_3=(0,& 0,& 0,& 1,& 0,& 0,& 0)\\
\beta_4=(1,& 1,& 2,& 3,& 3,& 2,& 1)\\
\beta_5=(0,& 0,& 0,& 0,& 0,& 1,& 0)\\
\beta_6=(0,& 0,& 0,& 0,& 0,& 0,& 1)\\
\beta_7=(1,& 2,& 2,& 3,& 2,& 1,& 0)\\
\end{array}$}
\\
$\Delta_1$-lowest $\alpha_8$&Linear combination\\
(-1,-1,-2,-3,-2,-2,-1)&$4\alpha_8+\beta_1+2\beta_2+3\beta_3+2\beta_4+3\beta_5+2\beta_6+\beta_7=0$\\
(-2,-2,-3,-4,-3,-2,-1)&$4\alpha_8+3\beta_1+2\beta_2+\beta_3+2\beta_4+\beta_5+2\beta_6+3\beta_7=0$\\
(-1,-1,-2,-2,-1,-1,-1)&$2\alpha_8+\beta_1+2\beta_2+\beta_3+\beta_5+2\beta_6+\beta_7=0$\\
\hline
\end{tabular}
\end{center}
\begin{center}
\begin{tabular}{|cc|}\hline
\multicolumn{2}{|c|}{Case $A_7$}\\
\begin{picture}(0,50)
\put(1,5){\circle{1}}
\put(3,5){\tiny{7}}
\put(1,5.5){\line(0,1){6}}
\put(1,12){\circle{1}}
\put(3,12){\tiny{6}}
\put(1,12.5){\line(0,1){6}}
\put(1,19){\circle{1}}
\put(3,19){\tiny{5}}
\put(1,19.5){\line(0,1){6}}
\put(1,26){\circle{1}}
\put(3,26){\tiny{4}}
\put(1,26.5){\line(0,1){6}}
\put(1,33){\circle{1}}
\put(3,33){\tiny{3}}
\put(1,33.5){\line(0,1){6}}
\put(1,40){\circle{1}}
\put(3,40){\tiny{2}}
\put(1,40.5){\line(0,1){6}}
\put(1,47){\circle{1}}
\put(3,47){\tiny{1}}
\end{picture}
&
\raisebox{39pt}{$\begin{array}{ccccccc}
\beta_1=(1,& 0,& 0,& 0,& 0,& 0,& 0)\\
\beta_2=(0,& 0,& 1,& 0,& 0,& 0,& 0)\\
\beta_3=(0,& 0,& 0,& 0,& 1,& 0,& 0)\\
\beta_4=(0,& 0,& 0,& 0,& 0,& 1,& 0)\\
\beta_5=(0,& 0,& 0,& 0,& 0,& 0,& 1)\\
\beta_6=(1,& 1,& 2,& 3,& 2,& 1,& 0)\\
\beta_7=(0,& 1,& 0,& 0,& 0,& 0,& 0)\\
\end{array}$}
\\
$\Delta_1$-lowest $\alpha_8$ &Linear combination\\
(-1,-1,-2,-3,-3,-2,-1)&$2\alpha_8+\beta_1+2\beta_2+3\beta_3+4\beta_4+3\beta_5+2\beta_6+\beta_7=0$\\\hline
\end{tabular}
\end{center}
\begin{center}
\setlength{\unitlength}{1.7 pt}
\begin{tabular}{|cc|}\hline
\multicolumn{2}{|c|}{Case $D_4+3A_1$}\\
\begin{picture}(0,41)
\put(8,40){\circle{1}}
\put(10,40){\tiny{2}}
\put(-6,40){\circle{1}}
\put(-4,40){\tiny{1}}
\put(1.35,33.35){\line(1,1){6}}
\put(0.65,33.35){\line(-1,1){6}}
\put(1,33){\circle{1}}
\put(3,33){\tiny{3}}
\put(1,32.5){\line(0,-1){6}}
\put(1,26){\circle{1}}
\put(3,26){\tiny{4}}
\put(1,19){\circle{1}}
\put(3,19){\tiny{5}}
\put(1,12){\circle{1}}
\put(3,12){\tiny{6}}
\put(1,5){\circle{1}}
\put(3,5){\tiny{7}}
\end{picture}& \raisebox{35pt}{$\begin{array}{ccccccc}
\beta_1=(1,& 0,& 0,& 0,& 0,& 0,& 0)\\
\beta_2=(0,& 0,& 0,& 1,& 0,& 0,& 0)\\
\beta_3=(0,& 0,& 1,& 0,& 0,& 0,& 0)\\
\beta_4=(0,& 1,& 0,& 1,& 1,& 1,& 1)\\
\beta_5=(1,& 1,& 2,& 3,& 3,& 2,& 1)\\
\beta_6=(0,& 0,& 0,& 0,& 0,& 1,& 0)\\
\beta_7=(1,& 2,& 2,& 3,& 2,& 1,& 0)\\
\end{array}$}
\\
$\Delta_1$-lowest $\alpha_8$&Linear combination\\
(-2,-2,-3,-4,-3,-2,-1)&$2\alpha_8+2\beta_1+\beta_2+2\beta_3+\beta_4+\beta_5+\beta_7=0$\\
(-1,-1,-2,-3,-2,-2,-1)&$2\alpha_8+\beta_1+2\beta_2+2\beta_3+\beta_4+\beta_5+\beta_6=0$\\
(-1,-2,-2,-3,-2,-2,-1)&$2\alpha_8+\beta_1+\beta_2+2\beta_3+2\beta_4+\beta_6+\beta_7=0$\\
\hline
\end{tabular}
\end{center}
\begin{center}
\begin{tabular}{|cc|}\hline
\multicolumn{2}{|c|}{Case $7A_1$}\\
\begin{picture}(0,50)
\put(1,5){\circle{1}}
\put(3,5){\tiny{7}}
\put(1,12){\circle{1}}
\put(3,12){\tiny{6}}
\put(1,19){\circle{1}}
\put(3,19){\tiny{5}}
\put(1,26){\circle{1}}
\put(3,26){\tiny{4}}
\put(1,33){\circle{1}}
\put(3,33){\tiny{3}}
\put(1,40){\circle{1}}
\put(3,40){\tiny{2}}
\put(1,47){\circle{1}}
\put(3,47){\tiny{1}}
\end{picture}
&
\raisebox{39pt}{$\begin{array}{ccccccc}
\beta_1=(1,& 0,& 0,& 0,& 0,& 0,& 0)\\
\beta_2=(0,& 0,& 0,& 1,& 0,& 0,& 0)\\
\beta_3=(0,& 1,& 0,& 1,& 1,& 1,& 1)\\
\beta_4=(1,& 1,& 2,& 3,& 3,& 2,& 1)\\
\beta_5=(0,& 0,& 0,& 0,& 0,& 1,& 0)\\
\beta_6=(1,& 2,& 2,& 3,& 2,& 1,& 0)\\
\beta_7=(1,& 1,& 2,& 2,& 1,& 1,& 1)\\
\end{array}$}
\\
$\Delta_1$-lowest $\alpha_8$ &Linear combination\\
(-1,-1,-1,-2,-1,-1,0)&$2\alpha_8+\beta_1+\beta_2+\beta_5+\beta_6=0$\\
(-1,-1,-1,-2,-1,-1,-1)&$2\alpha_8+\beta_1+\beta_2+\beta_3+\beta_7=0$\\
(-1,-1,-1,-2,-2,-2,-1)&$2\alpha_8+\beta_1+\beta_3+\beta_4+\beta_5=0$\\
(-1,-1,-2,-3,-2,-2,-1)&$2\alpha_8+\beta_2+\beta_4+\beta_5+\beta_7=0$\\
(-1,-2,-2,-3,-2,-2,-1)&$2\alpha_8+\beta_3+\beta_5+\beta_6+\beta_7=0$\\
(-1,-2,-2,-4,-3,-2,-1)&$2\alpha_8+\beta_2+\beta_3+\beta_4+\beta_6=0$\\
(-2,-2,-3,-4,-3,-2,-1)&$2\alpha_8+\beta_1+\beta_4+\beta_6+\beta_7=0$\\
\hline
\end{tabular}
\end{center}
If for some $\alpha_i$ the root system $\Delta_{9-i}$ generated by $\alpha_j, j\neq i$, is of type $A_2+A_5$ we are done by repeating the considerations from the very first case. Similarly we are done if $\Delta_{9-i}$ is of type $E_7$ - then $\alpha_i$ is an integral linear combination of the rest of the roots. Suppose finally $\Delta_{9-i}$ are all not of type $A_5+A_2$ or $E_7$ and assume $\alpha_i$ is $\Delta_{9-i}$-lowest. Then (possibly choosing different height functions) we see from the remaining 5 tables that all $a_i$ must be even, in contradiction with the minimality of \refeq{eqleCanChoose1_first}.  
\item  $\Delta$ is of type $E_8$. All possibilities for $\Delta_1\neq\Delta$ come in fourteen conjugacy classes (\cite[Table 10]{Dynkin:Semisimple}). Since all roots of $E_8$ are of equal length, we can choose the regular height function $h$ on the roots so that $\alpha_9$ is a $\Delta_1$-lowest weight (i.e. $\alpha_9-\beta$ is not a root for any $h$-positive root $\beta$). 
\begin{center}
\begin{tabular}{|cc|}\hline
\multicolumn{2}{|c|}{Case $A_8$}\\
\begin{picture}(0,57)
\put(1,5){\circle{1}}
\put(3,5){\tiny{8}}
\put(1,5.5){\line(0,1){6}}
\put(1,12){\circle{1}}
\put(3,12){\tiny{7}}
\put(1,12.5){\line(0,1){6}}
\put(1,19){\circle{1}}
\put(3,19){\tiny{6}}
\put(1,19.5){\line(0,1){6}}
\put(1,26){\circle{1}}
\put(3,26){\tiny{5}}
\put(1,26.5){\line(0,1){6}}
\put(1,33){\circle{1}}
\put(3,33){\tiny{4}}
\put(1,33.5){\line(0,1){6}}
\put(1,40){\circle{1}}
\put(3,40){\tiny{3}}
\put(1,40.5){\line(0,1){6}}
\put(1,47){\circle{1}}
\put(3,47){\tiny{2}}
\put(1,47.5){\line(0,1){6}}
\put(1,54){\circle{1}}
\put(3,54){\tiny{1}}
\end{picture}
&
\raisebox{47pt}{$\begin{array}{cccccccc}
\beta_1=(1,& 3,& 3,& 5,& 4,& 3,& 2,& 1)\\
\beta_2=(1,& 0,& 0,& 0,& 0,& 0,& 0,& 0)\\
\beta_3=(0,& 0,& 1,& 0,& 0,& 0,& 0,& 0)\\
\beta_4=(0,& 0,& 0,& 1,& 0,& 0,& 0,& 0)\\
\beta_5=(0,& 0,& 0,& 0,& 1,& 0,& 0,& 0)\\
\beta_6=(0,& 0,& 0,& 0,& 0,& 1,& 0,& 0)\\
\beta_7=(0,& 0,& 0,& 0,& 0,& 0,& 1,& 0)\\
\beta_8=(0,& 0,& 0,& 0,& 0,& 0,& 0,& 1)\\
\end{array}$}
\\
$\Delta_1$-lowest $\alpha_9$&Linear combination\\
(-1,-1,-2,-3,-3,-3,-2,-1)&$3\alpha_9+\beta_1+2\beta_2+3\beta_3+4\beta_4+5\beta_5+6\beta_6+4\beta_7+2\beta_8=0$\\
(-2,-2,-4,-5,-4,-3,-2,-1)&$3\alpha_9+2\beta_1+4\beta_2+6\beta_3+5\beta_4+4\beta_5+3\beta_6+2\beta_7+\beta_8=0$\\\hline
\end{tabular}
\end{center}
A$_8$ has 36 positive roots and there are only $\binomial{36}{8}=30260340$ possible choices of roots $\alpha_1, \dots \alpha_8$. Direct computation of all possibilities shows that for all such choices of maximal rank, there is at least one $\alpha_i$ with coefficient $\pm 1$\footnote{ All easy considerations, attempted by the author and similar to the ones given for the $E_7$ and E$_6$ cases, failed in the current case. We therefore enumerated all possible 30260340 combinations by brute computer force and verified the claim (the computation itself took less than 45 minutes). 
}.   
\begin{center}
\begin{tabular}{|cc|}\hline
\multicolumn{2}{|c|}{Case $2A_4 $}\\
\begin{picture}(0,57)
\put(1,5){\circle{1}}
\put(3,5){\tiny{8}}
\put(1,5.5){\line(0,1){6}}
\put(1,12){\circle{1}}
\put(3,12){\tiny{7}}
\put(1,12.5){\line(0,1){6}}
\put(1,19){\circle{1}}
\put(3,19){\tiny{6}}
\put(1,19.5){\line(0,1){6}}
\put(1,26){\circle{1}}
\put(3,26){\tiny{5}}
\put(1,33){\circle{1}}
\put(3,33){\tiny{4}}
\put(1,33.5){\line(0,1){6}}
\put(1,40){\circle{1}}
\put(3,40){\tiny{3}}
\put(1,40.5){\line(0,1){6}}
\put(1,47){\circle{1}}
\put(3,47){\tiny{2}}
\put(1,47.5){\line(0,1){6}}
\put(1,54){\circle{1}}
\put(3,54){\tiny{1}}
\end{picture}
&
\raisebox{47pt}{$\begin{array}{cccccccc}
\beta_1=(0,& 0,& 0,& 0,& 1,& 0,& 0,& 0)\\
\beta_2=(0,& 0,& 0,& 1,& 0,& 0,& 0,& 0)\\
\beta_3=(0,& 0,& 1,& 0,& 0,& 0,& 0,& 0)\\
\beta_4=(1,& 0,& 0,& 0,& 0,& 0,& 0,& 0)\\
\beta_5=(1,& 1,& 2,& 3,& 3,& 3,& 2,& 1)\\
\beta_6=(1,& 2,& 2,& 3,& 2,& 1,& 0,& 0)\\
\beta_7=(0,& 0,& 0,& 0,& 0,& 0,& 1,& 0)\\
\beta_8=(0,& 0,& 0,& 0,& 0,& 0,& 0,& 1)\\
\end{array}$}
\\
$\Delta_1$-lowest $\alpha_9$&Linear combination\\
(-1,-1,-2,-3,-2,-1,0,0)&$5\alpha_9+3\beta_1+6\beta_2+4\beta_3+2\beta_4+\beta_5+2\beta_6-2\beta_7-\beta_8=0$\\
(-1,-1,-2,-3,-2,-1,-1,-1)&$5\alpha_9+3\beta_1+6\beta_2+4\beta_3+2\beta_4+\beta_5+2\beta_6+3\beta_7+4\beta_8=0$\\
(-1,-1,-2,-3,-3,-2,-1,0)&$5\alpha_9+4\beta_1+3\beta_2+2\beta_3+\beta_4+3\beta_5+\beta_6-\beta_7-3\beta_8=0$\\
(-1,-1,-2,-3,-3,-2,-1,-1)&$5\alpha_9+4\beta_1+3\beta_2+2\beta_3+\beta_4+3\beta_5+\beta_6-\beta_7+2\beta_8=0$\\
(-2,-2,-3,-4,-3,-2,-1,-1)&$5\alpha_9+\beta_1+2\beta_2+3\beta_3+4\beta_4+2\beta_5+4\beta_6+\beta_7+3\beta_8=0$\\
(-2,-2,-3,-4,-3,-2,-2,-1)&$5\alpha_9+\beta_1+2\beta_2+3\beta_3+4\beta_4+2\beta_5+4\beta_6+6\beta_7+3\beta_8=0$\\
(-2,-2,-4,-5,-4,-3,-2,-1)&$5\alpha_9+2\beta_1+4\beta_2+6\beta_3+3\beta_4+4\beta_5+3\beta_6+2\beta_7+\beta_8=0$\\
(-2,-3,-4,-6,-5,-3,-2,-1)&$5\alpha_9+4\beta_1+3\beta_2+2\beta_3+\beta_4+3\beta_5+6\beta_6+4\beta_7+2\beta_8=0$\\\hline
\end{tabular}
\end{center}
$2A_4$ has 20 positive roots and there are only $\binomial{10}{4}^2=44100$ possible choices of roots $\alpha_1, \dots \alpha_8$. Direct computation of all possibilities shows that for all such choices of maximal rank, there is at least one $\alpha_i$ with coefficient $\pm 1$\footnote{We have, for a second time, summoned our computer computational power.}.
\begin{center}
\begin{tabular}{|cc|}\hline
\multicolumn{2}{|c|}{Case $A_5+A_1+A_2 $}\\
\begin{picture}(0,57)
\put(1,5){\circle{1}}
\put(3,5){\tiny{8}}
\put(1,5.5){\line(0,1){6}}
\put(1,12){\circle{1}}
\put(3,12){\tiny{7}}
\put(1,19){\circle{1}}
\put(3,19){\tiny{6}}
\put(1,26){\circle{1}}
\put(3,26){\tiny{5}}
\put(1,26.5){\line(0,1){6}}
\put(1,33){\circle{1}}
\put(3,33){\tiny{4}}
\put(1,33.5){\line(0,1){6}}
\put(1,40){\circle{1}}
\put(3,40){\tiny{3}}
\put(1,40.5){\line(0,1){6}}
\put(1,47){\circle{1}}
\put(3,47){\tiny{2}}
\put(1,47.5){\line(0,1){6}}
\put(1,54){\circle{1}}
\put(3,54){\tiny{1}}
\end{picture}
&
\raisebox{47pt}{$\begin{array}{cccccccc}
\beta_1=(0,& 0,& 0,& 0,& 0,& 1,& 0,& 0)\\
\beta_2=(0,& 0,& 0,& 0,& 1,& 0,& 0,& 0)\\
\beta_3=(0,& 0,& 0,& 1,& 0,& 0,& 0,& 0)\\
\beta_4=(0,& 0,& 1,& 0,& 0,& 0,& 0,& 0)\\
\beta_5=(1,& 0,& 0,& 0,& 0,& 0,& 0,& 0)\\
\beta_6=(1,& 2,& 2,& 3,& 2,& 1,& 0,& 0)\\
\beta_7=(2,& 3,& 4,& 6,& 5,& 4,& 3,& 1)\\
\beta_8=(0,& 0,& 0,& 0,& 0,& 0,& 0,& 1)\\
\end{array}$}
\\
$\Delta_1$-lowest $\alpha_9$&Linear combination\\
(-1,-1,-2,-3,-2,-1,0,0)&$2\alpha_9+\beta_1+2\beta_2+3\beta_3+2\beta_4+\beta_5+\beta_6=0$\\
(-1,-1,-2,-3,-3,-2,-1,0)&$3\alpha_9+2\beta_1+4\beta_2+3\beta_3+2\beta_4+\beta_5+\beta_7-\beta_8=0$\\
(-1,-1,-2,-3,-3,-2,-1,-1)&$3\alpha_9+2\beta_1+4\beta_2+3\beta_3+2\beta_4+\beta_5+\beta_7+2\beta_8=0$\\
(-2,-2,-3,-4,-3,-2,-1,-1)&$6\alpha_9+\beta_1+2\beta_2+3\beta_3+4\beta_4+5\beta_5+3\beta_6+2\beta_7+4\beta_8=0$\\
(-2,-2,-4,-5,-4,-3,-2,-1)&$3\alpha_9+\beta_1+2\beta_2+3\beta_3+4\beta_4+2\beta_5+2\beta_7+\beta_8=0$\\
(-2,-3,-4,-6,-5,-4,-2,-1)&$6\alpha_9+5\beta_1+4\beta_2+3\beta_3+2\beta_4+\beta_5+3\beta_6+4\beta_7+2\beta_8=0$\\
\hline
\end{tabular}
\end{center}
For the cases for which the coefficient in front of $\alpha_9$ is 3, the mini-lemma from the $E_6$ case applies to the roots $\beta_7$ and $\beta_8$ to produce a relation in the $\alpha_i$'s with coefficient $\pm 1$. The case when the coefficient in front of $\alpha_9$ is even is treated together with the most of the remaining cases around at the end of the $E_8$ proof. 
\begin{center}
\begin{tabular}{|cc|}\hline
\multicolumn{2}{|c|}{Case $A_2+E_6$}\\
\begin{picture}(0,50)
\put(1,5){\circle{1}}
\put(3,5){\tiny{8}}
\put(1,5.5){\line(0,1){6}}
\put(1,12){\circle{1}}
\put(3,12){\tiny{7}}
\put(1,12.5){\line(0,1){6}}
\put(1,19){\circle{1}}
\put(3,19){\tiny{6}}
\put(1,19.5){\line(0,1){6}}
\put(1,26){\circle{1}} \put(8,19){\circle{1}}
\put(3,26){\tiny{5}} \put(10,19){\tiny{3}}
\put(1,26.5){\line(0,1){6}}\put(1.5,19){\line(1,0){6}}
\put(1,33){\circle{1}}
\put(3,33){\tiny{4}}
\put(1,40){\circle{1}}
\put(3,40){\tiny{2}}
\put(1,40.5){\line(0,1){6}}
\put(1,47){\circle{1}}
\put(3,47){\tiny{1}}
\end{picture}
&
\raisebox{39pt}{$\begin{array}{cccccccc}
\beta_1=(2,& 3,& 4,& 6,& 5,& 4,& 3,& 1)\\
\beta_2=(0,& 0,& 0,& 0,& 0,& 0,& 0,& 1)\\
\beta_3=(1,& 0,& 0,& 0,& 0,& 0,& 0,& 0)\\
\beta_4=(0,& 1,& 0,& 0,& 0,& 0,& 0,& 0)\\
\beta_5=(0,& 0,& 1,& 0,& 0,& 0,& 0,& 0)\\
\beta_6=(0,& 0,& 0,& 1,& 0,& 0,& 0,& 0)\\
\beta_7=(0,& 0,& 0,& 0,& 1,& 0,& 0,& 0)\\
\beta_8=(0,& 0,& 0,& 0,& 0,& 1,& 0,& 0)\\
\end{array}$}
\\
$\Delta_1$-lowest $\alpha_9$&Linear combination\\
(-1,-1,-2,-3,-2,-1,-1,-1)&$3\alpha_9+\beta_1+2\beta_2+\beta_3+2\beta_5+3\beta_6+\beta_7-\beta_8=0$\\
(-1,-1,-2,-2,-2,-2,-1,-1)&$3\alpha_9+\beta_1+2\beta_2+\beta_3+2\beta_5+\beta_7+2\beta_8=0$\\
(-1,-1,-2,-3,-3,-2,-1,0)&$3\alpha_9+\beta_1-\beta_2+\beta_3+2\beta_5+3\beta_6+4\beta_7+2\beta_8=0$\\
(-1,-1,-2,-3,-3,-2,-1,-1)&$3\alpha_9+\beta_1+2\beta_2+\beta_3+2\beta_5+3\beta_6+4\beta_7+2\beta_8=0$\\
(-1,-2,-2,-3,-3,-2,-1,-1)&$3\alpha_9+\beta_1+2\beta_2+\beta_3+3\beta_4+2\beta_5+3\beta_6+4\beta_7+2\beta_8=0$\\
(-2,-2,-3,-4,-3,-2,-1,-1)&$3\alpha_9+\beta_1+2\beta_2+4\beta_3+3\beta_4+5\beta_5+6\beta_6+4\beta_7+2\beta_8=0$\\
(-2,-2,-3,-4,-3,-2,-2,-1)&$3\alpha_9+2\beta_1+\beta_2+2\beta_3+\beta_5-\beta_7-2\beta_8=0$\\
(-2,-2,-3,-4,-4,-3,-2,-1)&$3\alpha_9+2\beta_1+\beta_2+2\beta_3+\beta_5+2\beta_7+\beta_8=0$\\
(-2,-2,-4,-5,-4,-3,-2,-1)&$3\alpha_9+2\beta_1+\beta_2+2\beta_3+4\beta_5+3\beta_6+2\beta_7+\beta_8=0$\\
(-2,-3,-4,-5,-4,-3,-2,-1)&$3\alpha_9+2\beta_1+\beta_2+2\beta_3+3\beta_4+4\beta_5+3\beta_6+2\beta_7+\beta_8=0$\\
(-2,-3,-4,-6,-5,-3,-2,-1)&$3\alpha_9+2\beta_1+\beta_2+2\beta_3+3\beta_4+4\beta_5+6\beta_6+5\beta_7+\beta_8=0$\\
(-2,-3,-4,-6,-5,-4,-2,-1)&$3\alpha_9+2\beta_1+\beta_2+2\beta_3+3\beta_4+4\beta_5+6\beta_6+5\beta_7+4\beta_8=0$\\
\hline
\end{tabular}
\end{center}
The mini-lemma from the $E_6$ case applies to the roots $\beta_1$ and $\beta_2$ to produce a relation that contains a coefficient $\pm 1$.
\begin{center}
\begin{tabular}{|cc|}\hline
\multicolumn{2}{|c|}{Case $4A_2 $}\\
\begin{picture}(0,57)
\put(1,5){\circle{1}}
\put(3,5){\tiny{8}}
\put(1,5.5){\line(0,1){6}}
\put(1,12){\circle{1}}
\put(3,12){\tiny{7}}
\put(1,19){\circle{1}}
\put(3,19){\tiny{6}}
\put(1,19.5){\line(0,1){6}}
\put(1,26){\circle{1}}
\put(3,26){\tiny{5}}
\put(1,33){\circle{1}}
\put(3,33){\tiny{4}}
\put(1,33.5){\line(0,1){6}}
\put(1,40){\circle{1}}
\put(3,40){\tiny{3}}
\put(1,47){\circle{1}}
\put(3,47){\tiny{2}}
\put(1,47.5){\line(0,1){6}}
\put(1,54){\circle{1}}
\put(3,54){\tiny{1}}
\end{picture}
&
\raisebox{47pt}{$\begin{array}{cccccccc}
\beta_1=(0,& 0,& 0,& 0,& 0,& 1,& 0,& 0)\\
\beta_2=(0,& 0,& 0,& 0,& 1,& 0,& 0,& 0)\\
\beta_3=(1,& 1,& 2,& 3,& 2,& 1,& 0,& 0)\\
\beta_4=(0,& 1,& 0,& 0,& 0,& 0,& 0,& 0)\\
\beta_5=(0,& 0,& 1,& 0,& 0,& 0,& 0,& 0)\\
\beta_6=(1,& 0,& 0,& 0,& 0,& 0,& 0,& 0)\\
\beta_7=(2,& 3,& 4,& 6,& 5,& 4,& 3,& 1)\\
\beta_8=(0,& 0,& 0,& 0,& 0,& 0,& 0,& 1)\\
\end{array}$}
\\
$\Delta_1$-lowest $\alpha_9$&Linear combination\\
(-1,-1,-1,-1,-1,-1,0,0)&$3\alpha_9+2\beta_1+\beta_2+\beta_3+2\beta_4+\beta_5+2\beta_6=0$\\
(-1,-1,-2,-2,-2,-1,0,0)&$3\alpha_9+\beta_1+2\beta_2+2\beta_3+\beta_4+2\beta_5+\beta_6=0$\\
(-1,-1,-2,-3,-2,-1,-1,-1)&$3\alpha_9-2\beta_1-\beta_2+\beta_3-\beta_4+\beta_7+2\beta_8=0$\\
(-1,-1,-2,-2,-2,-2,-1,-1)&$3\alpha_9+2\beta_1+\beta_2+2\beta_5+\beta_6+\beta_7+2\beta_8=0$\\
(-1,-1,-2,-3,-3,-2,-1,0)&$3\alpha_9+\beta_1+2\beta_2+\beta_3-\beta_4+\beta_7-\beta_8=0$\\
(-1,-1,-2,-3,-3,-2,-1,-1)&$3\alpha_9+\beta_1+2\beta_2+\beta_3-\beta_4+\beta_7+2\beta_8=0$\\
(-1,-2,-2,-3,-3,-2,-1,-1)&$3\alpha_9+\beta_1+2\beta_2+\beta_3+2\beta_4+\beta_7+2\beta_8=0$\\
(-2,-2,-3,-4,-3,-2,-1,-1)&$3\alpha_9+2\beta_3+\beta_4+\beta_5+2\beta_6+\beta_7+2\beta_8=0$\\
(-2,-2,-3,-4,-3,-2,-2,-1)&$3\alpha_9-2\beta_1-\beta_2+\beta_5+2\beta_6+2\beta_7+\beta_8=0$\\
(-2,-2,-3,-4,-4,-3,-2,-1)&$3\alpha_9+\beta_1+2\beta_2+\beta_5+2\beta_6+2\beta_7+\beta_8=0$\\
(-2,-2,-4,-5,-4,-3,-2,-1)&$3\alpha_9+\beta_3-\beta_4+2\beta_5+\beta_6+2\beta_7+\beta_8=0$\\
(-2,-3,-4,-5,-4,-3,-2,-1)&$3\alpha_9+\beta_3+2\beta_4+2\beta_5+\beta_6+2\beta_7+\beta_8=0$\\
(-2,-3,-4,-6,-5,-3,-2,-1)&$3\alpha_9-\beta_1+\beta_2+2\beta_3+\beta_4+2\beta_7+\beta_8=0$\\
(-2,-3,-4,-6,-5,-4,-2,-1)&$3\alpha_9+2\beta_1+\beta_2+2\beta_3+\beta_4+2\beta_7+\beta_8=0$\\\hline
\end{tabular}
\end{center}
The mini-lemma from the $E_6$ case applies to the roots $\beta_1$ and $\beta_2$ to produce a relation that contains a coefficient $\pm 1$.

\begin{center}
\begin{tabular}{|cc|}\hline
\multicolumn{2}{|c|}{Case $D_8$}\\
\begin{picture}(0,50)
\put(1,5){\circle{1}}
\put(3,5){\tiny{8}}
\put(1,5.5){\line(0,1){6}}
\put(1,12){\circle{1}}
\put(3,12){\tiny{7}}
\put(1,12.5){\line(0,1){6}}
\put(1,19){\circle{1}}
\put(3,19){\tiny{6}}
\put(1,19.5){\line(0,1){6}}
\put(1,26){\circle{1}}
\put(3,26){\tiny{5}}
\put(1,26.5){\line(0,1){6}}
\put(1,33){\circle{1}}
\put(3,33){\tiny{4}}
\put(1,33.5){\line(0,1){6}}
\put(1,40){\circle{1}}
\put(3,40){\tiny{3}}
\put(1.35,40.35){\line(1,1){6}}
\put(0.65,40.35){\line(-1,1){6}}
\put(8,47){\circle{1}}
\put(10,47){\tiny{2}}
\put(-6,47){\circle{1}}
\put(-4,47){\tiny{1}}
\end{picture}
&
\raisebox{39pt}{$\begin{array}{cccccccc}
\beta_1=(1,& 2,& 2,& 3,& 2,& 1,& 0,& 0)\\
\beta_2=(0,& 0,& 0,& 0,& 0,& 0,& 0,& 1)\\
\beta_3=(0,& 0,& 0,& 0,& 0,& 0,& 1,& 0)\\
\beta_4=(0,& 0,& 0,& 0,& 0,& 1,& 0,& 0)\\
\beta_5=(0,& 0,& 0,& 0,& 1,& 0,& 0,& 0)\\
\beta_6=(0,& 0,& 0,& 1,& 0,& 0,& 0,& 0)\\
\beta_7=(0,& 0,& 1,& 0,& 0,& 0,& 0,& 0)\\
\beta_8=(1,& 0,& 0,& 0,& 0,& 0,& 0,& 0)\\
\end{array}$}
\\
$\Delta_1$-lowest $\alpha_9$&Linear combination\\
(-2,-3,-4,-6,-5,-4,-3,-2)&$2\alpha_9+3\beta_1+4\beta_2+6\beta_3+5\beta_4+4\beta_5+3\beta_6+2\beta_7+\beta_8=0$\\\hline
\end{tabular}
\end{center}

\begin{center}
\begin{tabular}{|cc|}\hline
\multicolumn{2}{|c|}{Case $ A_1+A_7$}\\
\begin{picture}(0,57)
\put(1,5){\circle{1}}
\put(3,5){\tiny{8}}
\put(1,5.5){\line(0,1){6}}
\put(1,12){\circle{1}}
\put(3,12){\tiny{7}}
\put(1,12.5){\line(0,1){6}}
\put(1,19){\circle{1}}
\put(3,19){\tiny{6}}
\put(1,19.5){\line(0,1){6}}
\put(1,26){\circle{1}}
\put(3,26){\tiny{5}}
\put(1,26.5){\line(0,1){6}}
\put(1,33){\circle{1}}
\put(3,33){\tiny{4}}
\put(1,33.5){\line(0,1){6}}
\put(1,40){\circle{1}}
\put(3,40){\tiny{3}}
\put(1,40.5){\line(0,1){6}}
\put(1,47){\circle{1}}
\put(3,47){\tiny{2}}
\put(1,54){\circle{1}}
\put(3,54){\tiny{1}}
\end{picture}
&
\raisebox{47pt}{$\begin{array}{cccccccc}
\beta_1=(2,& 3,& 4,& 6,& 5,& 4,& 3,& 2)\\
\beta_2=(1,& 2,& 2,& 3,& 2,& 1,& 0,& 0)\\
\beta_3=(0,& 0,& 0,& 0,& 0,& 0,& 1,& 0)\\
\beta_4=(0,& 0,& 0,& 0,& 0,& 1,& 0,& 0)\\
\beta_5=(0,& 0,& 0,& 0,& 1,& 0,& 0,& 0)\\
\beta_6=(0,& 0,& 0,& 1,& 0,& 0,& 0,& 0)\\
\beta_7=(0,& 0,& 1,& 0,& 0,& 0,& 0,& 0)\\
\beta_8=(1,& 0,& 0,& 0,& 0,& 0,& 0,& 0)\\
\end{array}$}
\\
$\Delta_1$-lowest $\alpha_9$&Linear combination\\
(-1,-1,-2,-3,-3,-2,-1,0)&$2\alpha_9+\beta_2+2\beta_3+3\beta_4+4\beta_5+3\beta_6+2\beta_7+\beta_8=0$\\
(-2,-2,-4,-5,-4,-3,-2,-1)&$4\alpha_9+2\beta_1+\beta_2+2\beta_3+3\beta_4+4\beta_5+5\beta_6+6\beta_7+3\beta_8=0$\\
(-2,-3,-4,-6,-5,-4,-3,-1)&$4\alpha_9+2\beta_1+3\beta_2+6\beta_3+5\beta_4+4\beta_5+3\beta_6+2\beta_7+\beta_8=0$\\
\hline
\end{tabular}
\end{center}

\begin{center}
\begin{tabular}{|cc|}\hline
\multicolumn{2}{|c|}{Case $A_1+E_7$}\\
\begin{picture}(0,50)
\put(1,5){\circle{1}}
\put(3,5){\tiny{8}}
\put(1,5.5){\line(0,1){6}}
\put(1,12){\circle{1}}
\put(3,12){\tiny{7}}
\put(1,12.5){\line(0,1){6}}
\put(1,19){\circle{1}}
\put(3,19){\tiny{6}}
\put(1,19.5){\line(0,1){6}}
\put(1,26){\circle{1}} \put(8,26){\circle{1}}
\put(3,26){\tiny{5}} \put(10,26){\tiny{3}}
\put(1,26.5){\line(0,1){6}}\put(1.5,26){\line(1,0){6}}
\put(1,33){\circle{1}}
\put(3,33){\tiny{4}}
\put(1,33.5){\line(0,1){6}}
\put(1,40){\circle{1}}
\put(3,40){\tiny{2}}
\put(1,47){\circle{1}}
\put(3,47){\tiny{1}}
\end{picture}
&
\raisebox{39pt}{$\begin{array}{cccccccc}
\beta_1=(2,& 3,& 4,& 6,& 5,& 4,& 3,& 1)\\
\beta_2=(1,& 0,& 0,& 0,& 0,& 0,& 0,& 0)\\
\beta_3=(0,& 1,& 0,& 0,& 0,& 0,& 0,& 0)\\
\beta_4=(0,& 0,& 1,& 0,& 0,& 0,& 0,& 0)\\
\beta_5=(0,& 0,& 0,& 1,& 0,& 0,& 0,& 0)\\
\beta_6=(0,& 0,& 0,& 0,& 1,& 0,& 0,& 0)\\
\beta_7=(0,& 0,& 0,& 0,& 0,& 1,& 0,& 0)\\
\beta_8=(0,& 0,& 0,& 0,& 0,& 0,& 1,& 1)\\
\end{array}$}
\\
$\Delta_1$-lowest $\alpha_9$&Linear combination\\
(-1,-1,-2,-3,-3,-2,-1,0)&$2\alpha_9+\beta_1-\beta_3+\beta_6-\beta_8=0$\\
(-2,-2,-3,-4,-3,-2,-2,-1)&$2\alpha_9+\beta_1+2\beta_2+\beta_3+2\beta_4+2\beta_5+\beta_6+\beta_8=0$\\
(-2,-2,-3,-4,-4,-3,-2,-1)&$2\alpha_9+\beta_1+2\beta_2+\beta_3+2\beta_4+2\beta_5+3\beta_6+2\beta_7+\beta_8=0$\\
(-2,-2,-4,-5,-4,-3,-2,-1)&$2\alpha_9+\beta_1+2\beta_2+\beta_3+4\beta_4+4\beta_5+3\beta_6+2\beta_7+\beta_8=0$\\
(-2,-3,-4,-5,-4,-3,-2,-1)&$2\alpha_9+\beta_1+2\beta_2+3\beta_3+4\beta_4+4\beta_5+3\beta_6+2\beta_7+\beta_8=0$\\
(-2,-3,-4,-6,-5,-3,-2,-1)&$2\alpha_9+\beta_1+2\beta_2+3\beta_3+4\beta_4+6\beta_5+5\beta_6+2\beta_7+\beta_8=0$\\
(-2,-3,-4,-6,-5,-4,-2,-1)&$2\alpha_9+\beta_1+2\beta_2+3\beta_3+4\beta_4+6\beta_5+5\beta_6+4\beta_7+\beta_8=0$\\
(-2,-3,-4,-6,-5,-4,-3,-2)&$2\alpha_9+\beta_1+2\beta_2+3\beta_3+4\beta_4+6\beta_5+5\beta_6+4\beta_7+3\beta_8=0$\\\hline
\end{tabular}
\end{center}

\begin{center}
\begin{tabular}{|cc|}\hline
\multicolumn{2}{|c|}{Case $D_6+2A_1$}\\
\begin{picture}(0,50)
\put(1,5){\circle{1}}
\put(3,5){\tiny{8}}

\put(1,12){\circle{1}}
\put(3,12){\tiny{7}}

\put(1,19){\circle{1}}
\put(3,19){\tiny{6}}
\put(1,19.5){\line(0,1){6}}
\put(1,26){\circle{1}}
\put(3,26){\tiny{5}}
\put(1,26.5){\line(0,1){6}}
\put(1,33){\circle{1}}
\put(3,33){\tiny{4}}
\put(1,33.5){\line(0,1){6}}
\put(1,40){\circle{1}}
\put(3,40){\tiny{3}}
\put(1.35,40.35){\line(1,1){6}}
\put(0.65,40.35){\line(-1,1){6}}
\put(8,47){\circle{1}}
\put(10,47){\tiny{2}}
\put(-6,47){\circle{1}}
\put(-4,47){\tiny{1}}
\end{picture}
&
\raisebox{39pt}{$\begin{array}{cccccccc}
\beta_1=(1,& 0,& 0,& 0,& 0,& 0,& 0,& 0)\\
\beta_2=(0,& 1,& 0,& 1,& 1,& 1,& 1,& 0)\\
\beta_3=(0,& 0,& 1,& 0,& 0,& 0,& 0,& 0)\\
\beta_4=(0,& 0,& 0,& 1,& 0,& 0,& 0,& 0)\\
\beta_5=(0,& 0,& 0,& 0,& 1,& 0,& 0,& 0)\\
\beta_6=(0,& 0,& 0,& 0,& 0,& 1,& 0,& 0)\\
\beta_7=(2,& 3,& 4,& 6,& 5,& 4,& 3,& 2)\\
\beta_8=(1,& 2,& 2,& 3,& 2,& 1,& 0,& 0)\\
\end{array}$}
\\
$\Delta_1$-lowest $\alpha_9$&Linear combination\\
(-1,-1,-1,-1,-1,-1,0,0)&$2\alpha_9+\beta_1-\beta_4+\beta_6+\beta_8=0$\\
(-1,-1,-2,-2,-2,-1,0,0)&$2\alpha_9+\beta_1+2\beta_3+\beta_4+2\beta_5+\beta_6+\beta_8=0$\\
(-1,-1,-2,-3,-2,-1,0,0)&$2\alpha_9+\beta_1+2\beta_3+3\beta_4+2\beta_5+\beta_6+\beta_8=0$\\
(-1,-1,-2,-3,-2,-1,-1,-1)&$2\alpha_9-\beta_2+\beta_4-\beta_6+\beta_7=0$\\
(-1,-1,-2,-2,-2,-2,-1,-1)&$2\alpha_9-\beta_2-\beta_4+\beta_6+\beta_7=0$\\
(-1,-1,-2,-3,-3,-2,-1,-1)&$2\alpha_9-\beta_2+\beta_4+2\beta_5+\beta_6+\beta_7=0$\\
(-1,-2,-2,-3,-3,-2,-1,-1)&$2\alpha_9-\beta_1-\beta_2-2\beta_3-2\beta_4+\beta_7+\beta_8=0$\\
(-2,-2,-3,-4,-3,-2,-1,0)&$2\alpha_9+3\beta_1+2\beta_2+4\beta_3+3\beta_4+2\beta_5+\beta_6+\beta_8=0$\\
(-2,-2,-3,-4,-3,-2,-1,-1)&$2\alpha_9+\beta_1-\beta_2+\beta_7+\beta_8=0$\\
(-2,-2,-3,-4,-3,-2,-2,-1)&$2\alpha_9+2\beta_1+\beta_2+2\beta_3+\beta_4-\beta_6+\beta_7=0$\\
(-2,-2,-3,-4,-4,-3,-2,-1)&$2\alpha_9+2\beta_1+\beta_2+2\beta_3+\beta_4+2\beta_5+\beta_6+\beta_7=0$\\
(-2,-2,-4,-5,-4,-3,-2,-1)&$2\alpha_9+2\beta_1+\beta_2+4\beta_3+3\beta_4+2\beta_5+\beta_6+\beta_7=0$\\
(-2,-3,-4,-5,-4,-3,-2,-1)&$2\alpha_9+\beta_1+\beta_2+2\beta_3+\beta_7+\beta_8=0$\\
(-2,-3,-4,-6,-5,-3,-2,-1)&$2\alpha_9+\beta_1+\beta_2+2\beta_3+2\beta_4+2\beta_5+\beta_7+\beta_8=0$\\
(-2,-3,-4,-6,-5,-4,-2,-1)&$2\alpha_9+\beta_1+\beta_2+2\beta_3+2\beta_4+2\beta_5+2\beta_6+\beta_7+\beta_8=0$\\
(-2,-3,-4,-6,-5,-4,-3,-1)&$2\alpha_9+2\beta_1+3\beta_2+4\beta_3+3\beta_4+2\beta_5+\beta_6+\beta_7=0$\\
\hline
\end{tabular}
\end{center}
\begin{center}
\begin{tabular}{|cc|}\hline
\multicolumn{2}{|c|}{Case $D_5+A_3$}\\
\begin{picture}(0,50)
\put(1,5){\circle{1}}
\put(3,5){\tiny{8}}
\put(1,12){\circle{1}}
\put(3,12){\tiny{7}}
\put(1,19){\circle{1}}
\put(3,19){\tiny{6}}
\put(1,26){\circle{1}}
\put(3,26){\tiny{5}}
\put(1,26.5){\line(0,1){6}}
\put(1,33){\circle{1}}
\put(3,33){\tiny{4}}
\put(1,33.5){\line(0,1){6}}
\put(1,40){\circle{1}}
\put(3,40){\tiny{3}}
\put(1.35,40.35){\line(1,1){6}}
\put(0.65,40.35){\line(-1,1){6}}
\put(8,47){\circle{1}}
\put(10,47){\tiny{2}}
\put(-6,47){\circle{1}}
\put(-4,47){\tiny{1}}
\end{picture}
&
\raisebox{39pt}{$\begin{array}{cccccccc}
\beta_1=(1,& 0,& 0,& 0,& 0,& 0,& 0,& 0)\\
\beta_2=(0,& 1,& 0,& 1,& 1,& 1,& 1,& 0)\\
\beta_3=(0,& 0,& 1,& 0,& 0,& 0,& 0,& 0)\\
\beta_4=(0,& 0,& 0,& 1,& 0,& 0,& 0,& 0)\\
\beta_5=(0,& 0,& 0,& 0,& 1,& 0,& 0,& 0)\\
\beta_6=(1,& 1,& 2,& 3,& 3,& 3,& 2,& 1)\\
\beta_7=(1,& 2,& 2,& 3,& 2,& 1,& 0,& 0)\\
\beta_8=(0,& 0,& 0,& 0,& 0,& 0,& 1,& 1)\\
\end{array}$}
\\
$\Delta_1$-lowest $\alpha_9$&Linear combination\\
(-1,-1,-1,-1,-1,-1,0,0)&$4\alpha_9+\beta_1-\beta_2-2\beta_3-4\beta_4-2\beta_5+\beta_6+2\beta_7-\beta_8=0$\\
(-1,-1,-2,-2,-2,-1,0,0)&$4\alpha_9+\beta_1-\beta_2+2\beta_3+2\beta_5+\beta_6+2\beta_7-\beta_8=0$\\
(-1,-1,-2,-3,-2,-1,0,0)&$4\alpha_9+\beta_1-\beta_2+2\beta_3+4\beta_4+2\beta_5+\beta_6+2\beta_7-\beta_8=0$\\
(-1,-1,-2,-3,-2,-1,-1,-1)&$4\alpha_9+\beta_1-\beta_2+2\beta_3+4\beta_4+2\beta_5+\beta_6+2\beta_7+3\beta_8=0$\\
(-1,-1,-2,-2,-2,-2,-1,-1)&$4\alpha_9-\beta_1-3\beta_2-2\beta_3-4\beta_4-2\beta_5+3\beta_6+2\beta_7+\beta_8=0$\\
(-1,-1,-2,-3,-3,-2,-1,0)&$2\alpha_9+\beta_1+\beta_2+2\beta_3+2\beta_4+2\beta_5+\beta_6-\beta_8=0$\\
(-1,-1,-2,-3,-3,-2,-1,-1)&$4\alpha_9-\beta_1-3\beta_2-2\beta_3+2\beta_5+3\beta_6+2\beta_7+\beta_8=0$\\
(-1,-2,-2,-3,-3,-2,-1,-1)&$2\alpha_9-\beta_1-\beta_2-2\beta_3-2\beta_4+\beta_6+2\beta_7+\beta_8=0$\\
(-2,-2,-3,-4,-3,-2,-1,0)&$4\alpha_9+5\beta_1+3\beta_2+6\beta_3+4\beta_4+2\beta_5+\beta_6+2\beta_7-\beta_8=0$\\
(-2,-2,-3,-4,-3,-2,-1,-1)&$2\alpha_9+\beta_1-\beta_2+\beta_6+2\beta_7+\beta_8=0$\\
(-2,-2,-3,-4,-3,-2,-2,-1)&$4\alpha_9+5\beta_1+3\beta_2+6\beta_3+4\beta_4+2\beta_5+\beta_6+2\beta_7+3\beta_8=0$\\
(-2,-2,-3,-4,-4,-3,-2,-1)&$4\alpha_9+3\beta_1+\beta_2+2\beta_3+2\beta_5+3\beta_6+2\beta_7+\beta_8=0$\\
(-2,-2,-4,-5,-4,-3,-2,-1)&$4\alpha_9+3\beta_1+\beta_2+6\beta_3+4\beta_4+2\beta_5+3\beta_6+2\beta_7+\beta_8=0$\\
(-2,-3,-4,-5,-4,-3,-2,-1)&$2\alpha_9+\beta_1+\beta_2+2\beta_3+\beta_6+2\beta_7+\beta_8=0$\\
(-2,-3,-4,-6,-5,-3,-2,-1)&$2\alpha_9+\beta_1+\beta_2+2\beta_3+2\beta_4+2\beta_5+\beta_6+2\beta_7+\beta_8=0$\\
(-2,-3,-4,-6,-5,-4,-3,-1)&$4\alpha_9+3\beta_1+5\beta_2+6\beta_3+4\beta_4+2\beta_5+3\beta_6+2\beta_7+\beta_8=0$\\
\hline
\end{tabular}
\end{center}
\begin{center}
\begin{tabular}{|cc|}\hline
\multicolumn{2}{|c|}{Case $2D_4 $}\\
\begin{picture}(0,43)

\put(1,5){\circle{1}}
\put(3,5){\tiny{8}}
\put(1,5.5){\line(0,1){6}}
\put(1,12){\circle{1}}
\put(3,12){\tiny{7}}
\put(1.35,12.35){\line(1,1){6}}
\put(0.65,12.35){\line(-1,1){6}}
\put(8,19){\circle{1}}
\put(10,19){\tiny{6}}
\put(-6,19){\circle{1}}
\put(-4,19){\tiny{5}}

\put(1,26){\circle{1}}
\put(3,26){\tiny{4}}
\put(1,26.5){\line(0,1){6}}
\put(1,33){\circle{1}}
\put(3,33){\tiny{3}}
\put(1.35,33.35){\line(1,1){6}}
\put(0.65,33.35){\line(-1,1){6}}
\put(8,40){\circle{1}}
\put(10,40){\tiny{2}}
\put(-6,40){\circle{1}}
\put(-4,40){\tiny{1}}
\end{picture}
&
\raisebox{32pt}{$\begin{array}{cccccccc}
\beta_1=(1,& 0,& 0,& 0,& 0,& 0,& 0,& 0)\\
\beta_2=(0,& 1,& 0,& 1,& 1,& 1,& 1,& 0)\\
\beta_3=(0,& 0,& 1,& 0,& 0,& 0,& 0,& 0)\\
\beta_4=(0,& 0,& 0,& 1,& 0,& 0,& 0,& 0)\\
\beta_5=(1,& 1,& 2,& 3,& 3,& 2,& 1,& 0)\\
\beta_6=(0,& 0,& 0,& 0,& 0,& 1,& 0,& 0)\\
\beta_7=(0,& 0,& 0,& 0,& 0,& 0,& 1,& 1)\\
\beta_8=(1,& 2,& 2,& 3,& 2,& 1,& 0,& 0)\\
\end{array}$}
\\
$\Delta_1$-lowest $\alpha_9$&Linear combination\\
(-1,-1,-1,-1,-1,-1,0,0)&$2\alpha_9+\beta_1-\beta_4+\beta_6+\beta_8=0$\\
(-1,-1,-2,-2,-2,-1,0,0)&$2\alpha_9-\beta_2-\beta_4+\beta_5+\beta_8=0$\\
(-1,-1,-2,-3,-2,-1,0,0)&$2\alpha_9-\beta_2+\beta_4+\beta_5+\beta_8=0$\\
(-1,-1,-2,-3,-2,-1,-1,-1)&$2\alpha_9-\beta_2+\beta_4+\beta_5+2\beta_7+\beta_8=0$\\
(-1,-1,-2,-2,-2,-2,-1,-1)&$2\alpha_9-\beta_2-\beta_4+\beta_5+2\beta_6+2\beta_7+\beta_8=0$\\
(-1,-1,-2,-3,-3,-2,-1,-1)&$2\alpha_9-\beta_1-2\beta_2-2\beta_3-\beta_4+2\beta_5+\beta_6+2\beta_7+\beta_8=0$\\
(-1,-2,-2,-3,-3,-2,-1,-1)&$2\alpha_9-\beta_1-\beta_2-2\beta_3-2\beta_4+\beta_5+\beta_6+2\beta_7+2\beta_8=0$\\
(-2,-2,-3,-4,-3,-2,-1,0)&$2\alpha_9+2\beta_1+\beta_2+2\beta_3+\beta_4+\beta_5+\beta_8=0$\\
(-2,-2,-3,-4,-3,-2,-1,-1)&$2\alpha_9+\beta_1-\beta_2+\beta_5+\beta_6+2\beta_7+2\beta_8=0$\\
(-2,-2,-3,-4,-3,-2,-2,-1)&$2\alpha_9+2\beta_1+\beta_2+2\beta_3+\beta_4+\beta_5+2\beta_7+\beta_8=0$\\
(-2,-2,-3,-4,-3,-3,-2,-1)&$2\alpha_9+2\beta_1+\beta_2+2\beta_3+\beta_4+\beta_5+2\beta_6+2\beta_7+\beta_8=0$\\
(-2,-2,-3,-4,-4,-3,-2,-1)&$2\alpha_9+\beta_1-\beta_4+2\beta_5+\beta_6+2\beta_7+\beta_8=0$\\
(-2,-2,-4,-5,-4,-3,-2,-1)&$2\alpha_9+\beta_1+2\beta_3+\beta_4+2\beta_5+\beta_6+2\beta_7+\beta_8=0$\\
(-2,-3,-4,-5,-4,-3,-2,-1)&$2\alpha_9+\beta_1+\beta_2+2\beta_3+\beta_5+\beta_6+2\beta_7+2\beta_8=0$\\
(-2,-3,-4,-6,-4,-3,-2,-1)&$2\alpha_9+\beta_1+\beta_2+2\beta_3+2\beta_4+\beta_5+\beta_6+2\beta_7+2\beta_8=0$\\
(-2,-3,-4,-6,-5,-4,-3,-1)&$2\alpha_9+\beta_1+2\beta_2+2\beta_3+\beta_4+2\beta_5+\beta_6+2\beta_7+\beta_8=0$\\
\hline
\end{tabular}
\end{center}
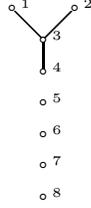
\begin{center}
\begin{tabular}{|cc|}\hline
\multicolumn{2}{|c|}{Case $D_4+4A_1 $}\\
\begin{picture}(0,50)
\put(1,5){\circle{1}}
\put(3,5){\tiny{8}}
\put(1,12){\circle{1}}
\put(3,12){\tiny{7}}
\put(1,19){\circle{1}}
\put(3,19){\tiny{6}}
\put(1,26){\circle{1}}
\put(3,26){\tiny{5}}
\put(1,33){\circle{1}}
\put(3,33){\tiny{4}}
\put(1,33.5){\line(0,1){6}}
\put(1,40){\circle{1}}
\put(3,40){\tiny{3}}
\put(1.35,40.35){\line(1,1){6}}
\put(0.65,40.35){\line(-1,1){6}}
\put(8,47){\circle{1}}
\put(10,47){\tiny{2}}
\put(-6,47){\circle{1}}
\put(-4,47){\tiny{1}}
\end{picture}
&
\raisebox{39pt}{$\begin{array}{cccccccc}
\beta_1=(0,& 1,& 0,& 1,& 1,& 1,& 1,& 0)\\
\beta_2=(1,& 1,& 2,& 3,& 3,& 2,& 1,& 0)\\
\beta_3=(0,& 0,& 0,& 0,& 0,& 0,& 0,& 1)\\
\beta_4=(1,& 1,& 2,& 2,& 1,& 1,& 1,& 0)\\
\beta_5=(0,& 0,& 0,& 1,& 0,& 0,& 0,& 0)\\
\beta_6=(1,& 0,& 0,& 0,& 0,& 0,& 0,& 0)\\
\beta_7=(1,& 2,& 2,& 3,& 2,& 1,& 0,& 0)\\
\beta_8=(0,& 0,& 0,& 0,& 0,& 1,& 0,& 0)\\
\end{array}$}
\\
$\Delta_1$-lowest $\alpha_9$&Linear combination\\
(-1,-1,-1,-1,-1,-1,0,0)&$2\alpha_9-\beta_5+\beta_6+\beta_7+\beta_8=0$\\
(-1,-1,-1,-2,-1,-1,0,0)&$2\alpha_9+\beta_5+\beta_6+\beta_7+\beta_8=0$\\
(-1,-1,-2,-2,-2,-1,0,0)&$2\alpha_9-\beta_1+\beta_2-\beta_5+\beta_7=0$\\
(-1,-1,-2,-3,-2,-1,0,0)&$2\alpha_9-\beta_1+\beta_2+\beta_5+\beta_7=0$\\
(-1,-1,-2,-3,-2,-1,-1,-1)&$2\alpha_9+\beta_2+2\beta_3+\beta_4+\beta_5-\beta_8=0$\\
(-1,-1,-2,-2,-2,-2,-1,-1)&$2\alpha_9+\beta_2+2\beta_3+\beta_4-\beta_5+\beta_8=0$\\
(-1,-2,-2,-3,-3,-2,-1,-1)&$2\alpha_9+\beta_1+\beta_2+2\beta_3-\beta_5+\beta_7=0$\\
(-2,-2,-3,-4,-3,-2,-1,0)&$2\alpha_9+\beta_2+\beta_4+\beta_6+\beta_7=0$\\
(-2,-2,-3,-4,-3,-2,-1,-1)&$2\alpha_9+\beta_2+2\beta_3+\beta_4+\beta_6+\beta_7=0$\\
(-1,-2,-2,-4,-3,-3,-2,-1)&$2\alpha_9+2\beta_1+\beta_2+2\beta_3+\beta_4+\beta_5+\beta_8=0$\\
(-2,-2,-3,-4,-3,-2,-2,-1)&$2\alpha_9+\beta_1+\beta_2+2\beta_3+2\beta_4+\beta_6-\beta_8=0$\\
(-2,-2,-3,-4,-3,-3,-2,-1)&$2\alpha_9+\beta_1+\beta_2+2\beta_3+2\beta_4+\beta_6+\beta_8=0$\\
(-2,-2,-3,-4,-4,-3,-2,-1)&$2\alpha_9+\beta_1+2\beta_2+2\beta_3+\beta_4-\beta_5+\beta_6=0$\\
(-2,-2,-3,-5,-4,-3,-2,-1)&$2\alpha_9+\beta_1+2\beta_2+2\beta_3+\beta_4+\beta_5+\beta_6=0$\\
(-2,-3,-3,-5,-4,-3,-2,-1)&$2\alpha_9+2\beta_1+\beta_2+2\beta_3+\beta_4+\beta_6+\beta_7=0$\\
(-2,-3,-4,-5,-4,-3,-2,-1)&$2\alpha_9+\beta_1+\beta_2+2\beta_3+2\beta_4-\beta_5+\beta_7=0$\\
(-2,-3,-4,-6,-4,-3,-2,-1)&$2\alpha_9+\beta_1+\beta_2+2\beta_3+2\beta_4+\beta_5+\beta_7=0$\\
(-2,-3,-4,-6,-5,-3,-2,-1)&$2\alpha_9+\beta_1+2\beta_2+2\beta_3+\beta_4+\beta_7-\beta_8=0$\\
(-2,-3,-4,-6,-5,-4,-2,-1)&$2\alpha_9+\beta_1+2\beta_2+2\beta_3+\beta_4+\beta_7+\beta_8=0$\\
\hline
\end{tabular}
\end{center}
\begin{center}
\begin{tabular}{|cc|}\hline
\multicolumn{2}{|c|}{Case $2A_3+2A_1 $}\\
\begin{picture}(0,57)
\put(1,5){\circle{1}}
\put(3,5){\tiny{8}}
\put(1,12){\circle{1}}
\put(3,12){\tiny{7}}
\put(1,19){\circle{1}}
\put(3,19){\tiny{6}}
\put(1,19.5){\line(0,1){6}}
\put(1,26){\circle{1}}
\put(3,26){\tiny{5}}
\put(1,26.5){\line(0,1){6}}
\put(1,33){\circle{1}}
\put(3,33){\tiny{4}}
\put(1,40){\circle{1}}
\put(3,40){\tiny{3}}
\put(1,40.5){\line(0,1){6}}
\put(1,47){\circle{1}}
\put(3,47){\tiny{2}}
\put(1,47.5){\line(0,1){6}}
\put(1,54){\circle{1}}
\put(3,54){\tiny{1}}
\end{picture}
&
\raisebox{47pt}{$\begin{array}{cccccccc}
\beta_1=(1,& 2,& 2,& 3,& 2,& 2,& 2,& 1)\\
\beta_2=(1,& 1,& 2,& 3,& 3,& 2,& 1,& 0)\\
\beta_3=(0,& 0,& 0,& 0,& 0,& 0,& 0,& 1)\\
\beta_4=(0,& 0,& 0,& 1,& 0,& 0,& 0,& 0)\\
\beta_5=(0,& 0,& 1,& 0,& 0,& 0,& 0,& 0)\\
\beta_6=(1,& 0,& 0,& 0,& 0,& 0,& 0,& 0)\\
\beta_7=(1,& 2,& 2,& 3,& 2,& 1,& 0,& 0)\\
\beta_8=(0,& 0,& 0,& 0,& 0,& 1,& 0,& 0)\\
\end{array}$}
\\
$\Delta_1$-lowest $\alpha_9$&Linear combination\\
(-1,-1,-1,-1,-1,-1,0,0)&$2\alpha_9-\beta_4+\beta_6+\beta_7+\beta_8=0$\\
(-1,-1,-1,-2,-1,-1,0,0)&$2\alpha_9+\beta_4+\beta_6+\beta_7+\beta_8=0$\\
(-1,-1,-2,-2,-1,-1,0,0)&$2\alpha_9+\beta_4+2\beta_5+\beta_6+\beta_7+\beta_8=0$\\
(-1,-1,-2,-2,-2,-1,0,0)&$4\alpha_9-\beta_1+2\beta_2+\beta_3-\beta_4+2\beta_5+\beta_6+2\beta_7=0$\\
(-1,-1,-2,-3,-2,-1,0,0)&$4\alpha_9-\beta_1+2\beta_2+\beta_3+3\beta_4+2\beta_5+\beta_6+2\beta_7=0$\\
(-1,-1,-2,-3,-2,-1,-1,-1)&$4\alpha_9+\beta_1+2\beta_2+3\beta_3+3\beta_4+2\beta_5+\beta_6-2\beta_8=0$\\
(-1,-1,-2,-2,-2,-2,-1,-1)&$4\alpha_9+\beta_1+2\beta_2+3\beta_3-\beta_4+2\beta_5+\beta_6+2\beta_8=0$\\
(-1,-1,-2,-3,-2,-2,-1,-1)&$4\alpha_9+\beta_1+2\beta_2+3\beta_3+3\beta_4+2\beta_5+\beta_6+2\beta_8=0$\\
(-1,-2,-2,-3,-3,-2,-1,-1)&$4\alpha_9+\beta_1+2\beta_2+3\beta_3-3\beta_4-2\beta_5-\beta_6+2\beta_7=0$\\
(-2,-2,-3,-4,-3,-2,-1,0)&$4\alpha_9+\beta_1+2\beta_2-\beta_3+\beta_4+2\beta_5+3\beta_6+2\beta_7=0$\\
(-2,-2,-3,-4,-3,-2,-1,-1)&$4\alpha_9+\beta_1+2\beta_2+3\beta_3+\beta_4+2\beta_5+3\beta_6+2\beta_7=0$\\
(-1,-2,-2,-4,-3,-3,-2,-1)&$4\alpha_9+3\beta_1+2\beta_2+\beta_3+\beta_4-2\beta_5-\beta_6+2\beta_8=0$\\
(-2,-2,-3,-4,-3,-2,-2,-1)&$4\alpha_9+3\beta_1+2\beta_2+\beta_3+\beta_4+2\beta_5+3\beta_6-2\beta_8=0$\\
(-2,-2,-3,-4,-3,-3,-2,-1)&$4\alpha_9+3\beta_1+2\beta_2+\beta_3+\beta_4+2\beta_5+3\beta_6+2\beta_8=0$\\
(-2,-2,-3,-4,-4,-3,-2,-1)&$2\alpha_9+\beta_1+2\beta_2+\beta_3-\beta_4+\beta_6=0$\\
(-2,-2,-3,-5,-4,-3,-2,-1)&$2\alpha_9+\beta_1+2\beta_2+\beta_3+\beta_4+\beta_6=0$\\
(-2,-3,-3,-5,-4,-3,-2,-1)&$4\alpha_9+3\beta_1+2\beta_2+\beta_3-\beta_4-2\beta_5+\beta_6+2\beta_7=0$\\
(-2,-2,-4,-5,-4,-3,-2,-1)&$2\alpha_9+\beta_1+2\beta_2+\beta_3+\beta_4+2\beta_5+\beta_6=0$\\
(-2,-3,-4,-5,-4,-3,-2,-1)&$4\alpha_9+3\beta_1+2\beta_2+\beta_3-\beta_4+2\beta_5+\beta_6+2\beta_7=0$\\
(-2,-3,-4,-6,-4,-3,-2,-1)&$4\alpha_9+3\beta_1+2\beta_2+\beta_3+3\beta_4+2\beta_5+\beta_6+2\beta_7=0$\\
(-2,-3,-4,-6,-5,-3,-2,-1)&$2\alpha_9+\beta_1+2\beta_2+\beta_3+\beta_7-\beta_8=0$\\
(-2,-3,-4,-6,-5,-4,-2,-1)&$2\alpha_9+\beta_1+2\beta_2+\beta_3+\beta_7+\beta_8=0$\\
\hline
\end{tabular}
\end{center}
\begin{center}
\begin{tabular}{|cc|}\hline
\multicolumn{2}{|c|}{Case $8A_1 $}\\
\begin{picture}(0,57)
\put(1,5){\circle{1}}
\put(3,5){\tiny{8}}
\put(1,12){\circle{1}}
\put(3,12){\tiny{7}}
\put(1,19){\circle{1}}
\put(3,19){\tiny{6}}
\put(1,26){\circle{1}}
\put(3,26){\tiny{5}}
\put(1,33){\circle{1}}
\put(3,33){\tiny{4}}
\put(1,40){\circle{1}}
\put(3,40){\tiny{3}}
\put(1,47){\circle{1}}
\put(3,47){\tiny{2}}
\put(1,54){\circle{1}}
\put(3,54){\tiny{1}}
\end{picture}
&
\raisebox{47pt}{$\begin{array}{cccccccc}
\beta_1=(0,& 1,& 0,& 1,& 1,& 1,& 1,& 0)\\
\beta_2=(1,& 1,& 2,& 3,& 3,& 2,& 1,& 0)\\
\beta_3=(1,& 1,& 2,& 2,& 1,& 1,& 1,& 0)\\
\beta_4=(0,& 0,& 0,& 1,& 0,& 0,& 0,& 0)\\
\beta_5=(1,& 0,& 0,& 0,& 0,& 0,& 0,& 0)\\
\beta_6=(1,& 2,& 2,& 3,& 2,& 1,& 0,& 0)\\
\beta_7=(0,& 0,& 0,& 0,& 0,& 1,& 0,& 0)\\
\beta_8=(2,& 3,& 4,& 6,& 5,& 4,& 3,& 2)\\
\end{array}$}
\\
$\Delta_1$-lowest $\alpha_9$&Linear combination\\
(-1,-1,-1,-1,-1,-1,0,0)&$2\alpha_9-\beta_4+\beta_5+\beta_6+\beta_7=0$\\
(-1,-1,-1,-2,-1,-1,0,0)&$2\alpha_9+\beta_4+\beta_5+\beta_6+\beta_7=0$\\
(-1,-1,-2,-2,-1,-1,0,0)&$2\alpha_9-\beta_1+\beta_3+\beta_6+\beta_7=0$\\
(-1,-1,-1,-2,-1,-1,-1,0)&$2\alpha_9+\beta_1+\beta_3+\beta_4+\beta_5=0$\\
(-1,-1,-2,-2,-2,-1,0,0)&$2\alpha_9-\beta_1+\beta_2-\beta_4+\beta_6=0$\\
(-1,-1,-2,-3,-2,-1,0,0)&$2\alpha_9-\beta_1+\beta_2+\beta_4+\beta_6=0$\\
(-1,-1,-1,-2,-2,-2,-1,0)&$2\alpha_9+\beta_1+\beta_2+\beta_5+\beta_7=0$\\
(-1,-1,-2,-3,-2,-2,-1,0)&$2\alpha_9+\beta_2+\beta_3+\beta_4+\beta_7=0$\\
(-1,-1,-2,-3,-2,-1,-1,-1)&$2\alpha_9-\beta_1+\beta_4-\beta_7+\beta_8=0$\\
(-1,-1,-2,-2,-2,-2,-1,-1)&$2\alpha_9-\beta_1-\beta_4+\beta_7+\beta_8=0$\\
(-1,-2,-2,-3,-2,-2,-1,0)&$2\alpha_9+\beta_1+\beta_3+\beta_6+\beta_7=0$\\
(-1,-1,-2,-3,-2,-2,-1,-1)&$2\alpha_9-\beta_1+\beta_4+\beta_7+\beta_8=0$\\
(-1,-1,-2,-3,-3,-2,-1,-1)&$2\alpha_9-\beta_1+\beta_2-\beta_3+\beta_8=0$\\
(-1,-2,-2,-4,-3,-2,-1,0)&$2\alpha_9+\beta_1+\beta_2+\beta_4+\beta_6=0$\\
(-1,-2,-2,-3,-3,-2,-1,-1)&$2\alpha_9-\beta_3-\beta_4+\beta_6+\beta_8=0$\\
(-2,-2,-3,-4,-3,-2,-1,0)&$2\alpha_9+\beta_2+\beta_3+\beta_5+\beta_6=0$\\
(-2,-2,-3,-4,-3,-2,-1,-1)&$2\alpha_9-\beta_1+\beta_5+\beta_6+\beta_8=0$\\
(-1,-2,-2,-4,-3,-3,-2,-1)&$2\alpha_9+\beta_1+\beta_4+\beta_7+\beta_8=0$\\
(-2,-2,-3,-4,-3,-2,-2,-1)&$2\alpha_9+\beta_3+\beta_5-\beta_7+\beta_8=0$\\
(-2,-2,-3,-4,-3,-3,-2,-1)&$2\alpha_9+\beta_3+\beta_5+\beta_7+\beta_8=0$\\
(-2,-2,-3,-4,-4,-3,-2,-1)&$2\alpha_9+\beta_2-\beta_4+\beta_5+\beta_8=0$\\
(-2,-2,-3,-5,-4,-3,-2,-1)&$2\alpha_9+\beta_2+\beta_4+\beta_5+\beta_8=0$\\
(-2,-3,-3,-5,-4,-3,-2,-1)&$2\alpha_9+\beta_1+\beta_5+\beta_6+\beta_8=0$\\
(-2,-2,-4,-5,-4,-3,-2,-1)&$2\alpha_9-\beta_1+\beta_2+\beta_3+\beta_8=0$\\
(-2,-3,-4,-5,-4,-3,-2,-1)&$2\alpha_9+\beta_3-\beta_4+\beta_6+\beta_8=0$\\
(-2,-3,-4,-6,-4,-3,-2,-1)&$2\alpha_9+\beta_3+\beta_4+\beta_6+\beta_8=0$\\
(-2,-3,-4,-6,-5,-3,-2,-1)&$2\alpha_9+\beta_2+\beta_6-\beta_7+\beta_8=0$\\
(-2,-3,-4,-6,-5,-4,-2,-1)&$2\alpha_9+\beta_2+\beta_6+\beta_7+\beta_8=0$\\
(-2,-3,-4,-6,-5,-4,-3,-1)&$2\alpha_9+\beta_1+\beta_2+\beta_3+\beta_8=0$\\\hline
\end{tabular}
\end{center}
Let $\Delta_{9-i}$ be the root systems generated by all $\alpha_j$ except $j=i$. Suppose that for some index $i$, the coefficient in front of $\alpha_i$ is odd. Then, choosing a height function for which $\alpha_i$ is $\Delta_{9-i}$-lowest, we see that we land in one of the already treated cases (that of $\Delta_{9-i}=A_8$, $2A_4$, $A_2+E_6$ $4A_2$ or $A_5+A_1+A_2$ with $a_i=3$) and we have a linear combination with a coefficient $\pm 1$ for some $\alpha_k$. On the other hand, the assumption that all coefficients in front of $\alpha_i$ are even for all $i$ immediately contradicts the minimality of the starting relation \refeq{eqleCanChoose1_first}.  
\item  $\Delta$ is of type $F_4$. All possibilities for $\Delta_1\neq\Delta$ come in seven conjugacy classes (\cite[Table 10]{Dynkin:Semisimple}).
\begin{center}
\begin{tabular}{|cc|}\hline
\multicolumn{2}{|c|}{Case $A_2+A_2$}\\
\begin{picture}(0,25)
\put(1,5){\circle{1}}
\put(1,5.5){\line(0,1){6}}
\put(3,5){\tiny{4}}
\put(1,12){\circle{1}}
\put(3,12){\tiny{3}}
\put(1,19){\circle*{1}}
\put(1,19.5){\line(0,1){6}}
\put(3,19){\tiny{2}}
\put(1,26){\circle*{1}}
\put(3,26){\tiny{1}}
\end{picture}
&
\raisebox{25pt}{$\begin{array}{cccccccc}
\beta_1=(1,& 0,& 0,& 0)\\
\beta_2=(0,& 1,& 0,& 0)\\
\beta_3=(0,& 0,& 0,& 1)\\
\beta_4=(2,& 4,& 3,& 1)\\
\end{array}$}
\\
All possible $\alpha_5$&Linear combination\\
(0,0,1,0)&$3\alpha_5+2\beta_1+4\beta_2+\beta_3-\beta_4=0$\\
(0,1,1,0)&$3\alpha_5+2\beta_1+\beta_2+\beta_3-\beta_4=0$\\
(0,2,1,0)&$3\alpha_5+2\beta_1-2\beta_2+\beta_3-\beta_4=0$\\
(0,0,-1,0)&$3\alpha_5-2\beta_1-4\beta_2-\beta_3+\beta_4=0$\\
(0,0,1,1)&$3\alpha_5+2\beta_1+4\beta_2-2\beta_3-\beta_4=0$\\
(1,1,1,0)&$3\alpha_5-\beta_1+\beta_2+\beta_3-\beta_4=0$\\
(0,-1,-1,0)&$3\alpha_5-2\beta_1-\beta_2-\beta_3+\beta_4=0$\\
(0,1,1,1)&$3\alpha_5+2\beta_1+\beta_2-2\beta_3-\beta_4=0$\\
(2,2,1,0)&$3\alpha_5-4\beta_1-2\beta_2+\beta_3-\beta_4=0$\\
(0,2,1,1)&$3\alpha_5+2\beta_1-2\beta_2-2\beta_3-\beta_4=0$\\
(0,-2,-1,0)&$3\alpha_5-2\beta_1+2\beta_2-\beta_3+\beta_4=0$\\
(0,0,-1,-1)&$3\alpha_5-2\beta_1-4\beta_2+2\beta_3+\beta_4=0$\\
(-1,-1,-1,0)&$3\alpha_5+\beta_1-\beta_2-\beta_3+\beta_4=0$\\
(1,2,1,0)&$3\alpha_5-\beta_1-2\beta_2+\beta_3-\beta_4=0$\\
(1,1,1,1)&$3\alpha_5-\beta_1+\beta_2-2\beta_3-\beta_4=0$\\
(0,-1,-1,-1)&$3\alpha_5-2\beta_1-\beta_2+2\beta_3+\beta_4=0$\\
(2,2,1,1)&$3\alpha_5-4\beta_1-2\beta_2-2\beta_3-\beta_4=0$\\
(0,2,2,1)&$3\alpha_5+4\beta_1+2\beta_2-\beta_3-2\beta_4=0$\\
(-2,-2,-1,0)&$3\alpha_5+4\beta_1+2\beta_2-\beta_3+\beta_4=0$\\
(0,-2,-1,-1)&$3\alpha_5-2\beta_1+2\beta_2+2\beta_3+\beta_4=0$\\
(-1,-2,-1,0)&$3\alpha_5+\beta_1+2\beta_2-\beta_3+\beta_4=0$\\
(-1,-1,-1,-1)&$3\alpha_5+\beta_1-\beta_2+2\beta_3+\beta_4=0$\\
(1,2,1,1)&$3\alpha_5-\beta_1-2\beta_2-2\beta_3-\beta_4=0$\\
(2,2,2,1)&$3\alpha_5-2\beta_1+2\beta_2-\beta_3-2\beta_4=0$\\
(-2,-2,-1,-1)&$3\alpha_5+4\beta_1+2\beta_2+2\beta_3+\beta_4=0$\\
(0,-2,-2,-1)&$3\alpha_5-4\beta_1-2\beta_2+\beta_3+2\beta_4=0$\\
(-1,-2,-1,-1)&$3\alpha_5+\beta_1+2\beta_2+2\beta_3+\beta_4=0$\\
(1,2,2,1)&$3\alpha_5+\beta_1+2\beta_2-\beta_3-2\beta_4=0$\\
(2,4,2,1)&$3\alpha_5-2\beta_1-4\beta_2-\beta_3-2\beta_4=0$\\
(-2,-2,-2,-1)&$3\alpha_5+2\beta_1-2\beta_2+\beta_3+2\beta_4=0$\\
(-1,-2,-2,-1)&$3\alpha_5-\beta_1-2\beta_2+\beta_3+2\beta_4=0$\\
(1,3,2,1)&$3\alpha_5+\beta_1-\beta_2-\beta_3-2\beta_4=0$\\
(-2,-4,-2,-1)&$3\alpha_5+2\beta_1+4\beta_2+\beta_3+2\beta_4=0$\\
(-1,-3,-2,-1)&$3\alpha_5-\beta_1+\beta_2+\beta_3+2\beta_4=0$\\
(2,3,2,1)&$3\alpha_5-2\beta_1-\beta_2-\beta_3-2\beta_4=0$\\
(-2,-3,-2,-1)&$3\alpha_5+2\beta_1+\beta_2+\beta_3+2\beta_4=0$\\
\multicolumn{2}{|c|}{Number of relations: 36 }\\\hline
\end{tabular}
\end{center}
By direct observation we see that the mini-lemma from case $E_6$ applies to $\beta_3$ and $\beta_4$ for each of 36 possible values for $\alpha_5$ and we get that there exists some $\alpha_i$ with coefficient $a_i=\pm 1$. 
\begin{center}
\begin{tabular}{|cc|}\hline
\multicolumn{2}{|c|}{Case $B_4$}\\
\begin{picture}(0,25)
\put(1,5){\circle{1}}
\put(1,5.5){\line(0,1){6}}
\put(3,5){\tiny{4}}
\put(1,12){\circle{1}}
\put(1,12.5){\line(0,1){6}}
\put(3,12){\tiny{3}}
\put(1,19){\circle{1}}
\put(0.65,19.35){\line(0,1){6.3}}
\put(1.35,19.35){\line(0,1){6.3}}
\put(3,19){\tiny{2}}
\put(1,26){\circle*{1}}
\put(3,26){\tiny{1}}
\end{picture}
&
\raisebox{25pt}{$\begin{array}{cccccccc}
\beta_1=(0,& 1,& 0,& 0)\\
\beta_2=(0,& 0,& 1,& 0)\\
\beta_3=(0,& 0,& 0,& 1)\\
\beta_4=(2,& 2,& 1,& 0)\\
\end{array}$}
\\
Possibilities for lowest $\alpha_5$&Linear combination\\
(-1,-3,-2,-1)&$2\alpha_5+4\beta_1+3\beta_2+2\beta_3+\beta_4=0$\\
\hline
\end{tabular}
\end{center}
\begin{center}
\begin{tabular}{|cc|}\hline
\multicolumn{2}{|c|}{Case $A_1+A_3$}\\
\begin{picture}(0,25)
\put(1,5){\circle{1}}
\put(1,5.5){\line(0,1){6}}
\put(3,5){\tiny{4}}
\put(1,12){\circle{1}}
\put(1,12.5){\line(0,1){6}}
\put(3,12){\tiny{3}}
\put(1,19){\circle{1}}
\put(3,19){\tiny{2}}
\put(1,26){\circle*{1}}
\put(3,26){\tiny{1}}
\end{picture}
&
\raisebox{25pt}{$\begin{array}{cccccccc}
\beta_1=(1,& 0,& 0,& 0)\\
\beta_2=(0,& 0,& 1,& 0)\\
\beta_3=(0,& 0,& 0,& 1)\\
\beta_4=(2,& 4,& 2,& 1)\\
\end{array}$}
\\
$\Delta_1$-lowest $\alpha_5$&Linear combination\\
(-1,-1,-1,-1)&$4\alpha_5+2\beta_1+2\beta_2+3\beta_3+\beta_4=0$\\
(-2,-2,-2,-1)&$2\alpha_5+2\beta_1+2\beta_2+\beta_3+\beta_4=0$\\
(-1,-3,-2,-1)&$4\alpha_5-2\beta_1+2\beta_2+\beta_3+3\beta_4=0$\\
(-2,-3,-2,-1)&$4\alpha_5+2\beta_1+2\beta_2+\beta_3+3\beta_4=0$\\
\hline
\end{tabular}
\end{center}
\begin{center}
\begin{tabular}{|cc|}\hline
\multicolumn{2}{|c|}{Case $C_3+A_1$}\\
\begin{picture}(0,25)
\put(1,5){\circle{1}}
\put(3,5){\tiny{4}}
\put(1,12){\circle{1}}
\put(0.65,12.35){\line(0,1){6.3}}
\put(1.35,12.35){\line(0,1){6.3}}
\put(3,12){\tiny{3}}
\put(1,19){\circle*{1}}
\put(1,19.5){\line(0,1){6}}
\put(3,19){\tiny{2}}
\put(1,26){\circle*{1}}
\put(3,26){\tiny{1}}
\end{picture}
&
\raisebox{25pt}{$\begin{array}{cccccccc}
\beta_1=(1,& 0,& 0,& 0)\\
\beta_2=(0,& 1,& 0,& 0)\\
\beta_3=(0,& 0,& 1,& 0)\\
\beta_4=(2,& 4,& 3,& 2)\\
\end{array}$}
\\
$\Delta_1$-lowest $\alpha_5$&Linear combination\\
(-1,-1,-1,-1)&$2\alpha_5-2\beta_2-\beta_3+\beta_4=0$\\
(-2,-2,-1,-1)&$2\alpha_5+2\beta_1-\beta_3+\beta_4=0$\\
(-2,-2,-2,-1)&$2\alpha_5+2\beta_1+\beta_3+\beta_4=0$\\
(-2,-4,-2,-1)&$2\alpha_5+2\beta_1+4\beta_2+\beta_3+\beta_4=0$\\
(-1,-3,-2,-1)&$2\alpha_5+2\beta_2+\beta_3+\beta_4=0$\\
(-2,-4,-3,-1)&$2\alpha_5+2\beta_1+4\beta_2+3\beta_3+\beta_4=0$\\
(-2,-3,-2,-1)&$2\alpha_5+2\beta_1+2\beta_2+\beta_3+\beta_4=0$\\
\hline
\end{tabular}
\end{center}
\begin{center}
\begin{tabular}{|cc|}\hline
\multicolumn{2}{|c|}{Case $D_4$}\\
\begin{picture}(0,25)
\put(1,5){\circle{1}}
\put(3,5){\tiny{4}}
\put(1,5.5){\line(0,1){6}}
\put(1,12){\circle{1}}
\put(3,12){\tiny{3}}
\put(1.35,12.35){\line(1,1){6}}
\put(0.65,12.35){\line(-1,1){6}}
\put(8,19){\circle{1}}
\put(10,19){\tiny{2}}
\put(-6,19){\circle{1}}
\put(-4,19){\tiny{1}}
\end{picture}
&
\raisebox{25pt}{$\begin{array}{cccccccc}
\beta_1=(0,& 0,& 1,& 0)\\
\beta_2=(0,& 0,& 0,& 1)\\
\beta_3=(0,& 2,& 1,& 0)\\
\beta_4=(2,& 2,& 1,& 0)\\
\end{array}$}
\\
$\Delta_1$-lowest $\alpha_5$&Linear combination\\
(-1,-1,-1,-1)&$2\alpha_5+\beta_1+2\beta_2+\beta_4=0$\\
(-1,-2,-2,-1)&$2\alpha_5+2\beta_1+2\beta_2+\beta_3+\beta_4=0$\\
(-1,-3,-2,-1)&$2\alpha_5+\beta_1+2\beta_2+2\beta_3+\beta_4=0$\\
(-2,-3,-2,-1)&$2\alpha_5+\beta_1+2\beta_2+\beta_3+2\beta_4=0$\\
\hline
\end{tabular}
\end{center}

\begin{center}
\begin{tabular}{|cc|}\hline
\multicolumn{2}{|c|}{Case $B_2+2A_1$}\\
\begin{picture}(0,25)
\put(1,5){\circle{1}}
\put(3,5){\tiny{4}}
\put(1,12){\circle{1}}
\put(3,12){\tiny{3}}
\put(1,19){\circle{1}}
\put(0.65,19.35){\line(0,1){6.3}}
\put(1.35,19.35){\line(0,1){6.3}}
\put(3,19){\tiny{2}}
\put(1,26){\circle*{1}}
\put(3,26){\tiny{1}}
\end{picture}
&
\raisebox{25pt}{$\begin{array}{cccccccc}
\beta_1=(1,& 0,& 0,& 0)\\
\beta_2=(0,& 2,& 1,& 0)\\
\beta_3=(0,& 0,& 1,& 0)\\
\beta_4=(2,& 4,& 3,& 2)\\
\end{array}$}
\\
$\Delta_1$-lowest $\alpha_5$&Linear combination\\
(-1,-1,-1,0)&$2\alpha_5+2\beta_1+\beta_2+\beta_3=0$\\
(-1,-1,-1,-1)&$2\alpha_5-\beta_2+\beta_4=0$\\
(-2,-2,-1,-1)&$2\alpha_5+2\beta_1-\beta_3+\beta_4=0$\\
(-2,-2,-2,-1)&$2\alpha_5+2\beta_1+\beta_3+\beta_4=0$\\
(-1,-2,-2,-1)&$2\alpha_5+\beta_3+\beta_4=0$\\
(-2,-4,-2,-1)&$2\alpha_5+2\beta_1+2\beta_2-\beta_3+\beta_4=0$\\
(-1,-3,-2,-1)&$2\alpha_5+\beta_2+\beta_4=0$\\
(-2,-4,-3,-1)&$2\alpha_5+2\beta_1+2\beta_2+\beta_3+\beta_4=0$\\
(-2,-3,-2,-1)&$2\alpha_5+2\beta_1+\beta_2+\beta_4=0$\\
\hline
\end{tabular}
\end{center}

\begin{center}
\begin{tabular}{|cc|}\hline
\multicolumn{2}{|c|}{Case $4A_1$}\\
\begin{picture}(0,25)
\put(1,5){\circle{1}}
\put(3,5){\tiny{4}}
\put(1,12){\circle{1}}
\put(3,12){\tiny{3}}
\put(1,19){\circle{1}}
\put(3,19){\tiny{2}}
\put(1,26){\circle{1}}
\put(3,26){\tiny{1}}
\end{picture}
&
\raisebox{25pt}{$\begin{array}{cccccccc}
\beta_1=(2,& 2,& 1,& 0)\\
\beta_2=(0,& 2,& 1,& 0)\\
\beta_3=(0,& 0,& 1,& 0)\\
\beta_4=(2,& 4,& 3,& 2)\\
\end{array}$}
\\
$\Delta_1$-lowest $\alpha_5$&Linear combination\\
(0,-1,-1,0)&$2\alpha_5+\beta_2+\beta_3=0$\\
(-1,-1,-1,0)&$2\alpha_5+\beta_1+\beta_3=0$\\
(-1,-2,-1,0)&$2\alpha_5+\beta_1+\beta_2=0$\\
(-1,-1,-1,-1)&$2\alpha_5-\beta_2+\beta_4=0$\\
(-2,-2,-1,-1)&$2\alpha_5+\beta_1-\beta_2-\beta_3+\beta_4=0$\\
(-2,-2,-2,-1)&$2\alpha_5+\beta_1-\beta_2+\beta_3+\beta_4=0$\\
(-1,-2,-2,-1)&$2\alpha_5+\beta_3+\beta_4=0$\\
(-2,-4,-2,-1)&$2\alpha_5+\beta_1+\beta_2-\beta_3+\beta_4=0$\\
(-1,-3,-2,-1)&$2\alpha_5+\beta_2+\beta_4=0$\\
(-2,-4,-3,-1)&$2\alpha_5+\beta_1+\beta_2+\beta_3+\beta_4=0$\\
(-2,-3,-2,-1)&$2\alpha_5+\beta_1+\beta_4=0$\\
\hline
\end{tabular}
\end{center}
If the coefficient $a_5$ of $\alpha_5$ is odd then $\Delta_1=A_2+A_2$ or $F_4$, and we get the desired relation with coefficient $a_i=\pm 1$ for $i\neq 5$ in the former case, and $a_5=1$ in the latter. On the other hand, the assumption that all coefficients in front of $\alpha_i$ are even for all $i$ immediately contradicts the minimality of the starting relation \refeq{eqleCanChoose1_first}.
\item $\Delta$ is of type $G_2$. If the three roots $\alpha_i$ are of equal length and positive then one of them is the sum of the other two. If two of the $\alpha_i$ are long, say $\alpha_1$ and $\alpha_2$, and the third, $\alpha_3$, is short, then by direct check, after eventually switching $\alpha_1$ and $\alpha_2$, we get that $3\alpha_3=\pm 2\alpha_1\pm\alpha_2$. If two of the $\alpha_i$ are short then the third root is their integral linear combination. The general case for arbitrary positivity of the roots $\alpha_i$'s follows directly. 
\end{itemize}
\end{proof}

\section{On the actual implementation}\label{secImplementation}
\subsection{Modifications to \refsec{secTheTheoreticalAlgorithm} used in the implementation} \label{secTheModifications}

In this section $\Delta^+$ will stand for the root system of a classical Lie algebra and will be extended to a set $\bar\Delta$. In addition to the data given for the algorithm in \refsec{secTheTheoreticalAlgorithm}, we will give a total order on $\bar\Delta$ for which $\alpha_1\succ \dots \succ \alpha_N$. In the context of computing the vector partition function, the total orders $\succ$ were defined in \cite{BBCV}.  

We are using the notation from \refsec{secRootSystems}. For A$_n$ define $\bar\Delta$ to be the root system of A$_n$, define $\prec$ to be the order induced by $\eta_i-\eta_{i+1}\succ\eta_{i+1}-\eta_{i+2}$ and by $\alpha+\beta\succ \alpha,\beta$. Equivalently, define $\succ$ to be the order naturally induced from the lexicographic order on weights $\eta_i$ of $\gl(n)$, projected down to $\sL(n)$. For systems D$_n$, respectively B$_n$, define $\bar\Delta$ to be the root system of B$_n$, respectively D$_n$, $\cup\{2\eta_i\}$ and define $\prec$ in the same way as for A$_n$. For system C$_n$ define $\bar\Delta$ to be the root system of C$_n$ and define the order in the same way as for A$_n$.
\subsubsection{An algorithm for the partial fraction decomposition of the generating function for the vector partition function of the classical root systems }\label{secTheActualAlgorithm}

\begin{itemize}
\item[Step 1] (Reduction step). We are given a fraction $p(x)\frac{1}{(1-x^{\alpha_{i_1}})^{l_1}} \dots\frac{ 1}{( 1- x^{ \alpha_{i_k} })^{l_k}}$, $i_1<\dots<i_k$.
\begin{itemize}
\item[Step 1.1] Find the smallest index $a$, and the smallest index $b$ depending on $a$, and the smaller of the  numbers $n=1,2$, for which $\alpha_{i_a}- n\alpha_{i_b}\in \bar\Delta$ and $\alpha_{i_a}- n\alpha_{i_b}\prec \alpha_{i_b}$. If no such couple of indices $(i_a,i_b)$ exists label the fraction as ``reduced''. 
\item[Step 1.2] For the couple of indices $(i_a,i_b)$ perform the reduction formula \refle{leFlas}(e). Label any newly created fractions as ``non-reduced''.
\end{itemize}
\item[Step 2] (Termination step). Find the first fraction labeled ``non-reduced'' in the memory and go to Step 1. If all fractions are labeled as reduced proceed to Step 3.
\item[Step 3] (Getting rid of the redundant short roots in the denominators). Transform the output of Step 2 to a partial fraction decomposition (see \refdef{defPartFracDecomposition}). One needs only to apply the formula  $\frac{1}{1-x^{\alpha}}= \frac{1 +x^\alpha} {1-x^{ 2\alpha }}$ where it is needed (see \refle{le1} below).
\end{itemize}

We leave the following lemma without proof\footnote{The ``vector partition function'' program is designed to display an error message and crash if it doesn't reach a partial fraction decomposition.}. 
\begin{lemma}\label{le1}
Let $\Delta^+\subset\bar\Delta$ be the set of positive roots of a classical root systems of rank $n$. Steps 1 and 2 of the algorithm in \refsec{secTheActualAlgorithm} produce a decomposition of the generating function of the partition function of any subset $I\subset\Delta^+$ of the form
\begin{eqnarray*}
\prod_{\alpha\in I}\frac{1}{1-x^{\alpha}}&=& \sum_{\substack{\alpha_{i_j}\in\Delta\\ l_j,m_j\in\ZZ_{\geq 0}, \sum_j l_j+\sum_j m_j= |I|\\ 2\alpha_{i_j}\notin \bar\Delta\Rightarrow m_j=0}} 
p_{l_1,\dots,l_n, m_1,\dots,m_n,\alpha_{i_1},\dots,\alpha_{i_n}}(x)  \\&& \frac{1} {(1-x^{\alpha_1})^{l_1} (1- x^{ 2\alpha_1} )^{ m_1}} \dots \frac{1}{(1-x^{\alpha_{i_n}})^{l_n}(1-x^{2\alpha_{i_n}})^{m_n}},
\end{eqnarray*}
where $n$ is exactly the rank of the root system and $p_{l_1,\dots,l_n,m_1,\dots,m_n,\alpha_{i_1},\dots,\alpha_{i_n}}$ are polynomials in $x_1,\dots, x_n$ and $x_1^{-1},\dots x_n^{-1} $. 
\end{lemma}
\subsection{Short outline of program functionality}\label{secProgramFunctionality}
The ``vector partition function'' program is being written in C++. It uses template classes for polynomials over arbitrary commutative rings with identity, a realization of the ring of quasi-numbers (multiplicative relations realized, additive not), classes for operations with large rational numbers (size limit dictated only by system resources), and hash-table array management template classes. Additionally, there are classes for polyhedral cone divisions of Euclidean space, classes for computing orbits of Weyl groups and Kazhdan-Lusztig polynomials and coefficients. The entire code is written from scratch by the author; no external packages are used. The standard C++ library for stream/string manipulations is used only for the final printouts of the algorithm. 

Before and after each partial fraction decomposition runs, the ``vector partition function'' program is set to substitute rational numbers in the expressions it transforms. Comparing the two numbers has been our main tool for catching programming mistakes. Each step of the program has successfully passed such ``checksum'' tests. 

\nocite{LatE} \nocite{SphericalExplorer}
\newcommand{\etalchar}[1]{$^{#1}$}
\providecommand{\bysame}{\leavevmode\hbox to3em{\hrulefill}\thinspace}
\providecommand{\MR}{\relax\ifhmode\unskip\space\fi MR }
\providecommand{\MRhref}[2]{%
  \href{http://www.ams.org/mathscinet-getitem?mr=#1}{#2}
}
\providecommand{\href}[2]{#2}

\appendix
\section{Appendix. Complex exponents and quasi-numbers}
Let $x\in\ZZ$. Then
\begin{eqnarray*}
Q_{l,D}(x):=\doublebrace{1}{\mathrm{if~}x\equiv l (\mod D)}{0}{\mathrm{otherwise}}&=& \frac{e^{2\pi i \frac{l(1-n)}{D}}}{D} \prod_{\substack{0\leq k\leq D\\k\neq l}}\left(e^{2\pi i\frac{x}{D}}-e^{2\pi i\frac{k}{D}}\right).
\end{eqnarray*}
\newpage
\section{Appendix. Partial fraction decomposition tables for $A_2$, $A_3$, $A_4$, $B_2$, $B_3$, $C_2$, $C_3$ and $G_2$}
\addtolength{\hoffset}{-3.5cm}
\addtolength{\textwidth}{6.8cm}
\addtolength{\voffset}{-3.3cm}
\addtolength{\textheight}{6.3cm}
\begin{eqnarray*}
\prod_{\alpha\in G_2}\frac{1}{1-x^\alpha} \\ 
\end{eqnarray*}
=
\begin{eqnarray*}
&& \frac{1}{(1-x_1)} \frac{1}{(1-x_2)} \frac{1}{(1-x_1x_2^{3})} \frac{1}{(1-x_1x_2)} \frac{1}{(1-x_1^{2}x_2^{3})} \frac{1}{(1-x_1x_2^{2})} \\ 
\end{eqnarray*}
=
\begin{eqnarray*}
&&+(x_1x_2^{-4}+x_1x_2^{-5}+x_1x_2^{-6}+4x_1^{2}x_2^{-4}+4x_1^{2}x_2^{-5}+4x_1^{2}x_2^{-6}+6x_1^{3}x_2^{-4}+6x_1^{3}x_2^{-5}+6x_1^{3}x_2^{-6}\\
&& 
+4x_1^{4}x_2^{-4}+4x_1^{4}x_2^{-5}+4x_1^{4}x_2^{-6}+x_1^{5}x_2^{-4}+x_1^{5}x_2^{-5}+x_1^{5}x_2^{-6})\\
&& \frac{1}{(1-x_1^{2})^5} \frac{1}{(1-x_1x_2)} \\ 
&&+(13x_2^{-8}+15x_2^{-9}+x_1^{-3}x_2^{-10}+x_1^{-3}x_2^{-11}+2x_1^{-3}x_2^{-12}+2x_1^{-2}x_2^{-9}+5x_1^{-2}x_2^{-10}+6x_1^{-2}x_2^{-11}+6x_1^{-1}x_2^{-8}\\
&& 
+9x_1^{-1}x_2^{-9}+11x_1^{-1}x_2^{-10}+x_1^{-4}x_2^{-13}+x_1^{-4}x_2^{-14}+x_1^{-4}x_2^{-15}+5x_1^{-3}x_2^{-13}+5x_1^{-3}x_2^{-14}+9x_1^{-2}x_2^{-12}\\
&& 
+11x_1^{-2}x_2^{-13}+5x_1^{-1}x_2^{-7}+2x_1^{-1}x_2^{-6}+x_1^{-2}x_2^{-7}+x_1^{-2}x_2^{-8}+x_1^{-1}x_2^{-4}+x_1^{-1}x_2^{-5}+x_2^{-3}+5x_2^{-4}\\
&& 
+5x_2^{-5}+2x_1x_2^{-2}+6x_1x_2^{-3}+11x_1x_2^{-4}+x_1x_2^{-1}+x_1^{2}+5x_1^{2}x_2^{-1}+9x_1^{2}x_2^{-2}+x_1^{3}x_2+5x_1^{3}+11x_1^{3}x_2^{-1}\\
&& 
+12x_1x_2^{-8}+10x_1x_2^{-9}+14x_1^{-1}x_2^{-11}+14x_2^{-10}+4x_1^{-3}x_2^{-15}+10x_1^{-2}x_2^{-14}+16x_1^{-1}x_2^{-12}+14x_1^{-1}x_2^{-13}\\
&& 
+10x_2^{-7}+8x_2^{-6}+10x_1x_2^{-5}+14x_1^{2}x_2^{-3}+14x_1^{2}x_2^{-4}+16x_1^{3}x_2^{-2}+4x_1^{4}x_2+10x_1^{4}+14x_1^{4}x_2^{-1}+3x_1^{2}x_2^{-8}\\
&& 
-x_1^{4}x_2^{-9}+16x_2^{-11}+11x_1x_2^{-10}+6x_1^{-2}x_2^{-15}+10x_1^{-1}x_2^{-14}+14x_2^{-12}+11x_2^{-13}+10x_1x_2^{-7}+12x_1x_2^{-6}\\
&& 
+10x_1^{2}x_2^{-5}+16x_1^{3}x_2^{-3}+11x_1^{3}x_2^{-4}+14x_1^{4}x_2^{-2}+6x_1^{5}x_2+10x_1^{5}+11x_1^{5}x_2^{-1}-2x_1^{3}x_2^{-8}-3x_1^{3}x_2^{-9}\\
&& 
+9x_1x_2^{-11}+5x_1^{2}x_2^{-10}+4x_1^{-1}x_2^{-15}+5x_2^{-14}+6x_1x_2^{-12}+5x_1x_2^{-13}+5x_1^{2}x_2^{-7}+8x_1^{2}x_2^{-6}+5x_1^{3}x_2^{-5}\\
&& 
+9x_1^{4}x_2^{-3}+5x_1^{4}x_2^{-4}+6x_1^{5}x_2^{-2}+4x_1^{6}x_2+5x_1^{6}+5x_1^{6}x_2^{-1}-x_1^{4}x_2^{-8}+2x_1^{2}x_2^{-11}+x_1^{3}x_2^{-10}\\
&& 
+x_2^{-15}+x_1x_2^{-14}+x_1^{2}x_2^{-12}+x_1^{2}x_2^{-13}+x_1^{3}x_2^{-7}+2x_1^{3}x_2^{-6}+x_1^{4}x_2^{-5}+2x_1^{5}x_2^{-3}+x_1^{5}x_2^{-4}\\
&& 
+x_1^{6}x_2^{-2}+x_1^{7}x_2+x_1^{7}+x_1^{7}x_2^{-1})\\
&& \frac{1}{(1-x_1^{2})^5} \frac{1}{(1-x_1x_2^{2})} \\ 
\end{eqnarray*}\begin{eqnarray*}
&&-x_1^{2}x_2^{-4} \frac{1}{(1-x_1)^5} \frac{1}{(1-x_1x_2)} \\ 
&&+(-24x_1^{2}x_2^{-4}-15x_1^{2}x_2^{-3}-8x_1^{2}x_2^{-2}-3x_1x_2^{-3}-8x_1x_2^{-4}-11x_1x_2^{-5}-13x_1x_2^{-6}-2x_2^{-5}-3x_2^{-6}\\
&& 
-6x_2^{-7}-5x_2^{-8}-x_1^{-1}x_2^{-7}-x_1^{-1}x_2^{-8}-x_1^{-1}x_2^{-9}-x_1^{2}-2x_1^{2}x_2^{-1}-x_1^{3}x_2-6x_1^{3}-10x_1^{3}x_2^{-1}\\
&& 
-x_1^{4}x_2^{2}-5x_1^{4}x_2-15x_1^{4}-x_1x_2^{-2}-23x_1^{3}x_2^{-2}-x_2^{-4}-36x_1^{3}x_2^{-4}-30x_1^{3}x_2^{-3}-24x_1^{2}x_2^{-5}\\
&& 
-22x_1^{2}x_2^{-6}-14x_1x_2^{-7}-10x_1x_2^{-8}-4x_2^{-9}-20x_1^{4}x_2^{-1}-4x_1^{5}x_2^{2}-10x_1^{5}x_2-20x_1^{5}-32x_1^{4}x_2^{-2}\\
&& 
-29x_1^{4}x_2^{-4}-30x_1^{4}x_2^{-3}-26x_1^{3}x_2^{-5}-18x_1^{3}x_2^{-6}-16x_1^{2}x_2^{-7}-10x_1^{2}x_2^{-8}-6x_1x_2^{-9}-20x_1^{5}x_2^{-1}\\
&& 
-6x_1^{6}x_2^{2}-10x_1^{6}x_2-15x_1^{6}-23x_1^{5}x_2^{-2}-12x_1^{5}x_2^{-4}-15x_1^{5}x_2^{-3}-14x_1^{4}x_2^{-5}-7x_1^{4}x_2^{-6}\\
&& 
-9x_1^{3}x_2^{-7}-5x_1^{3}x_2^{-8}-4x_1^{2}x_2^{-9}-10x_1^{6}x_2^{-1}-4x_1^{7}x_2^{2}-5x_1^{7}x_2-6x_1^{7}-8x_1^{6}x_2^{-2}-2x_1^{6}x_2^{-4}\\
&& 
-3x_1^{6}x_2^{-3}-3x_1^{5}x_2^{-5}-x_1^{5}x_2^{-6}-2x_1^{4}x_2^{-7}-x_1^{4}x_2^{-8}-x_1^{3}x_2^{-9}-2x_1^{7}x_2^{-1}-x_1^{8}x_2^{2}\\
&& 
-x_1^{8}x_2-x_1^{8}-x_1^{7}x_2^{-2})\\
&& \frac{1}{(1-x_1^{2})^5} \frac{1}{(1-x_1^{2}x_2^{3})} \\ 
&&+(-x_1^{-4}x_2^{-13}-x_1^{-4}x_2^{-14}-x_1^{-4}x_2^{-15}-2x_1^{-3}x_2^{-12}-6x_1^{-3}x_2^{-13}-5x_1^{-3}x_2^{-14}-7x_1^{-2}x_2^{-11}\\
&& 
-10x_1^{-2}x_2^{-12}-15x_1^{-2}x_2^{-13}-x_1^{-3}x_2^{-11}-2x_1^{-2}x_2^{-10}-6x_1^{-1}x_2^{-9}-10x_1^{-1}x_2^{-10}-19x_1^{-1}x_2^{-11}\\
&& 
-x_1^{-2}x_2^{-9}-x_1^{-1}x_2^{-8}-x_2^{-7}-5x_2^{-8}-15x_2^{-9}-4x_1^{-3}x_2^{-15}-10x_1^{-2}x_2^{-14}-20x_1^{-1}x_2^{-12}-20x_1^{-1}x_2^{-13}\\
&& 
-20x_2^{-10}-26x_2^{-11}-4x_1x_2^{-7}-10x_1x_2^{-8}-20x_1x_2^{-9}-6x_1^{-2}x_2^{-15}-10x_1^{-1}x_2^{-14}-20x_2^{-12}-15x_2^{-13}\\
&& 
-20x_1x_2^{-10}-19x_1x_2^{-11}-6x_1^{2}x_2^{-7}-10x_1^{2}x_2^{-8}-15x_1^{2}x_2^{-9}-4x_1^{-1}x_2^{-15}-5x_2^{-14}-10x_1x_2^{-12}\\
&& 
-6x_1x_2^{-13}-10x_1^{2}x_2^{-10}-7x_1^{2}x_2^{-11}-4x_1^{3}x_2^{-7}-5x_1^{3}x_2^{-8}-6x_1^{3}x_2^{-9}-x_2^{-15}-x_1x_2^{-14}-2x_1^{2}x_2^{-12}\\
&& 
-x_1^{2}x_2^{-13}-2x_1^{3}x_2^{-10}-x_1^{3}x_2^{-11}-x_1^{4}x_2^{-7}-x_1^{4}x_2^{-8}-x_1^{4}x_2^{-9})\\
&& \frac{1}{(1-x_1^{2})^5} \frac{1}{(1-x_1x_2^{3})} \\ 
\end{eqnarray*}\begin{eqnarray*}
&&+(x_2^{-9}+4x_1x_2^{-9}+6x_1^{2}x_2^{-9}+4x_1^{3}x_2^{-9}+x_1^{4}x_2^{-9}) \frac{1}{(1-x_1^{2})^5} \frac{1}{(1-x_2)} \\ 
&&+(x_1^{3}+x_1^{4}x_2+x_1^{5}x_2^{2}+x_1^{2}x_2^{-2}+x_1^{3}x_2^{-1}+x_1^{4}+x_1x_2^{-4}+x_1^{2}x_2^{-3}+x_1^{3}x_2^{-2}) \frac{1}{(1-x_1)^5} \frac{1}{(1-x_1^{2}x_2^{3})} \\ 
&&+(x_1x_2^{-4}+x_1^{2}x_2^{-3}-x_2^{-4}-x_1x_2^{-3}-x_1^{2}x_2^{-2}-x_1^{2}x_2^{-1}-x_1^{3}-x_1^{4}x_2) \frac{1}{(1-x_1)^5} \frac{1}{(1-x_1x_2^{2})} \\ 
\end{eqnarray*}

\begin{eqnarray*}
\prod_{\alpha\in A_2}\frac{1}{1-x^\alpha}&=&-x_2 \frac{1}{(1-x_1x_2)} \frac{1}{(1-x_2)^2} \\ 
&&+ \frac{1}{(1-x_1)} \frac{1}{(1-x_2)^2} \\ 
\end{eqnarray*}
\begin{eqnarray*}
\prod_{\alpha\in A_3}\frac{1}{1-x^\alpha}&=&-x_2x_3^{3} \frac{1}{(1-x_1x_2x_3)} \frac{1}{(1-x_2x_3)} \frac{1}{(1-x_3)^4} \\ 
&&-x_2x_3^{2} \frac{1}{(1-x_1x_2)} \frac{1}{(1-x_2x_3)} \frac{1}{(1-x_3)^4} \\ 
&&+2x_3^{2} \frac{1}{(1-x_1)} \frac{1}{(1-x_2x_3)} \frac{1}{(1-x_3)^4} \\ 
&&+x_2x_3^{2} \frac{1}{(1-x_1x_2x_3)} \frac{1}{(1-x_2)} \frac{1}{(1-x_3)^4} \\ 
&&-x_2x_3^{3} \frac{1}{(1-x_1x_2x_3)} \frac{1}{(1-x_2x_3)^2} \frac{1}{(1-x_3)^3} \\ 
&&+ \frac{1}{(1-x_1)} \frac{1}{(1-x_2)^2} \frac{1}{(1-x_3)^3} \\ 
&&-x_2 \frac{1}{(1-x_1x_2)} \frac{1}{(1-x_2)^2} \frac{1}{(1-x_3)^3} \\ 
&&+x_2x_3 \frac{1}{(1-x_1x_2)} \frac{1}{(1-x_2)} \frac{1}{(1-x_3)^4} \\ 
&&-2x_3 \frac{1}{(1-x_1)} \frac{1}{(1-x_2)} \frac{1}{(1-x_3)^4} \\ 
&&+x_3^{2} \frac{1}{(1-x_1)} \frac{1}{(1-x_2x_3)^2} \frac{1}{(1-x_3)^3} \\ 
\end{eqnarray*}
\begin{eqnarray*}
\prod_{\alpha\in A_4}\frac{1}{1-x^\alpha}&=&+2x_2x_3^{3}x_4^{6} \frac{1}{(1-x_1x_2x_3x_4)} \frac{1}{(1-x_2x_3x_4)} \frac{1}{(1-x_3x_4)} \frac{1}{(1-x_4)^7} \\ 
&&-2x_2x_3^{3}x_4^{5} \frac{1}{(1-x_1x_2x_3)} \frac{1}{(1-x_2x_3x_4)} \frac{1}{(1-x_3x_4)} \frac{1}{(1-x_4)^7} \\ 
&&+4x_2x_3^{2}x_4^{5} \frac{1}{(1-x_1x_2)} \frac{1}{(1-x_2x_3x_4)} \frac{1}{(1-x_3x_4)} \frac{1}{(1-x_4)^7} \\ 
&&-4x_3^{2}x_4^{5} \frac{1}{(1-x_1)} \frac{1}{(1-x_2x_3x_4)} \frac{1}{(1-x_3x_4)} \frac{1}{(1-x_4)^7} \\ 
&&+2x_2x_3^{3}x_4^{5} \frac{1}{(1-x_1x_2x_3x_4)} \frac{1}{(1-x_2x_3)} \frac{1}{(1-x_3x_4)} \frac{1}{(1-x_4)^7} \\ 
&&+x_2x_3^{3}x_4^{6} \frac{1}{(1-x_1x_2x_3x_4)} \frac{1}{(1-x_2x_3x_4)} \frac{1}{(1-x_3x_4)^4} \frac{1}{(1-x_4)^4} \\ 
&&-4x_2x_3^{2}x_4^{5} \frac{1}{(1-x_1x_2x_3x_4)} \frac{1}{(1-x_2)} \frac{1}{(1-x_3x_4)} \frac{1}{(1-x_4)^7} \\ 
&&-2x_2x_3^{3}x_4^{5} \frac{1}{(1-x_1x_2x_3x_4)} \frac{1}{(1-x_2x_3x_4)} \frac{1}{(1-x_3)} \frac{1}{(1-x_4)^7} \\ 
&&+2x_2x_3^{3}x_4^{6} \frac{1}{(1-x_1x_2x_3x_4)} \frac{1}{(1-x_2x_3x_4)} \frac{1}{(1-x_3x_4)^3} \frac{1}{(1-x_4)^5} \\ 
&&+2x_2x_3^{3}x_4^{6} \frac{1}{(1-x_1x_2x_3x_4)} \frac{1}{(1-x_2x_3x_4)} \frac{1}{(1-x_3x_4)^2} \frac{1}{(1-x_4)^6} \\ 
&&+x_2x_3^{2}x_4^{3} \frac{1}{(1-x_1x_2)} \frac{1}{(1-x_2x_3x_4)} \frac{1}{(1-x_3)^2} \frac{1}{(1-x_4)^6} \\ 
&&+x_3^{2}x_4^{2} \frac{1}{(1-x_1)} \frac{1}{(1-x_2x_3)^2} \frac{1}{(1-x_3)} \frac{1}{(1-x_4)^6} \\ 
&&-2x_2x_3^{3}x_4^{4} \frac{1}{(1-x_1x_2x_3)} \frac{1}{(1-x_2x_3)} \frac{1}{(1-x_3x_4)} \frac{1}{(1-x_4)^7} \\ 
&&+4x_2x_3^{2}x_4^{4} \frac{1}{(1-x_1x_2x_3)} \frac{1}{(1-x_2)} \frac{1}{(1-x_3x_4)} \frac{1}{(1-x_4)^7} \\ 
&&+x_2x_3^{2}x_4^{4} \frac{1}{(1-x_1x_2x_3)} \frac{1}{(1-x_2)} \frac{1}{(1-x_3x_4)^2} \frac{1}{(1-x_4)^6} \\ 
&&-x_2x_3^{3}x_4^{3} \frac{1}{(1-x_1x_2x_3)} \frac{1}{(1-x_2x_3x_4)} \frac{1}{(1-x_3)^2} \frac{1}{(1-x_4)^6} \\ 
\end{eqnarray*}\begin{eqnarray*}
&&+2x_2x_3^{3}x_4^{4} \frac{1}{(1-x_1x_2x_3)} \frac{1}{(1-x_2x_3x_4)} \frac{1}{(1-x_3)} \frac{1}{(1-x_4)^7} \\ 
&&-x_2x_3^{3}x_4^{5} \frac{1}{(1-x_1x_2x_3)} \frac{1}{(1-x_2x_3x_4)} \frac{1}{(1-x_3x_4)^2} \frac{1}{(1-x_4)^6} \\ 
&&-2x_3^{2}x_4^{5} \frac{1}{(1-x_1)} \frac{1}{(1-x_2x_3x_4)} \frac{1}{(1-x_3x_4)^4} \frac{1}{(1-x_4)^4} \\ 
&&-4x_2x_3^{2}x_4^{4} \frac{1}{(1-x_1x_2)} \frac{1}{(1-x_2x_3)} \frac{1}{(1-x_3x_4)} \frac{1}{(1-x_4)^7} \\ 
&&+6x_2x_4^{3} \frac{1}{(1-x_1x_2)} \frac{1}{(1-x_2)^2} \frac{1}{(1-x_3x_4)} \frac{1}{(1-x_4)^6} \\ 
&&-x_3^{2}x_4 \frac{1}{(1-x_1)} \frac{1}{(1-x_2x_3)^2} \frac{1}{(1-x_3)^2} \frac{1}{(1-x_4)^5} \\ 
&&-4x_2x_3^{2}x_4^{4} \frac{1}{(1-x_1x_2)} \frac{1}{(1-x_2x_3x_4)} \frac{1}{(1-x_3)} \frac{1}{(1-x_4)^7} \\ 
&&+x_2x_3^{2}x_4^{5} \frac{1}{(1-x_1x_2)} \frac{1}{(1-x_2x_3x_4)} \frac{1}{(1-x_3x_4)^4} \frac{1}{(1-x_4)^4} \\ 
&&+2x_2x_3^{2}x_4^{5} \frac{1}{(1-x_1x_2)} \frac{1}{(1-x_2x_3x_4)} \frac{1}{(1-x_3x_4)^3} \frac{1}{(1-x_4)^5} \\ 
&&+3x_2x_3^{2}x_4^{5} \frac{1}{(1-x_1x_2)} \frac{1}{(1-x_2x_3x_4)} \frac{1}{(1-x_3x_4)^2} \frac{1}{(1-x_4)^6} \\ 
&&+4x_3^{2}x_4^{4} \frac{1}{(1-x_1)} \frac{1}{(1-x_2x_3)} \frac{1}{(1-x_3x_4)} \frac{1}{(1-x_4)^7} \\ 
&&-x_3^{2}x_4^{5} \frac{1}{(1-x_1)} \frac{1}{(1-x_2x_3x_4)^2} \frac{1}{(1-x_3x_4)} \frac{1}{(1-x_4)^6} \\ 
&&+x_3^{2} \frac{1}{(1-x_1)} \frac{1}{(1-x_2x_3)^2} \frac{1}{(1-x_3)^3} \frac{1}{(1-x_4)^4} \\ 
&&+4x_3^{2}x_4^{4} \frac{1}{(1-x_1)} \frac{1}{(1-x_2x_3x_4)} \frac{1}{(1-x_3)} \frac{1}{(1-x_4)^7} \\ 
&&-4x_3^{2}x_4^{5} \frac{1}{(1-x_1)} \frac{1}{(1-x_2x_3x_4)} \frac{1}{(1-x_3x_4)^3} \frac{1}{(1-x_4)^5} \\ 
&&-4x_3^{2}x_4^{5} \frac{1}{(1-x_1)} \frac{1}{(1-x_2x_3x_4)} \frac{1}{(1-x_3x_4)^2} \frac{1}{(1-x_4)^6} \\ 
\end{eqnarray*}\begin{eqnarray*}
&&-x_2x_3^{2}x_4^{3} \frac{1}{(1-x_1x_2x_3x_4)} \frac{1}{(1-x_2)} \frac{1}{(1-x_3)^2} \frac{1}{(1-x_4)^6} \\ 
&&+x_2x_3^{3}x_4^{3} \frac{1}{(1-x_1x_2x_3x_4)} \frac{1}{(1-x_2x_3)} \frac{1}{(1-x_3)^2} \frac{1}{(1-x_4)^6} \\ 
&&-2x_2x_3^{3}x_4^{4} \frac{1}{(1-x_1x_2x_3x_4)} \frac{1}{(1-x_2x_3)} \frac{1}{(1-x_3)} \frac{1}{(1-x_4)^7} \\ 
&&+x_2x_3^{3}x_4^{5} \frac{1}{(1-x_1x_2x_3x_4)} \frac{1}{(1-x_2x_3)} \frac{1}{(1-x_3x_4)^2} \frac{1}{(1-x_4)^6} \\ 
&&-x_2x_3^{2}x_4^{5} \frac{1}{(1-x_1x_2x_3x_4)} \frac{1}{(1-x_2)} \frac{1}{(1-x_3x_4)^4} \frac{1}{(1-x_4)^4} \\ 
&&+x_2x_3^{3}x_4^{6} \frac{1}{(1-x_1x_2x_3x_4)} \frac{1}{(1-x_2x_3x_4)^2} \frac{1}{(1-x_3x_4)} \frac{1}{(1-x_4)^6} \\ 
&&+4x_2x_3^{2}x_4^{4} \frac{1}{(1-x_1x_2x_3x_4)} \frac{1}{(1-x_2)} \frac{1}{(1-x_3)} \frac{1}{(1-x_4)^7} \\ 
&&-2x_2x_3^{2}x_4^{5} \frac{1}{(1-x_1x_2x_3x_4)} \frac{1}{(1-x_2)} \frac{1}{(1-x_3x_4)^3} \frac{1}{(1-x_4)^5} \\ 
&&-3x_2x_3^{2}x_4^{5} \frac{1}{(1-x_1x_2x_3x_4)} \frac{1}{(1-x_2)} \frac{1}{(1-x_3x_4)^2} \frac{1}{(1-x_4)^6} \\ 
&&-4x_3x_4^{2} \frac{1}{(1-x_1)} \frac{1}{(1-x_2)} \frac{1}{(1-x_3)^2} \frac{1}{(1-x_4)^6} \\ 
&&-x_2x_3^{2} \frac{1}{(1-x_1x_2)} \frac{1}{(1-x_2x_3)} \frac{1}{(1-x_3)^4} \frac{1}{(1-x_4)^4} \\ 
&&+3x_2x_4 \frac{1}{(1-x_1x_2)} \frac{1}{(1-x_2)^2} \frac{1}{(1-x_3)^2} \frac{1}{(1-x_4)^5} \\ 
&&+2x_2x_3x_4^{2} \frac{1}{(1-x_1x_2)} \frac{1}{(1-x_2)} \frac{1}{(1-x_3)^2} \frac{1}{(1-x_4)^6} \\ 
&&-4x_3^{2}x_4 \frac{1}{(1-x_1)} \frac{1}{(1-x_2x_3)} \frac{1}{(1-x_3)^3} \frac{1}{(1-x_4)^5} \\ 
&&+4x_3x_4 \frac{1}{(1-x_1)} \frac{1}{(1-x_2)} \frac{1}{(1-x_3)^3} \frac{1}{(1-x_4)^5} \\ 
&&+2x_3^{2} \frac{1}{(1-x_1)} \frac{1}{(1-x_2x_3)} \frac{1}{(1-x_3)^4} \frac{1}{(1-x_4)^4} \\ 
\end{eqnarray*}\begin{eqnarray*}
&&+x_2x_3^{2} \frac{1}{(1-x_1x_2x_3)} \frac{1}{(1-x_2)} \frac{1}{(1-x_3)^4} \frac{1}{(1-x_4)^4} \\ 
&&+x_2x_3^{3}x_4^{3} \frac{1}{(1-x_1x_2x_3)} \frac{1}{(1-x_2x_3)^2} \frac{1}{(1-x_3x_4)} \frac{1}{(1-x_4)^6} \\ 
&&-x_2x_3^{3} \frac{1}{(1-x_1x_2x_3)} \frac{1}{(1-x_2x_3)} \frac{1}{(1-x_3)^4} \frac{1}{(1-x_4)^4} \\ 
&&+2x_2x_3^{3}x_4 \frac{1}{(1-x_1x_2x_3)} \frac{1}{(1-x_2x_3)} \frac{1}{(1-x_3)^3} \frac{1}{(1-x_4)^5} \\ 
&&-2x_2x_3^{3}x_4^{2} \frac{1}{(1-x_1x_2x_3)} \frac{1}{(1-x_2x_3)} \frac{1}{(1-x_3)^2} \frac{1}{(1-x_4)^6} \\ 
&&+2x_2x_3^{3}x_4^{3} \frac{1}{(1-x_1x_2x_3)} \frac{1}{(1-x_2x_3)} \frac{1}{(1-x_3)} \frac{1}{(1-x_4)^7} \\ 
&&-2x_2x_3^{2}x_4 \frac{1}{(1-x_1x_2x_3)} \frac{1}{(1-x_2)} \frac{1}{(1-x_3)^3} \frac{1}{(1-x_4)^5} \\ 
&&+3x_2x_3^{2}x_4^{2} \frac{1}{(1-x_1x_2x_3)} \frac{1}{(1-x_2)} \frac{1}{(1-x_3)^2} \frac{1}{(1-x_4)^6} \\ 
&&-4x_2x_3^{2}x_4^{3} \frac{1}{(1-x_1x_2x_3)} \frac{1}{(1-x_2)} \frac{1}{(1-x_3)} \frac{1}{(1-x_4)^7} \\ 
&&-x_2x_3^{3}x_4^{2} \frac{1}{(1-x_1x_2x_3)} \frac{1}{(1-x_2x_3)^2} \frac{1}{(1-x_3)} \frac{1}{(1-x_4)^6} \\ 
&&-6x_4^{3} \frac{1}{(1-x_1)} \frac{1}{(1-x_2)^2} \frac{1}{(1-x_3x_4)} \frac{1}{(1-x_4)^6} \\ 
&&+2x_3x_4^{4} \frac{1}{(1-x_1)} \frac{1}{(1-x_2)} \frac{1}{(1-x_3x_4)^4} \frac{1}{(1-x_4)^4} \\ 
&&-2x_2x_3x_4 \frac{1}{(1-x_1x_2)} \frac{1}{(1-x_2)} \frac{1}{(1-x_3)^3} \frac{1}{(1-x_4)^5} \\ 
&&+2x_2x_3^{2}x_4 \frac{1}{(1-x_1x_2)} \frac{1}{(1-x_2x_3)} \frac{1}{(1-x_3)^3} \frac{1}{(1-x_4)^5} \\ 
&&-3x_2x_3^{2}x_4^{2} \frac{1}{(1-x_1x_2)} \frac{1}{(1-x_2x_3)} \frac{1}{(1-x_3)^2} \frac{1}{(1-x_4)^6} \\ 
&&+4x_2x_3^{2}x_4^{3} \frac{1}{(1-x_1x_2)} \frac{1}{(1-x_2x_3)} \frac{1}{(1-x_3)} \frac{1}{(1-x_4)^7} \\ 
\end{eqnarray*}\begin{eqnarray*}
&&-x_2x_3^{2}x_4^{4} \frac{1}{(1-x_1x_2)} \frac{1}{(1-x_2x_3)} \frac{1}{(1-x_3x_4)^2} \frac{1}{(1-x_4)^6} \\ 
&&-6x_2x_4^{2} \frac{1}{(1-x_1x_2)} \frac{1}{(1-x_2)^2} \frac{1}{(1-x_3)} \frac{1}{(1-x_4)^6} \\ 
&&+x_2x_4^{3} \frac{1}{(1-x_1x_2)} \frac{1}{(1-x_2)^2} \frac{1}{(1-x_3x_4)^3} \frac{1}{(1-x_4)^4} \\ 
&&+3x_2x_4^{3} \frac{1}{(1-x_1x_2)} \frac{1}{(1-x_2)^2} \frac{1}{(1-x_3x_4)^2} \frac{1}{(1-x_4)^5} \\ 
&&-3x_4^{3} \frac{1}{(1-x_1)} \frac{1}{(1-x_2)^2} \frac{1}{(1-x_3x_4)^2} \frac{1}{(1-x_4)^5} \\ 
&&-x_2x_3x_4^{4} \frac{1}{(1-x_1x_2)} \frac{1}{(1-x_2)} \frac{1}{(1-x_3x_4)^4} \frac{1}{(1-x_4)^4} \\ 
&&-2x_2x_3x_4^{4} \frac{1}{(1-x_1x_2)} \frac{1}{(1-x_2)} \frac{1}{(1-x_3x_4)^3} \frac{1}{(1-x_4)^5} \\ 
&&-2x_2x_3x_4^{4} \frac{1}{(1-x_1x_2)} \frac{1}{(1-x_2)} \frac{1}{(1-x_3x_4)^2} \frac{1}{(1-x_4)^6} \\ 
&&+4x_3^{2}x_4^{2} \frac{1}{(1-x_1)} \frac{1}{(1-x_2x_3)} \frac{1}{(1-x_3)^2} \frac{1}{(1-x_4)^6} \\ 
&&-4x_3^{2}x_4^{3} \frac{1}{(1-x_1)} \frac{1}{(1-x_2x_3)} \frac{1}{(1-x_3)} \frac{1}{(1-x_4)^7} \\ 
&&-x_4^{3} \frac{1}{(1-x_1)} \frac{1}{(1-x_2)^2} \frac{1}{(1-x_3x_4)^3} \frac{1}{(1-x_4)^4} \\ 
&&+x_3^{2}x_4^{4} \frac{1}{(1-x_1)} \frac{1}{(1-x_2x_3x_4)^2} \frac{1}{(1-x_3)} \frac{1}{(1-x_4)^6} \\ 
&&-x_3^{2}x_4^{5} \frac{1}{(1-x_1)} \frac{1}{(1-x_2x_3x_4)^2} \frac{1}{(1-x_3x_4)^3} \frac{1}{(1-x_4)^4} \\ 
&&-x_3^{2}x_4^{5} \frac{1}{(1-x_1)} \frac{1}{(1-x_2x_3x_4)^2} \frac{1}{(1-x_3x_4)^2} \frac{1}{(1-x_4)^5} \\ 
&&+6x_4^{2} \frac{1}{(1-x_1)} \frac{1}{(1-x_2)^2} \frac{1}{(1-x_3)} \frac{1}{(1-x_4)^6} \\ 
&&+4x_3x_4^{4} \frac{1}{(1-x_1)} \frac{1}{(1-x_2)} \frac{1}{(1-x_3x_4)^3} \frac{1}{(1-x_4)^5} \\ 
\end{eqnarray*}\begin{eqnarray*}
&&+4x_3x_4^{4} \frac{1}{(1-x_1)} \frac{1}{(1-x_2)} \frac{1}{(1-x_3x_4)^2} \frac{1}{(1-x_4)^6} \\ 
&&-x_2x_3^{3}x_4^{5} \frac{1}{(1-x_1x_2x_3x_4)} \frac{1}{(1-x_2x_3x_4)^2} \frac{1}{(1-x_3)} \frac{1}{(1-x_4)^6} \\ 
&&+x_2x_3^{3}x_4^{6} \frac{1}{(1-x_1x_2x_3x_4)} \frac{1}{(1-x_2x_3x_4)^2} \frac{1}{(1-x_3x_4)^3} \frac{1}{(1-x_4)^4} \\ 
&&+x_2x_3^{3}x_4^{6} \frac{1}{(1-x_1x_2x_3x_4)} \frac{1}{(1-x_2x_3x_4)^2} \frac{1}{(1-x_3x_4)^2} \frac{1}{(1-x_4)^5} \\ 
&&-3x_4 \frac{1}{(1-x_1)} \frac{1}{(1-x_2)^2} \frac{1}{(1-x_3)^2} \frac{1}{(1-x_4)^5} \\ 
&&+ \frac{1}{(1-x_1)} \frac{1}{(1-x_2)^2} \frac{1}{(1-x_3)^3} \frac{1}{(1-x_4)^4} \\ 
&&-x_2 \frac{1}{(1-x_1x_2)} \frac{1}{(1-x_2)^2} \frac{1}{(1-x_3)^3} \frac{1}{(1-x_4)^4} \\ 
&&+x_2x_3 \frac{1}{(1-x_1x_2)} \frac{1}{(1-x_2)} \frac{1}{(1-x_3)^4} \frac{1}{(1-x_4)^4} \\ 
&&-2x_3 \frac{1}{(1-x_1)} \frac{1}{(1-x_2)} \frac{1}{(1-x_3)^4} \frac{1}{(1-x_4)^4} \\ 
&&-x_3^{2}x_4^{3} \frac{1}{(1-x_1)} \frac{1}{(1-x_2x_3)^2} \frac{1}{(1-x_3x_4)} \frac{1}{(1-x_4)^6} \\ 
&&-x_2x_3^{3} \frac{1}{(1-x_1x_2x_3)} \frac{1}{(1-x_2x_3)^2} \frac{1}{(1-x_3)^3} \frac{1}{(1-x_4)^4} \\ 
&&+x_2x_3^{3}x_4 \frac{1}{(1-x_1x_2x_3)} \frac{1}{(1-x_2x_3)^2} \frac{1}{(1-x_3)^2} \frac{1}{(1-x_4)^5} \\ 
\end{eqnarray*}

\begin{eqnarray*}
\prod_{\alpha\in B_2}\frac{1}{1-x^\alpha}&=&+(-x_2^{3}-2x_2^{2}-x_2) \frac{1}{(1-x_1)} \frac{1}{(1-x_2^{2})^3} \\ 
&&-x_2 \frac{1}{(1-x_1x_2)} \frac{1}{(1-x_2)^3} \\ 
&&+(x_2^{5}+2x_2^{4}+x_2^{3}) \frac{1}{(1-x_1x_2^{2})} \frac{1}{(1-x_2^{2})^3} \\ 
&&+ \frac{1}{(1-x_1)} \frac{1}{(1-x_2)^3} \\ 
\end{eqnarray*}

\begin{eqnarray*}
\prod_{\alpha\in C_2}\frac{1}{1-x^\alpha}&=&+(-x_1^{3}x_2^{2}-2x_1^{2}x_2^{2}-x_1x_2^{2}) \frac{1}{(1-x_1^{2}x_2)} \frac{1}{(1-x_2)^3} \\ 
&&+x_2^{2} \frac{1}{(1-x_1x_2)} \frac{1}{(1-x_2)^3} \\ 
&&+(-x_1+1) \frac{1}{(1-x_1)^2} \frac{1}{(1-x_2)^2} \\ 
&&+(x_1^{2}x_2+x_1x_2-x_2) \frac{1}{(1-x_1)} \frac{1}{(1-x_2)^3} \\ 
\end{eqnarray*}

\begin{eqnarray*}
\prod_{\alpha\in B_3} \frac{1}{1-x^\alpha}&= &+(-x_2^{2}x_3^{12}+x_2^{3}x_3^{14}-x_3^{11}-3x_2^{2}x_3^{11}+3x_2^{3}x_3^{13}-3x_3^{10}-3x_2^{2}x_3^{10}+3x_2^{3}x_3^{12}-3x_3^{9}-x_2^{2}x_3^{9}\\
&& 
+x_2^{3}x_3^{11}-x_3^{8})\\
&& \frac{1}{(1-x_1)} \frac{1}{(1-x_2x_3^{2})^2} \frac{1}{(1-x_3^{2})^6} \\ 
&&+(x_2x_3^{13}+3x_2x_3^{12}+3x_2x_3^{11}+x_2x_3^{10}) \frac{1}{(1-x_1x_2x_3^{2})} \frac{1}{(1-x_2x_3^{2})^2} \frac{1}{(1-x_3^{2})^6} \\ 
&&+(-x_2x_3^{10}-6x_2x_3^{9}-15x_2x_3^{8}-20x_2x_3^{7}-15x_2x_3^{6}-6x_2x_3^{5}-x_2x_3^{4}) \frac{1}{(1-x_1x_2)} \frac{1}{(1-x_2x_3)} \frac{1}{(1-x_3^{2})^7} \\ 
&&+(x_3^{10}+x_2x_3^{11}+10x_2^{6}x_3^{12}+10x_2^{5}x_3^{12}+4x_3^{9}+4x_2x_3^{10}+40x_2^{6}x_3^{11}+50x_2^{5}x_3^{11}+x_2^{8}x_3^{9}+6x_2x_3^{9}\\
&& 
+63x_2^{6}x_3^{10}+105x_2^{5}x_3^{10}-18x_3^{7}+4x_2x_3^{8}+49x_2^{6}x_3^{9}+119x_2^{5}x_3^{9}-30x_3^{6}+x_2x_3^{7}+19x_2^{6}x_3^{8}+77x_2^{5}x_3^{8}\\
&& 
+2x_2^{7}x_3^{8}+11x_2^{7}x_3^{9}+3x_2^{6}x_3^{7}+27x_2^{5}x_3^{7}+2x_2^{2}x_3^{11}-2x_2^{3}x_3^{13}+5x_2^{7}x_3^{12}+11x_2^{2}x_3^{10}\\
&& 
-11x_2^{3}x_3^{12}+17x_2^{7}x_3^{11}+21x_2^{2}x_3^{9}-21x_2^{3}x_3^{11}+21x_2^{7}x_3^{10}+17x_2^{2}x_3^{8}-17x_2^{3}x_3^{10}+5x_2^{2}x_3^{7}\\
&& 
-5x_2^{3}x_3^{9}-x_2^{4}x_3^{14}+x_2^{8}x_3^{12}-3x_2^{4}x_3^{13}+3x_2^{8}x_3^{11}-3x_2^{4}x_3^{12}+3x_2^{8}x_3^{10}-x_2^{4}x_3^{11}\\
&& 
-21x_3^{5}-7x_3^{4}+4x_2^{5}x_3^{6}-x_3^{3})\\
&& \frac{1}{(1-x_1)} \frac{1}{(1-x_2)} \frac{1}{(1-x_3^{2})^7} \\ 
&&+(-4x_3^{7}-x_2x_3^{9}-6x_2^{5}x_3^{10}-3x_3^{6}-4x_2x_3^{8}-26x_2^{5}x_3^{9}+6x_3^{5}-6x_2x_3^{7}-45x_2^{5}x_3^{8}+11x_3^{4}-4x_2x_3^{6}\\
&& 
-39x_2^{5}x_3^{7}+6x_3^{3}-x_2x_3^{5}-17x_2^{5}x_3^{6}+x_3^{2}-10x_2^{6}x_3^{7}-3x_2^{5}x_3^{5}-x_3^{8}-2x_2^{2}x_3^{9}+2x_2^{3}x_3^{11}\\
&& 
-4x_2^{6}x_3^{10}-10x_2^{2}x_3^{8}+10x_2^{3}x_3^{10}-14x_2^{6}x_3^{9}-18x_2^{2}x_3^{7}+18x_2^{3}x_3^{9}-18x_2^{6}x_3^{8}-14x_2^{2}x_3^{6}\\
&& 
+14x_2^{3}x_3^{8}-4x_2^{2}x_3^{5}+4x_2^{3}x_3^{7}-2x_2^{6}x_3^{6}+x_2^{4}x_3^{12}-x_2^{7}x_3^{10}+3x_2^{4}x_3^{11}-3x_2^{7}x_3^{9}\\
&& 
+3x_2^{4}x_3^{10}-3x_2^{7}x_3^{8}+x_2^{4}x_3^{9}-x_2^{7}x_3^{7})\\
&& \frac{1}{(1-x_1)} \frac{1}{(1-x_2)^2} \frac{1}{(1-x_3^{2})^6} \\ 
\end{eqnarray*}\begin{eqnarray*}
&&+(-x_2^{8}x_3^{12}-4x_2^{7}x_3^{12}-6x_2^{6}x_3^{12}-4x_2^{5}x_3^{12}-2x_2^{8}x_3^{11}-9x_2^{7}x_3^{11}-16x_2^{6}x_3^{11}-14x_2^{5}x_3^{11}\\
&& 
-x_2^{8}x_3^{10}-6x_2^{7}x_3^{10}-15x_2^{6}x_3^{10}-20x_2^{5}x_3^{10}-x_2^{7}x_3^{9}-6x_2^{6}x_3^{9}-15x_2^{5}x_3^{9}-x_2^{6}x_3^{8}\\
&& 
-6x_2^{5}x_3^{8}-x_2^{5}x_3^{7})\\
&& \frac{1}{(1-x_1)} \frac{1}{(1-x_2x_3)^2} \frac{1}{(1-x_3^{2})^6} \\ 
&&+(4x_2^{7}x_3^{14}+14x_2^{7}x_3^{13}+20x_2^{7}x_3^{12}+15x_2^{7}x_3^{11}+6x_2^{7}x_3^{10}+x_2^{7}x_3^{9}+6x_2^{8}x_3^{14}+16x_2^{8}x_3^{13}\\
&& 
+15x_2^{8}x_3^{12}+6x_2^{8}x_3^{11}+x_2^{8}x_3^{10}+4x_2^{9}x_3^{14}+9x_2^{9}x_3^{13}+6x_2^{9}x_3^{12}+x_2^{9}x_3^{11}+x_2^{10}x_3^{14}\\
&& 
+2x_2^{10}x_3^{13}+x_2^{10}x_3^{12})\\
&& \frac{1}{(1-x_1x_2^{2}x_3^{2})} \frac{1}{(1-x_2x_3)^2} \frac{1}{(1-x_3^{2})^6} \\ 
&&+(-x_2^{6}x_3^{10}+x_2^{7}x_3^{8}-3x_2^{6}x_3^{9}+3x_2^{7}x_3^{7}-3x_2^{6}x_3^{8}+3x_2^{7}x_3^{6}-x_2^{6}x_3^{7}+x_2^{7}x_3^{5})\\
&& \frac{1}{(1-x_1x_2^{2}x_3^{2})} \frac{1}{(1-x_2)^4} \frac{1}{(1-x_3^{2})^4} \\ 
&&+(2x_2^{5}x_3^{15}-5x_2^{9}x_3^{14}-10x_2^{8}x_3^{14}-10x_2^{7}x_3^{14}+11x_2^{5}x_3^{14}-17x_2^{9}x_3^{13}-40x_2^{8}x_3^{13}-50x_2^{7}x_3^{13}\\
&& 
+21x_2^{5}x_3^{13}-21x_2^{9}x_3^{12}-63x_2^{8}x_3^{12}-105x_2^{7}x_3^{12}+17x_2^{5}x_3^{12}-11x_2^{9}x_3^{11}-49x_2^{8}x_3^{11}-119x_2^{7}x_3^{11}\\
&& 
+5x_2^{5}x_3^{11}-2x_2^{9}x_3^{10}-19x_2^{8}x_3^{10}-77x_2^{7}x_3^{10}+x_2^{6}x_3^{16}-x_2^{10}x_3^{14}+3x_2^{6}x_3^{15}-3x_2^{10}x_3^{13}\\
&& 
+3x_2^{6}x_3^{14}-3x_2^{10}x_3^{12}+x_2^{6}x_3^{13}-x_2^{10}x_3^{11}-27x_2^{7}x_3^{9}-4x_2^{7}x_3^{8}-3x_2^{8}x_3^{9})\\
&& \frac{1}{(1-x_1x_2^{2}x_3^{2})} \frac{1}{(1-x_2)} \frac{1}{(1-x_3^{2})^7} \\ 
&&+(-2x_2^{5}x_3^{13}+4x_2^{8}x_3^{12}+6x_2^{7}x_3^{12}-10x_2^{5}x_3^{12}+14x_2^{8}x_3^{11}+26x_2^{7}x_3^{11}-18x_2^{5}x_3^{11}+18x_2^{8}x_3^{10}\\
&& 
+45x_2^{7}x_3^{10}-14x_2^{5}x_3^{10}+10x_2^{8}x_3^{9}+39x_2^{7}x_3^{9}-4x_2^{5}x_3^{9}+2x_2^{8}x_3^{8}+17x_2^{7}x_3^{8}-x_2^{6}x_3^{14}\\
&& 
+x_2^{9}x_3^{12}-3x_2^{6}x_3^{13}+3x_2^{9}x_3^{11}-3x_2^{6}x_3^{12}+3x_2^{9}x_3^{10}-x_2^{6}x_3^{11}+x_2^{9}x_3^{9}+3x_2^{7}x_3^{7})\\
&& \frac{1}{(1-x_1x_2^{2}x_3^{2})} \frac{1}{(1-x_2)^2} \frac{1}{(1-x_3^{2})^6} \\ 
\end{eqnarray*}\begin{eqnarray*}
&&+(2x_2^{5}x_3^{11}-3x_2^{7}x_3^{10}+9x_2^{5}x_3^{10}-11x_2^{7}x_3^{9}+15x_2^{5}x_3^{9}-15x_2^{7}x_3^{8}+11x_2^{5}x_3^{8}-9x_2^{7}x_3^{7}\\
&& 
+3x_2^{5}x_3^{7}-2x_2^{7}x_3^{6}+x_2^{6}x_3^{12}-x_2^{8}x_3^{10}+3x_2^{6}x_3^{11}-3x_2^{8}x_3^{9}+3x_2^{6}x_3^{10}-3x_2^{8}x_3^{8}\\
&& 
+x_2^{6}x_3^{9}-x_2^{8}x_3^{7})\\
&& \frac{1}{(1-x_1x_2^{2}x_3^{2})} \frac{1}{(1-x_2)^3} \frac{1}{(1-x_3^{2})^5} \\ 
&&+(-x_2^{5}x_3^{16}-3x_2^{5}x_3^{15}-3x_2^{5}x_3^{14}-x_2^{5}x_3^{13}) \frac{1}{(1-x_1x_2^{2}x_3^{2})} \frac{1}{(1-x_2x_3^{2})^2} \frac{1}{(1-x_3^{2})^6} \\ 
&&+(-10x_2^{5}x_3^{13}+7x_3^{10}+x_3^{11}-50x_2^{5}x_3^{12}+21x_3^{9}-105x_2^{5}x_3^{11}+35x_3^{8}-119x_2^{5}x_3^{10}+35x_3^{7}-77x_2^{5}x_3^{9}\\
&& 
+21x_3^{6}-27x_2^{5}x_3^{8}+7x_3^{5}-4x_2^{5}x_3^{7}+x_3^{4}-10x_2^{6}x_3^{13}-40x_2^{6}x_3^{12}-63x_2^{6}x_3^{11}-49x_2^{6}x_3^{10}\\
&& 
-19x_2^{6}x_3^{9}-3x_2^{6}x_3^{8}-5x_2^{7}x_3^{13}-17x_2^{7}x_3^{12}-21x_2^{7}x_3^{11}-11x_2^{7}x_3^{10}-2x_2^{7}x_3^{9}-x_2^{8}x_3^{13}\\
&& 
-3x_2^{8}x_3^{12}-3x_2^{8}x_3^{11}-x_2^{8}x_3^{10})\\
&& \frac{1}{(1-x_1)} \frac{1}{(1-x_2x_3)} \frac{1}{(1-x_3^{2})^7} \\ 
&&+(x_2^{4}x_3^{8}-x_2^{5}x_3^{6}+3x_2^{4}x_3^{7}-3x_2^{5}x_3^{5}+3x_2^{4}x_3^{6}-3x_2^{5}x_3^{4}+x_2^{4}x_3^{5}-x_2^{5}x_3^{3}) \frac{1}{(1-x_1)} \frac{1}{(1-x_2)^4} \frac{1}{(1-x_3^{2})^4} \\ 
&&+(-2x_2^{3}x_3^{11}+x_2^{2}x_3^{10}+16x_2x_3^{5}-11x_2^{3}x_3^{10}+4x_2^{2}x_3^{9}+7x_2x_3^{8}-21x_2^{3}x_3^{9}+6x_2^{2}x_3^{8}+18x_2x_3^{7}\\
&& 
-17x_2^{3}x_3^{8}+4x_2^{2}x_3^{7}+23x_2x_3^{6}-5x_2^{3}x_3^{7}+x_2^{2}x_3^{6}+x_2x_3^{9}+6x_2x_3^{4}+x_2x_3^{3})\\
&& \frac{1}{(1-x_1x_2)} \frac{1}{(1-x_2)} \frac{1}{(1-x_3^{2})^7} \\ 
&&+(2x_2^{3}x_3^{9}-x_2^{2}x_3^{8}-6x_2x_3^{3}+10x_2^{3}x_3^{8}-4x_2^{2}x_3^{7}-6x_2x_3^{6}+18x_2^{3}x_3^{7}-6x_2^{2}x_3^{6}-13x_2x_3^{5}\\
&& 
+14x_2^{3}x_3^{6}-4x_2^{2}x_3^{5}-13x_2x_3^{4}+4x_2^{3}x_3^{5}-x_2^{2}x_3^{4}-x_2x_3^{7}-x_2x_3^{2})\\
&& \frac{1}{(1-x_1x_2)} \frac{1}{(1-x_2)^2} \frac{1}{(1-x_3^{2})^6} \\ 
\end{eqnarray*}\begin{eqnarray*}
&&+(-2x_2^{3}x_3^{7}+x_2^{2}x_3^{6}+x_2x_3^{5}-9x_2^{3}x_3^{6}+4x_2^{2}x_3^{5}+5x_2x_3^{4}-15x_2^{3}x_3^{5}+6x_2^{2}x_3^{4}+9x_2x_3^{3}-11x_2^{3}x_3^{4}\\
&& 
+4x_2^{2}x_3^{3}+7x_2x_3^{2}-3x_2^{3}x_3^{3}+x_2^{2}x_3^{2}+2x_2x_3)\\
&& \frac{1}{(1-x_1x_2)} \frac{1}{(1-x_2)^3} \frac{1}{(1-x_3^{2})^5} \\ 
&&+(x_2^{3}x_3^{12}+3x_2^{3}x_3^{11}+3x_2^{3}x_3^{10}+x_2^{3}x_3^{9}) \frac{1}{(1-x_1x_2)} \frac{1}{(1-x_2x_3^{2})^2} \frac{1}{(1-x_3^{2})^6} \\ 
&&+(-x_2^{2}x_3^{12}-x_2x_3^{11}-5x_2^{2}x_3^{11}-5x_2x_3^{10}-10x_2^{2}x_3^{10}-10x_2x_3^{9}-10x_2^{2}x_3^{9}-10x_2x_3^{8}-5x_2^{2}x_3^{8}\\
&& 
-5x_2x_3^{7}-x_2^{2}x_3^{7}-x_2x_3^{6})\\
&& \frac{1}{(1-x_1x_2x_3)} \frac{1}{(1-x_2)} \frac{1}{(1-x_3^{2})^7} \\ 
&&+(6x_2x_3^{3}-3x_3^{3}-3x_2x_3^{4}+x_3^{4}+x_2x_3^{5}-3x_2^{2}x_3^{4}+3x_2^{3}x_3^{6}+2x_2^{2}x_3^{5}-2x_2^{3}x_3^{7}+4x_2^{4}x_3^{7}\\
&& 
-x_2^{4}x_3^{8})\\
&& \frac{1}{(1-x_1)} \frac{1}{(1-x_2x_3)} \frac{1}{(1-x_3)^7} \\ 
&&+(-3x_2^{2}x_3^{4}+2x_2^{2}x_3^{5}+2x_2x_3^{4}-x_2^{2}x_3^{6}-x_2x_3^{5}) \frac{1}{(1-x_1x_2x_3)} \frac{1}{(1-x_2x_3)} \frac{1}{(1-x_3)^7} \\ 
&&+(-2x_2^{2}x_3^{4}+x_2^{2}x_3^{5}+x_2x_3^{4}) \frac{1}{(1-x_1x_2x_3)} \frac{1}{(1-x_2x_3)^2} \frac{1}{(1-x_3)^6} \\ 
&&-x_2^{2}x_3^{4} \frac{1}{(1-x_1x_2x_3)} \frac{1}{(1-x_2x_3)^3} \frac{1}{(1-x_3)^5} \\ 
&&+(x_2^{2}x_3^{10}+x_2x_3^{9}+5x_2^{2}x_3^{9}+5x_2x_3^{8}+10x_2^{2}x_3^{8}+10x_2x_3^{7}+10x_2^{2}x_3^{7}+10x_2x_3^{6}+5x_2^{2}x_3^{6}+5x_2x_3^{5}\\
&& 
+x_2^{2}x_3^{5}+x_2x_3^{4})\\
&& \frac{1}{(1-x_1x_2x_3)} \frac{1}{(1-x_2)^2} \frac{1}{(1-x_3^{2})^6} \\ 
\end{eqnarray*}\begin{eqnarray*}
&&+(-x_2x_3^{6}+2x_2^{2}x_3^{7}-2x_2^{3}x_3^{9}+3x_2^{5}x_3^{8}-x_3^{5}-4x_2x_3^{5}+9x_2^{2}x_3^{6}-9x_2^{3}x_3^{8}+11x_2^{5}x_3^{7}-5x_3^{4}\\
&& 
-6x_2x_3^{4}+15x_2^{2}x_3^{5}-15x_2^{3}x_3^{7}+15x_2^{5}x_3^{6}-9x_3^{3}-4x_2x_3^{3}+11x_2^{2}x_3^{4}-11x_2^{3}x_3^{6}+9x_2^{5}x_3^{5}\\
&& 
-7x_3^{2}-x_2x_3^{2}+3x_2^{2}x_3^{3}-3x_2^{3}x_3^{5}+2x_2^{5}x_3^{4}-2x_3-x_2^{4}x_3^{10}+x_2^{6}x_3^{8}-3x_2^{4}x_3^{9}+3x_2^{6}x_3^{7}\\
&& 
-3x_2^{4}x_3^{8}+3x_2^{6}x_3^{6}-x_2^{4}x_3^{7}+x_2^{6}x_3^{5})\\
&& \frac{1}{(1-x_1)} \frac{1}{(1-x_2)^3} \frac{1}{(1-x_3^{2})^5} \\ 
&&+(12x_2x_3^{8}+x_2x_3^{6}+3x_2x_3^{10}+10x_2x_3^{9}+6x_2x_3^{7}) \frac{1}{(1-x_1x_2x_3^{2})} \frac{1}{(1-x_2)} \frac{1}{(1-x_3^{2})^7} \\ 
&&+(3x_2^{3}x_3^{4}-3x_2^{2}x_3^{3}-2x_2^{3}x_3^{5}+x_2^{2}x_3^{4}+x_2x_3^{3}) \frac{1}{(1-x_1x_2)} \frac{1}{(1-x_2x_3)} \frac{1}{(1-x_3)^7} \\ 
&&+(x_2^{3}x_3^{4}-x_2^{2}x_3^{3}) \frac{1}{(1-x_1x_2)} \frac{1}{(1-x_2x_3)^2} \frac{1}{(1-x_3)^6} \\ 
&&+(-3x_2x_3^{6}-x_2x_3^{5}-x_2x_3^{8}-3x_2x_3^{7}) \frac{1}{(1-x_1x_2x_3^{2})} \frac{1}{(1-x_2)^2} \frac{1}{(1-x_3^{2})^6} \\ 
&&+(-x_2x_3^{13}-6x_2x_3^{12}-15x_2x_3^{11}-20x_2x_3^{10}-15x_2x_3^{9}-6x_2x_3^{8}-x_2x_3^{7}) \frac{1}{(1-x_1x_2x_3^{2})} \frac{1}{(1-x_2x_3)} \frac{1}{(1-x_3^{2})^7} \\ 
&&+(-2x_3^{12}-x_2x_3^{13}-11x_3^{11}-4x_2x_3^{12}-21x_3^{10}-6x_2x_3^{11}-17x_3^{9}-4x_2x_3^{10}-5x_3^{8}-x_2x_3^{9}+5x_2^{3}x_3^{11}\\
&& 
+17x_2^{3}x_3^{12}-2x_2^{2}x_3^{13}+2x_2^{3}x_3^{15}-11x_2^{2}x_3^{12}+11x_2^{3}x_3^{14}-21x_2^{2}x_3^{11}+21x_2^{3}x_3^{13}-17x_2^{2}x_3^{10}\\
&& 
-5x_2^{2}x_3^{9}+x_2^{4}x_3^{16}+3x_2^{4}x_3^{15}+3x_2^{4}x_3^{14}+x_2^{4}x_3^{13})\\
&& \frac{1}{(1-x_1)} \frac{1}{(1-x_2x_3^{2})} \frac{1}{(1-x_3^{2})^7} \\ 
&&+(10x_2^{7}x_3^{15}+50x_2^{7}x_3^{14}+105x_2^{7}x_3^{13}+119x_2^{7}x_3^{12}+77x_2^{7}x_3^{11}+27x_2^{7}x_3^{10}+4x_2^{7}x_3^{9}+10x_2^{8}x_3^{15}\\
&& 
+40x_2^{8}x_3^{14}+63x_2^{8}x_3^{13}+49x_2^{8}x_3^{12}+19x_2^{8}x_3^{11}+3x_2^{8}x_3^{10}+5x_2^{9}x_3^{15}+17x_2^{9}x_3^{14}+21x_2^{9}x_3^{13}\\
&& 
+11x_2^{9}x_3^{12}+2x_2^{9}x_3^{11}+x_2^{10}x_3^{15}+3x_2^{10}x_3^{14}+3x_2^{10}x_3^{13}+x_2^{10}x_3^{12})\\
&& \frac{1}{(1-x_1x_2^{2}x_3^{2})} \frac{1}{(1-x_2x_3)} \frac{1}{(1-x_3^{2})^7} \\ 
\end{eqnarray*}\begin{eqnarray*}
&&+(-2x_2^{5}x_3^{17}-11x_2^{5}x_3^{16}-21x_2^{5}x_3^{15}-17x_2^{5}x_3^{14}-5x_2^{5}x_3^{13}-x_2^{6}x_3^{18}-3x_2^{6}x_3^{17}-3x_2^{6}x_3^{16}\\
&& 
-x_2^{6}x_3^{15})\\
&& \frac{1}{(1-x_1x_2^{2}x_3^{2})} \frac{1}{(1-x_2x_3^{2})} \frac{1}{(1-x_3^{2})^7} \\ 
&&+(-x_2^{11}x_3^{14}-4x_2^{10}x_3^{14}-6x_2^{9}x_3^{14}-4x_2^{8}x_3^{14}-x_2^{7}x_3^{14}-6x_2^{10}x_3^{13}-14x_2^{9}x_3^{13}-16x_2^{8}x_3^{13}\\
&& 
-9x_2^{7}x_3^{13}-2x_2^{6}x_3^{13}-9x_2^{9}x_3^{12}-16x_2^{8}x_3^{12}-14x_2^{7}x_3^{12}-6x_2^{6}x_3^{12}-x_2^{5}x_3^{12}-x_2^{11}x_3^{13}\\
&& 
-2x_2^{10}x_3^{12}-x_2^{9}x_3^{11}-4x_2^{8}x_3^{11}-6x_2^{7}x_3^{11}-4x_2^{6}x_3^{11}-x_2^{5}x_3^{11})\\
&& \frac{1}{(1-x_1)} \frac{1}{(1-x_2^{2}x_3^{2})^3} \frac{1}{(1-x_3^{2})^5} \\ 
&&+(2x_2^{3}x_3^{13}-x_2^{2}x_3^{12}-x_2x_3^{7}+11x_2^{3}x_3^{12}-4x_2^{2}x_3^{11}-x_2x_3^{10}+21x_2^{3}x_3^{11}-6x_2^{2}x_3^{10}-3x_2x_3^{9}\\
&& 
+17x_2^{3}x_3^{10}-4x_2^{2}x_3^{9}-3x_2x_3^{8}+5x_2^{3}x_3^{9}-x_2^{2}x_3^{8})\\
&& \frac{1}{(1-x_1x_2)} \frac{1}{(1-x_2x_3^{2})} \frac{1}{(1-x_3^{2})^7} \\ 
&&+(x_2^{2}x_3^{14}+x_2x_3^{13}+5x_2^{2}x_3^{13}+5x_2x_3^{12}+10x_2^{2}x_3^{12}+10x_2x_3^{11}+10x_2^{2}x_3^{11}+10x_2x_3^{10}+5x_2^{2}x_3^{10}\\
&& 
+5x_2x_3^{9}+x_2^{2}x_3^{9}+x_2x_3^{8})\\
&& \frac{1}{(1-x_1x_2x_3)} \frac{1}{(1-x_2x_3^{2})} \frac{1}{(1-x_3^{2})^7} \\ 
&&+(x_2x_3^{14}+6x_2x_3^{13}+12x_2x_3^{12}+10x_2x_3^{11}+3x_2x_3^{10}) \frac{1}{(1-x_1x_2x_3^{2})} \frac{1}{(1-x_2x_3^{2})} \frac{1}{(1-x_3^{2})^7} \\ 
&&+(-3x_2^{5}x_3^{8}+2x_2^{5}x_3^{9}-4x_2^{6}x_3^{9}+x_2^{6}x_3^{10}) \frac{1}{(1-x_1x_2^{2}x_3^{2})} \frac{1}{(1-x_2x_3)} \frac{1}{(1-x_3)^7} \\ 
&&+(-x_2^{5}x_3^{8}-x_2^{6}x_3^{9}) \frac{1}{(1-x_1x_2^{2}x_3^{2})} \frac{1}{(1-x_2x_3)^2} \frac{1}{(1-x_3)^6} \\ 
\end{eqnarray*}\begin{eqnarray*}
&&+(x_2^{13}x_3^{16}+4x_2^{12}x_3^{16}+6x_2^{11}x_3^{16}+4x_2^{10}x_3^{16}+x_2^{9}x_3^{16}+6x_2^{12}x_3^{15}+14x_2^{11}x_3^{15}+16x_2^{10}x_3^{15}\\
&& 
+9x_2^{9}x_3^{15}+2x_2^{8}x_3^{15}+9x_2^{11}x_3^{14}+16x_2^{10}x_3^{14}+14x_2^{9}x_3^{14}+6x_2^{8}x_3^{14}+x_2^{7}x_3^{14}+x_2^{13}x_3^{15}\\
&& 
+2x_2^{12}x_3^{14}+x_2^{11}x_3^{13}+4x_2^{10}x_3^{13}+6x_2^{9}x_3^{13}+4x_2^{8}x_3^{13}+x_2^{7}x_3^{13})\\
&& \frac{1}{(1-x_1x_2^{2}x_3^{2})} \frac{1}{(1-x_2^{2}x_3^{2})^3} \frac{1}{(1-x_3^{2})^5} \\ 
&&+(x_2^{4}x_3^{3}-x_2^{4}x_3^{4}) \frac{1}{(1-x_1)} \frac{1}{(1-x_2)^4} \frac{1}{(1-x_3)^4} \\ 
&&+(-3x_2^{3}x_3^{3}+3x_2^{2}x_3^{2}+2x_2^{3}x_3^{4}-x_2^{2}x_3^{3}-x_2x_3^{2}) \frac{1}{(1-x_1x_2)} \frac{1}{(1-x_2)} \frac{1}{(1-x_3)^7} \\ 
&&+(2x_2^{3}x_3^{2}-2x_2^{2}x_3-2x_2^{3}x_3^{3}+x_2^{2}x_3^{2}+x_2x_3) \frac{1}{(1-x_1x_2)} \frac{1}{(1-x_2)^2} \frac{1}{(1-x_3)^6} \\ 
&&+(-x_2^{6}x_3^{5}+x_2^{6}x_3^{6}) \frac{1}{(1-x_1x_2^{2}x_3^{2})} \frac{1}{(1-x_2)^4} \frac{1}{(1-x_3)^4} \\ 
&&+(3x_2^{5}x_3^{7}-2x_2^{5}x_3^{8}+4x_2^{6}x_3^{8}-x_2^{6}x_3^{9}) \frac{1}{(1-x_1x_2^{2}x_3^{2})} \frac{1}{(1-x_2)} \frac{1}{(1-x_3)^7} \\ 
&&+(-x_2^{3}x_3+x_2^{2}+2x_2^{3}x_3^{2}-x_2^{2}x_3-x_2) \frac{1}{(1-x_1x_2)} \frac{1}{(1-x_2)^3} \frac{1}{(1-x_3)^5} \\ 
&&+(3x_2^{2}x_3^{3}-2x_2^{2}x_3^{4}-2x_2x_3^{3}+x_2^{2}x_3^{5}+x_2x_3^{4}) \frac{1}{(1-x_1x_2x_3)} \frac{1}{(1-x_2)} \frac{1}{(1-x_3)^7} \\ 
&&+(-x_2^{2}x_3^{2}+x_2^{2}x_3^{3}+x_2x_3^{2}-x_2^{2}x_3^{4}-x_2x_3^{3}) \frac{1}{(1-x_1x_2x_3)} \frac{1}{(1-x_2)^2} \frac{1}{(1-x_3)^6} \\ 
&&+(3x_2x_3^{3}-x_3^{3}-x_2x_3^{4}-x_2^{2}x_3^{4}+x_2^{3}x_3^{6}+x_2^{4}x_3^{7}) \frac{1}{(1-x_1)} \frac{1}{(1-x_2x_3)^2} \frac{1}{(1-x_3)^6} \\ 
&&+x_2x_3^{3} \frac{1}{(1-x_1)} \frac{1}{(1-x_2x_3)^3} \frac{1}{(1-x_3)^5} \\ 
&&+(-6x_2x_3^{2}+3x_3^{2}+3x_2x_3^{3}-x_3^{3}-x_2x_3^{4}+3x_2^{2}x_3^{3}-3x_2^{3}x_3^{5}-2x_2^{2}x_3^{4}+2x_2^{3}x_3^{6}-4x_2^{4}x_3^{6}\\
&& 
+x_2^{4}x_3^{7})\\
&& \frac{1}{(1-x_1)} \frac{1}{(1-x_2)} \frac{1}{(1-x_3)^7} \\ 
\end{eqnarray*}\begin{eqnarray*}
&&+(3x_2x_3-2x_3-2x_2x_3^{2}+x_3^{2}+x_2x_3^{3}-2x_2^{2}x_3^{2}+2x_2^{3}x_3^{4}+2x_2^{2}x_3^{3}-2x_2^{3}x_3^{5}+3x_2^{4}x_3^{5}-x_2^{4}x_3^{6})\\
&& \frac{1}{(1-x_1)} \frac{1}{(1-x_2)^2} \frac{1}{(1-x_3)^6} \\ 
&&+(-x_2+x_2^{2}x_3-x_2^{3}x_3^{3}+x_2x_3-2x_2^{2}x_3^{2}+2x_2^{3}x_3^{4}+1-2x_2^{4}x_3^{4}+x_2^{4}x_3^{5}) \frac{1}{(1-x_1)} \frac{1}{(1-x_2)^3} \frac{1}{(1-x_3)^5} \\ 
&&+(x_2^{5}x_3^{5}-2x_2^{5}x_3^{6}+2x_2^{6}x_3^{6}-x_2^{6}x_3^{7}) \frac{1}{(1-x_1x_2^{2}x_3^{2})} \frac{1}{(1-x_2)^3} \frac{1}{(1-x_3)^5} \\ 
&&+(-2x_2^{5}x_3^{6}+2x_2^{5}x_3^{7}-3x_2^{6}x_3^{7}+x_2^{6}x_3^{8}) \frac{1}{(1-x_1x_2^{2}x_3^{2})} \frac{1}{(1-x_2)^2} \frac{1}{(1-x_3)^6} \\ 
\end{eqnarray*}
\begin{eqnarray*}
\prod_{\alpha\in C_3} \frac{1}{1-x^\alpha}&= &+(-x_1^{6}x_2^{12}x_3^{8}-10x_1^{5}x_2^{11}x_3^{8}-23x_1^{4}x_2^{10}x_3^{8}-2x_1^{5}x_2^{12}x_3^{8}-8x_1^{4}x_2^{11}x_3^{8}-10x_1^{3}x_2^{10}x_3^{8}\\
&& 
-x_1^{4}x_2^{12}x_3^{8}-2x_1^{3}x_2^{11}x_3^{8}-x_1^{2}x_2^{10}x_3^{8}-4x_1^{6}x_2^{11}x_3^{8}-20x_1^{5}x_2^{10}x_3^{8}-32x_1^{4}x_2^{9}x_3^{8}\\
&& 
-20x_1^{3}x_2^{9}x_3^{8}-4x_1^{2}x_2^{9}x_3^{8}-6x_1^{6}x_2^{10}x_3^{8}-20x_1^{5}x_2^{9}x_3^{8}-23x_1^{4}x_2^{8}x_3^{8}-20x_1^{3}x_2^{8}x_3^{8}\\
&& 
-6x_1^{2}x_2^{8}x_3^{8}-4x_1^{6}x_2^{9}x_3^{8}-10x_1^{5}x_2^{8}x_3^{8}-8x_1^{4}x_2^{7}x_3^{8}-10x_1^{3}x_2^{7}x_3^{8}-4x_1^{2}x_2^{7}x_3^{8}\\
&& 
-x_1^{6}x_2^{8}x_3^{8}-2x_1^{5}x_2^{7}x_3^{8}-x_1^{4}x_2^{6}x_3^{8}-2x_1^{3}x_2^{6}x_3^{8}-x_1^{2}x_2^{6}x_3^{8}-3x_1^{6}x_2^{10}x_3^{6}\\
&& 
-18x_1^{5}x_2^{9}x_3^{6}-30x_1^{4}x_2^{8}x_3^{6}-6x_1^{5}x_2^{10}x_3^{6}-18x_1^{4}x_2^{9}x_3^{6}-18x_1^{3}x_2^{8}x_3^{6}-3x_1^{4}x_2^{10}x_3^{6}\\
&& 
-6x_1^{3}x_2^{9}x_3^{6}-3x_1^{2}x_2^{8}x_3^{6}-6x_1^{6}x_2^{9}x_3^{6}-18x_1^{5}x_2^{8}x_3^{6}-18x_1^{4}x_2^{7}x_3^{6}-18x_1^{3}x_2^{7}x_3^{6}\\
&& 
-6x_1^{2}x_2^{7}x_3^{6}-3x_1^{6}x_2^{8}x_3^{6}-6x_1^{5}x_2^{7}x_3^{6}-3x_1^{4}x_2^{6}x_3^{6}-6x_1^{3}x_2^{6}x_3^{6}-3x_1^{2}x_2^{6}x_3^{6}\\
&& 
+3x_1^{6}x_2^{11}x_3^{7}+24x_1^{5}x_2^{10}x_3^{7}+48x_1^{4}x_2^{9}x_3^{7}+6x_1^{5}x_2^{11}x_3^{7}+21x_1^{4}x_2^{10}x_3^{7}+24x_1^{3}x_2^{9}x_3^{7}\\
&& 
+3x_1^{4}x_2^{11}x_3^{7}+6x_1^{3}x_2^{10}x_3^{7}+3x_1^{2}x_2^{9}x_3^{7}+9x_1^{6}x_2^{10}x_3^{7}+36x_1^{5}x_2^{9}x_3^{7}+48x_1^{4}x_2^{8}x_3^{7}\\
&& 
+36x_1^{3}x_2^{8}x_3^{7}+9x_1^{2}x_2^{8}x_3^{7}+9x_1^{6}x_2^{9}x_3^{7}+24x_1^{5}x_2^{8}x_3^{7}+21x_1^{4}x_2^{7}x_3^{7}+24x_1^{3}x_2^{7}x_3^{7}\\
&& 
+9x_1^{2}x_2^{7}x_3^{7}+3x_1^{6}x_2^{8}x_3^{7}+6x_1^{5}x_2^{7}x_3^{7}+3x_1^{4}x_2^{6}x_3^{7}+6x_1^{3}x_2^{6}x_3^{7}+3x_1^{2}x_2^{6}x_3^{7}
)\\
&& \frac{1}{(1-x_1^{2}x_2^{2}x_3)} \frac{1}{(1-x_2^{2}x_3)} \frac{1}{(1-x_3)^7} \\ 
&&+(3x_2^{14}x_3^{9}+18x_2^{13}x_3^{9}+45x_2^{12}x_3^{9}+60x_2^{11}x_3^{9}+45x_2^{10}x_3^{9}+18x_2^{9}x_3^{9}+3x_2^{8}x_3^{9}+10x_2^{12}x_3^{7}\\
&& 
+40x_2^{11}x_3^{7}+60x_2^{10}x_3^{7}+40x_2^{9}x_3^{7}+10x_2^{8}x_3^{7}-10x_2^{13}x_3^{8}-50x_2^{12}x_3^{8}-100x_2^{11}x_3^{8}-100x_2^{10}x_3^{8}\\
&& 
-50x_2^{9}x_3^{8}-10x_2^{8}x_3^{8})\\
&& \frac{1}{(1-x_1x_2^{2}x_3)} \frac{1}{(1-x_2^{2}x_3)} \frac{1}{(1-x_3)^7} \\ 
\end{eqnarray*}\begin{eqnarray*}
&&+(x_1x_2^{9}x_3^{8}+4x_1x_2^{8}x_3^{8}+6x_1x_2^{7}x_3^{8}+4x_1x_2^{6}x_3^{8}+x_1x_2^{5}x_3^{8}-x_2^{10}x_3^{8}-5x_2^{9}x_3^{8}-10x_2^{8}x_3^{8}\\
&& 
-10x_2^{7}x_3^{8}-5x_2^{6}x_3^{8}-x_2^{5}x_3^{8})\\
&& \frac{1}{(1-x_1x_2x_3)} \frac{1}{(1-x_2^{2}x_3)} \frac{1}{(1-x_3)^7} \\ 
&&+(2x_1^{2}x_2^{9}x_3^{5}-2x_1x_2^{8}x_3^{5}+11x_1^{2}x_2^{8}x_3^{5}-11x_1x_2^{7}x_3^{5}+24x_1^{2}x_2^{7}x_3^{5}-24x_1x_2^{6}x_3^{5}\\
&& 
+26x_1^{2}x_2^{6}x_3^{5}-26x_1x_2^{5}x_3^{5}+14x_1^{2}x_2^{5}x_3^{5}-14x_1x_2^{4}x_3^{5}+3x_1^{2}x_2^{4}x_3^{5}-3x_1x_2^{3}x_3^{5}\\
&& 
+x_1^{3}x_2^{9}x_3^{5}+4x_1^{3}x_2^{8}x_3^{5}+6x_1^{3}x_2^{7}x_3^{5}+4x_1^{3}x_2^{6}x_3^{5}+x_1^{3}x_2^{5}x_3^{5}-x_2^{11}x_3^{6}-6x_2^{10}x_3^{6}\\
&& 
-15x_2^{9}x_3^{6}-20x_2^{8}x_3^{6}-15x_2^{7}x_3^{6}-6x_2^{6}x_3^{6}-x_2^{5}x_3^{6}+x_2^{8}x_3^{5}+5x_2^{7}x_3^{5}+10x_2^{6}x_3^{5}\\
&& 
+10x_2^{5}x_3^{5}+5x_2^{4}x_3^{5}+x_2^{3}x_3^{5})\\
&& \frac{1}{(1-x_1x_2)} \frac{1}{(1-x_2^{2}x_3)} \frac{1}{(1-x_3)^7} \\ 
&&+(-x_1^{3}x_2^{4}x_3^{5}-3x_1^{2}x_2^{3}x_3^{5}+3x_1x_2^{2}x_3^{5}-x_2^{2}x_3^{5}+x_2^{4}x_3^{6}+5x_2^{3}x_3^{5}) \frac{1}{(1-x_1x_2)} \frac{1}{(1-x_2x_3)} \frac{1}{(1-x_3)^7} \\ 
&&+(x_1^{4}x_2^{7}x_3^{5}+2x_1^{3}x_2^{6}x_3^{5}+x_1^{2}x_2^{5}x_3^{5}+3x_1^{4}x_2^{6}x_3^{5}+6x_1^{3}x_2^{5}x_3^{5}+3x_1^{2}x_2^{4}x_3^{5}\\
&& 
+3x_1^{4}x_2^{5}x_3^{5}+6x_1^{3}x_2^{4}x_3^{5}+3x_1^{2}x_2^{3}x_3^{5}+x_1^{4}x_2^{4}x_3^{5}+2x_1^{3}x_2^{3}x_3^{5}+x_1^{2}x_2^{2}x_3^{5}\\
&& 
-2x_2^{9}x_3^{6}-10x_2^{8}x_3^{6}-20x_2^{7}x_3^{6}-20x_2^{6}x_3^{6}-10x_2^{5}x_3^{6}-2x_2^{4}x_3^{6}+x_1x_2^{5}x_3^{5}+3x_1x_2^{4}x_3^{5}\\
&& 
+3x_1x_2^{3}x_3^{5}+x_1x_2^{2}x_3^{5}-x_1^{2}x_2^{6}x_3^{4}+x_1x_2^{5}x_3^{4}-5x_1^{2}x_2^{5}x_3^{4}+5x_1x_2^{4}x_3^{4}-11x_1^{2}x_2^{4}x_3^{4}\\
&& 
+9x_1x_2^{3}x_3^{4}-11x_1^{2}x_2^{3}x_3^{4}+7x_1x_2^{2}x_3^{4}-4x_1^{2}x_2^{2}x_3^{4}+2x_1x_2x_3^{4}-x_1^{3}x_2^{6}x_3^{4}-7x_1^{3}x_2^{5}x_3^{4}\\
&& 
-11x_1^{3}x_2^{4}x_3^{4}-5x_1^{3}x_2^{3}x_3^{4}-2x_1^{4}x_2^{6}x_3^{4}-4x_1^{4}x_2^{5}x_3^{4}-2x_1^{4}x_2^{4}x_3^{4}-x_2^{5}x_3^{4}\\
&& 
-4x_2^{4}x_3^{4}-6x_2^{3}x_3^{4}-4x_2^{2}x_3^{4}-x_2x_3^{4}+33x_2^{6}x_3^{5}+22x_2^{5}x_3^{5}+3x_2^{4}x_3^{5}-3x_2^{3}x_3^{5}-x_2^{2}x_3^{5}\\
&& 
+5x_2^{8}x_3^{5}+21x_2^{7}x_3^{5})\\
&& \frac{1}{(1-x_1)^2} \frac{1}{(1-x_2^{2}x_3)} \frac{1}{(1-x_3)^6} \\ 
\end{eqnarray*}\begin{eqnarray*}
&&+(x_1^{5}x_2^{10}x_3^{7}+6x_1^{4}x_2^{9}x_3^{7}+9x_1^{3}x_2^{8}x_3^{7}+x_1^{4}x_2^{10}x_3^{7}+2x_1^{3}x_2^{9}x_3^{7}+x_1^{2}x_2^{8}x_3^{7}\\
&& 
+4x_1^{5}x_2^{9}x_3^{7}+14x_1^{4}x_2^{8}x_3^{7}+16x_1^{3}x_2^{7}x_3^{7}+4x_1^{2}x_2^{7}x_3^{7}+6x_1^{5}x_2^{8}x_3^{7}+16x_1^{4}x_2^{7}x_3^{7}\\
&& 
+14x_1^{3}x_2^{6}x_3^{7}+6x_1^{2}x_2^{6}x_3^{7}+4x_1^{5}x_2^{7}x_3^{7}+9x_1^{4}x_2^{6}x_3^{7}+6x_1^{3}x_2^{5}x_3^{7}+4x_1^{2}x_2^{5}x_3^{7}\\
&& 
+x_1^{5}x_2^{6}x_3^{7}+2x_1^{4}x_2^{5}x_3^{7}+x_1^{3}x_2^{4}x_3^{7}+x_1^{2}x_2^{4}x_3^{7}-3x_2^{12}x_3^{8}-18x_2^{11}x_3^{8}-45x_2^{10}x_3^{8}\\
&& 
-60x_2^{9}x_3^{8}-45x_2^{8}x_3^{8}-18x_2^{7}x_3^{8}-3x_2^{6}x_3^{8}-x_1x_2^{8}x_3^{7}-4x_1x_2^{7}x_3^{7}-6x_1x_2^{6}x_3^{7}-4x_1x_2^{5}x_3^{7}\\
&& 
-x_1x_2^{4}x_3^{7}-2x_1^{2}x_2^{8}x_3^{5}+2x_1x_2^{7}x_3^{5}-11x_1^{2}x_2^{7}x_3^{5}+11x_1x_2^{6}x_3^{5}-21x_1^{2}x_2^{6}x_3^{5}\\
&& 
+24x_1x_2^{5}x_3^{5}-20x_1^{2}x_2^{5}x_3^{5}+26x_1x_2^{4}x_3^{5}-11x_1^{2}x_2^{4}x_3^{5}+14x_1x_2^{3}x_3^{5}-3x_1^{2}x_2^{3}x_3^{5}\\
&& 
+3x_1x_2^{2}x_3^{5}-x_1^{3}x_2^{8}x_3^{5}+2x_1^{3}x_2^{7}x_3^{5}+9x_1^{3}x_2^{6}x_3^{5}+8x_1^{3}x_2^{5}x_3^{5}+2x_1^{3}x_2^{4}x_3^{5}\\
&& 
+3x_1^{5}x_2^{8}x_3^{5}+12x_1^{4}x_2^{7}x_3^{5}+3x_1^{4}x_2^{8}x_3^{5}+6x_1^{5}x_2^{7}x_3^{5}+15x_1^{4}x_2^{6}x_3^{5}+3x_1^{5}x_2^{6}x_3^{5}\\
&& 
+6x_1^{4}x_2^{5}x_3^{5}-3x_1^{5}x_2^{9}x_3^{6}-15x_1^{4}x_2^{8}x_3^{6}-21x_1^{3}x_2^{7}x_3^{6}-3x_1^{4}x_2^{9}x_3^{6}-6x_1^{3}x_2^{8}x_3^{6}\\
&& 
-3x_1^{2}x_2^{7}x_3^{6}-9x_1^{5}x_2^{8}x_3^{6}-27x_1^{4}x_2^{7}x_3^{6}-27x_1^{3}x_2^{6}x_3^{6}-9x_1^{2}x_2^{6}x_3^{6}-9x_1^{5}x_2^{7}x_3^{6}\\
&& 
-21x_1^{4}x_2^{6}x_3^{6}-15x_1^{3}x_2^{5}x_3^{6}-9x_1^{2}x_2^{5}x_3^{6}-3x_1^{5}x_2^{6}x_3^{6}-6x_1^{4}x_2^{5}x_3^{6}-3x_1^{3}x_2^{4}x_3^{6}\\
&& 
-3x_1^{2}x_2^{4}x_3^{6}-x_2^{7}x_3^{5}-5x_2^{6}x_3^{5}-10x_2^{5}x_3^{5}-10x_2^{4}x_3^{5}-5x_2^{3}x_3^{5}-x_2^{2}x_3^{5}+101x_2^{9}x_3^{7}\\
&& 
+105x_2^{8}x_3^{7}+60x_2^{7}x_3^{7}+20x_2^{6}x_3^{7}+5x_2^{5}x_3^{7}+x_2^{4}x_3^{7}+10x_2^{11}x_3^{7}+50x_2^{10}x_3^{7}-9x_2^{10}x_3^{6}\\
&& 
-34x_2^{9}x_3^{6}-45x_2^{8}x_3^{6}-20x_2^{7}x_3^{6}+5x_2^{6}x_3^{6}+6x_2^{5}x_3^{6}+x_2^{4}x_3^{6})\\
&& \frac{1}{(1-x_1)} \frac{1}{(1-x_2^{2}x_3)} \frac{1}{(1-x_3)^7} \\ 
\end{eqnarray*}\begin{eqnarray*}
&&+(x_1^{6}x_2^{7}x_3^{8}+2x_1^{5}x_2^{6}x_3^{8}+x_1^{4}x_2^{5}x_3^{8}+2x_1^{5}x_2^{7}x_3^{8}+4x_1^{4}x_2^{6}x_3^{8}+2x_1^{3}x_2^{5}x_3^{8}\\
&& 
+x_1^{4}x_2^{7}x_3^{8}+2x_1^{3}x_2^{6}x_3^{8}+x_1^{2}x_2^{5}x_3^{8})\\
&& \frac{1}{(1-x_1^{2}x_2^{2}x_3)} \frac{1}{(1-x_2x_3)} \frac{1}{(1-x_3)^7} \\ 
&&+(-2x_1^{6}x_2^{10}x_3^{6}-12x_1^{5}x_2^{9}x_3^{6}-20x_1^{4}x_2^{8}x_3^{6}-4x_1^{5}x_2^{10}x_3^{6}-12x_1^{4}x_2^{9}x_3^{6}-12x_1^{3}x_2^{8}x_3^{6}\\
&& 
-2x_1^{4}x_2^{10}x_3^{6}-4x_1^{3}x_2^{9}x_3^{6}-2x_1^{2}x_2^{8}x_3^{6}-4x_1^{6}x_2^{9}x_3^{6}-12x_1^{5}x_2^{8}x_3^{6}-12x_1^{4}x_2^{7}x_3^{6}\\
&& 
-12x_1^{3}x_2^{7}x_3^{6}-4x_1^{2}x_2^{7}x_3^{6}-2x_1^{6}x_2^{8}x_3^{6}-4x_1^{5}x_2^{7}x_3^{6}-2x_1^{4}x_2^{6}x_3^{6}-4x_1^{3}x_2^{6}x_3^{6}\\
&& 
-2x_1^{2}x_2^{6}x_3^{6}+x_1^{6}x_2^{11}x_3^{7}+8x_1^{5}x_2^{10}x_3^{7}+16x_1^{4}x_2^{9}x_3^{7}+2x_1^{5}x_2^{11}x_3^{7}+7x_1^{4}x_2^{10}x_3^{7}\\
&& 
+8x_1^{3}x_2^{9}x_3^{7}+x_1^{4}x_2^{11}x_3^{7}+2x_1^{3}x_2^{10}x_3^{7}+x_1^{2}x_2^{9}x_3^{7}+3x_1^{6}x_2^{10}x_3^{7}+12x_1^{5}x_2^{9}x_3^{7}\\
&& 
+16x_1^{4}x_2^{8}x_3^{7}+12x_1^{3}x_2^{8}x_3^{7}+3x_1^{2}x_2^{8}x_3^{7}+3x_1^{6}x_2^{9}x_3^{7}+8x_1^{5}x_2^{8}x_3^{7}+7x_1^{4}x_2^{7}x_3^{7}\\
&& 
+8x_1^{3}x_2^{7}x_3^{7}+3x_1^{2}x_2^{7}x_3^{7}+x_1^{6}x_2^{8}x_3^{7}+2x_1^{5}x_2^{7}x_3^{7}+x_1^{4}x_2^{6}x_3^{7}+2x_1^{3}x_2^{6}x_3^{7}\\
&& 
+x_1^{2}x_2^{6}x_3^{7})\\
&& \frac{1}{(1-x_1^{2}x_2^{2}x_3)} \frac{1}{(1-x_2^{2}x_3)^2} \frac{1}{(1-x_3)^6} \\ 
&&+(-x_1^{6}x_2^{10}x_3^{6}-6x_1^{5}x_2^{9}x_3^{6}-10x_1^{4}x_2^{8}x_3^{6}-2x_1^{5}x_2^{10}x_3^{6}-6x_1^{4}x_2^{9}x_3^{6}-6x_1^{3}x_2^{8}x_3^{6}\\
&& 
-x_1^{4}x_2^{10}x_3^{6}-2x_1^{3}x_2^{9}x_3^{6}-x_1^{2}x_2^{8}x_3^{6}-2x_1^{6}x_2^{9}x_3^{6}-6x_1^{5}x_2^{8}x_3^{6}-6x_1^{4}x_2^{7}x_3^{6}\\
&& 
-6x_1^{3}x_2^{7}x_3^{6}-2x_1^{2}x_2^{7}x_3^{6}-x_1^{6}x_2^{8}x_3^{6}-2x_1^{5}x_2^{7}x_3^{6}-x_1^{4}x_2^{6}x_3^{6}-2x_1^{3}x_2^{6}x_3^{6}\\
&& 
-x_1^{2}x_2^{6}x_3^{6})\\
&& \frac{1}{(1-x_1^{2}x_2^{2}x_3)} \frac{1}{(1-x_2^{2}x_3)^3} \frac{1}{(1-x_3)^5} \\ 
\end{eqnarray*}\begin{eqnarray*}
&&+(-x_1^{6}x_2^{8}x_3^{4}+x_1^{4}x_2^{5}x_3^{4}+3x_1^{4}x_2^{6}x_3^{4}-2x_1^{5}x_2^{8}x_3^{4}-3x_1^{4}x_2^{7}x_3^{4}+2x_1^{3}x_2^{5}x_3^{4}\\
&& 
-x_1^{4}x_2^{8}x_3^{4}-2x_1^{3}x_2^{7}x_3^{4}-x_1^{2}x_2^{6}x_3^{4}+x_1^{6}x_2^{7}x_3^{4}+2x_1^{5}x_2^{6}x_3^{4}+x_1^{2}x_2^{5}x_3^{4})\\
&& \frac{1}{(1-x_1^{2}x_2^{2}x_3)} \frac{1}{(1-x_2)^4} \frac{1}{(1-x_3)^4} \\ 
&&+(x_1^{6}x_2^{9}x_3^{5}+4x_1^{5}x_2^{8}x_3^{5}+4x_1^{4}x_2^{7}x_3^{5}+2x_1^{5}x_2^{9}x_3^{5}+5x_1^{4}x_2^{8}x_3^{5}+4x_1^{3}x_2^{7}x_3^{5}\\
&& 
+x_1^{4}x_2^{9}x_3^{5}+2x_1^{3}x_2^{8}x_3^{5}+x_1^{2}x_2^{7}x_3^{5}+x_1^{6}x_2^{8}x_3^{5}-x_1^{4}x_2^{5}x_3^{5}-3x_1^{4}x_2^{6}x_3^{5}\\
&& 
-2x_1^{3}x_2^{5}x_3^{5}+x_1^{2}x_2^{6}x_3^{5}-x_1^{6}x_2^{7}x_3^{5}-2x_1^{5}x_2^{6}x_3^{5}-x_1^{2}x_2^{5}x_3^{5}-x_1^{6}x_2^{8}x_3^{4}\\
&& 
-2x_1^{5}x_2^{7}x_3^{4}-x_1^{4}x_2^{6}x_3^{4}-2x_1^{5}x_2^{8}x_3^{4}-4x_1^{4}x_2^{7}x_3^{4}-2x_1^{3}x_2^{6}x_3^{4}-x_1^{4}x_2^{8}x_3^{4}\\
&& 
-2x_1^{3}x_2^{7}x_3^{4}-x_1^{2}x_2^{6}x_3^{4})\\
&& \frac{1}{(1-x_1^{2}x_2^{2}x_3)} \frac{1}{(1-x_2)^3} \frac{1}{(1-x_3)^5} \\ 
&&+(-x_1^{6}x_2^{10}x_3^{6}-6x_1^{5}x_2^{9}x_3^{6}-10x_1^{4}x_2^{8}x_3^{6}-2x_1^{5}x_2^{10}x_3^{6}-6x_1^{4}x_2^{9}x_3^{6}-6x_1^{3}x_2^{8}x_3^{6}\\
&& 
-x_1^{4}x_2^{10}x_3^{6}-2x_1^{3}x_2^{9}x_3^{6}-x_1^{2}x_2^{8}x_3^{6}-2x_1^{6}x_2^{9}x_3^{6}-6x_1^{5}x_2^{8}x_3^{6}-5x_1^{4}x_2^{7}x_3^{6}\\
&& 
-6x_1^{3}x_2^{7}x_3^{6}-2x_1^{2}x_2^{7}x_3^{6}-x_1^{6}x_2^{8}x_3^{6}+x_1^{4}x_2^{5}x_3^{6}+3x_1^{4}x_2^{6}x_3^{6}+2x_1^{3}x_2^{5}x_3^{6}\\
&& 
-x_1^{2}x_2^{6}x_3^{6}+x_1^{6}x_2^{7}x_3^{6}+2x_1^{5}x_2^{6}x_3^{6}+x_1^{2}x_2^{5}x_3^{6}-x_1^{6}x_2^{8}x_3^{4}-2x_1^{5}x_2^{7}x_3^{4}\\
&& 
-x_1^{4}x_2^{6}x_3^{4}-2x_1^{5}x_2^{8}x_3^{4}-4x_1^{4}x_2^{7}x_3^{4}-2x_1^{3}x_2^{6}x_3^{4}-x_1^{4}x_2^{8}x_3^{4}-2x_1^{3}x_2^{7}x_3^{4}\\
&& 
-x_1^{2}x_2^{6}x_3^{4}+2x_1^{6}x_2^{9}x_3^{5}+8x_1^{5}x_2^{8}x_3^{5}+10x_1^{4}x_2^{7}x_3^{5}+4x_1^{5}x_2^{9}x_3^{5}+10x_1^{4}x_2^{8}x_3^{5}\\
&& 
+8x_1^{3}x_2^{7}x_3^{5}+2x_1^{4}x_2^{9}x_3^{5}+4x_1^{3}x_2^{8}x_3^{5}+2x_1^{2}x_2^{7}x_3^{5}+2x_1^{6}x_2^{8}x_3^{5}+4x_1^{5}x_2^{7}x_3^{5}\\
&& 
+2x_1^{4}x_2^{6}x_3^{5}+4x_1^{3}x_2^{6}x_3^{5}+2x_1^{2}x_2^{6}x_3^{5})\\
&& \frac{1}{(1-x_1^{2}x_2^{2}x_3)} \frac{1}{(1-x_2)^2} \frac{1}{(1-x_3)^6} \\ 
\end{eqnarray*}\begin{eqnarray*}
&&+(x_1^{6}x_2^{11}x_3^{7}+8x_1^{5}x_2^{10}x_3^{7}+16x_1^{4}x_2^{9}x_3^{7}+2x_1^{5}x_2^{11}x_3^{7}+7x_1^{4}x_2^{10}x_3^{7}+8x_1^{3}x_2^{9}x_3^{7}\\
&& 
+x_1^{4}x_2^{11}x_3^{7}+2x_1^{3}x_2^{10}x_3^{7}+x_1^{2}x_2^{9}x_3^{7}+3x_1^{6}x_2^{10}x_3^{7}+12x_1^{5}x_2^{9}x_3^{7}+16x_1^{4}x_2^{8}x_3^{7}\\
&& 
+12x_1^{3}x_2^{8}x_3^{7}+3x_1^{2}x_2^{8}x_3^{7}+3x_1^{6}x_2^{9}x_3^{7}+8x_1^{5}x_2^{8}x_3^{7}+6x_1^{4}x_2^{7}x_3^{7}+8x_1^{3}x_2^{7}x_3^{7}\\
&& 
+3x_1^{2}x_2^{7}x_3^{7}+x_1^{6}x_2^{8}x_3^{7}-x_1^{4}x_2^{5}x_3^{7}-3x_1^{4}x_2^{6}x_3^{7}-2x_1^{3}x_2^{5}x_3^{7}+x_1^{2}x_2^{6}x_3^{7}\\
&& 
-x_1^{6}x_2^{7}x_3^{7}-2x_1^{5}x_2^{6}x_3^{7}-x_1^{2}x_2^{5}x_3^{7}+3x_1^{6}x_2^{9}x_3^{5}+12x_1^{5}x_2^{8}x_3^{5}+15x_1^{4}x_2^{7}x_3^{5}\\
&& 
+6x_1^{5}x_2^{9}x_3^{5}+15x_1^{4}x_2^{8}x_3^{5}+12x_1^{3}x_2^{7}x_3^{5}+3x_1^{4}x_2^{9}x_3^{5}+6x_1^{3}x_2^{8}x_3^{5}+3x_1^{2}x_2^{7}x_3^{5}\\
&& 
+3x_1^{6}x_2^{8}x_3^{5}+6x_1^{5}x_2^{7}x_3^{5}+3x_1^{4}x_2^{6}x_3^{5}+6x_1^{3}x_2^{6}x_3^{5}+3x_1^{2}x_2^{6}x_3^{5}-3x_1^{6}x_2^{10}x_3^{6}\\
&& 
-18x_1^{5}x_2^{9}x_3^{6}-30x_1^{4}x_2^{8}x_3^{6}-6x_1^{5}x_2^{10}x_3^{6}-18x_1^{4}x_2^{9}x_3^{6}-18x_1^{3}x_2^{8}x_3^{6}-3x_1^{4}x_2^{10}x_3^{6}\\
&& 
-6x_1^{3}x_2^{9}x_3^{6}-3x_1^{2}x_2^{8}x_3^{6}-6x_1^{6}x_2^{9}x_3^{6}-18x_1^{5}x_2^{8}x_3^{6}-18x_1^{4}x_2^{7}x_3^{6}-18x_1^{3}x_2^{7}x_3^{6}\\
&& 
-6x_1^{2}x_2^{7}x_3^{6}-3x_1^{6}x_2^{8}x_3^{6}-6x_1^{5}x_2^{7}x_3^{6}-3x_1^{4}x_2^{6}x_3^{6}-6x_1^{3}x_2^{6}x_3^{6}-3x_1^{2}x_2^{6}x_3^{6})\\
&& \frac{1}{(1-x_1^{2}x_2^{2}x_3)} \frac{1}{(1-x_2)} \frac{1}{(1-x_3)^7} \\ 
&&+(x_2^{5}x_3^{3}-x_2^{4}x_3^{3}) \frac{1}{(1-x_1x_2x_3)} \frac{1}{(1-x_2)^5} \frac{1}{(1-x_3)^3} \\ 
&&+(-x_2^{5}+x_2^{4}) \frac{1}{(1-x_1x_2)} \frac{1}{(1-x_2)^6} \frac{1}{(1-x_3)^2} \\ 
&&+(-2x_2^{5}x_3^{2}-2x_2^{4}x_3^{2}-x_1^{2}x_2^{2}+x_1x_2+x_1^{2}x_2-x_1-2x_2+1+2x_2^{3}x_3^{2}+x_2x_3+2x_2^{4}x_3+x_2^{3}x_3-x_2^{2}x_3
)\\
&& \frac{1}{(1-x_1)^2} \frac{1}{(1-x_2)^4} \frac{1}{(1-x_3)^3} \\ 
\end{eqnarray*}\begin{eqnarray*}
&&+(3x_2^{7}x_3^{3}+3x_2^{6}x_3^{3}+2x_1^{2}x_2^{3}-2x_1x_2^{2}-2x_1^{2}x_2^{2}+2x_1x_2+2x_2^{2}-x_2-2x_2^{4}x_3^{2}+x_2^{3}x_3^{2}-2x_2^{6}x_3^{2}\\
&& 
-3x_2^{5}x_3^{3}-x_2^{5}x_3-x_2^{4}x_3+x_2^{3}x_3)\\
&& \frac{1}{(1-x_1)} \frac{1}{(1-x_2)^5} \frac{1}{(1-x_3)^3} \\ 
&&+(-5x_2^{6}x_3^{8}-3x_2^{7}x_3^{9}) \frac{1}{(1-x_1x_2^{2}x_3)} \frac{1}{(1-x_2x_3)} \frac{1}{(1-x_3)^7} \\ 
&&-x_2^{6}x_3^{8} \frac{1}{(1-x_1x_2^{2}x_3)} \frac{1}{(1-x_2x_3)^2} \frac{1}{(1-x_3)^6} \\ 
&&+(4x_2^{12}x_3^{7}+16x_2^{11}x_3^{7}+24x_2^{10}x_3^{7}+16x_2^{9}x_3^{7}+4x_2^{8}x_3^{7}-2x_2^{13}x_3^{8}-10x_2^{12}x_3^{8}-20x_2^{11}x_3^{8}\\
&& 
-20x_2^{10}x_3^{8}-10x_2^{9}x_3^{8}-2x_2^{8}x_3^{8})\\
&& \frac{1}{(1-x_1x_2^{2}x_3)} \frac{1}{(1-x_2^{2}x_3)^2} \frac{1}{(1-x_3)^6} \\ 
&&+(x_2^{12}x_3^{7}+4x_2^{11}x_3^{7}+6x_2^{10}x_3^{7}+4x_2^{9}x_3^{7}+x_2^{8}x_3^{7}) \frac{1}{(1-x_1x_2^{2}x_3)} \frac{1}{(1-x_2^{2}x_3)^3} \frac{1}{(1-x_3)^5} \\ 
&&+(3x_2^{8}x_3^{3}-3x_2^{7}x_3^{3}) \frac{1}{(1-x_1x_2^{2}x_3)} \frac{1}{(1-x_2)^6} \frac{1}{(1-x_3)^2} \\ 
&&+(-3x_2^{9}x_3^{4}-3x_2^{8}x_3^{4}+x_2^{6}x_3^{3}+3x_2^{7}x_3^{4}+2x_2^{8}x_3^{3}) \frac{1}{(1-x_1x_2^{2}x_3)} \frac{1}{(1-x_2)^5} \frac{1}{(1-x_3)^3} \\ 
&&+(3x_2^{10}x_3^{5}+6x_2^{9}x_3^{5}+3x_2^{8}x_3^{5}-2x_2^{6}x_3^{4}-3x_2^{7}x_3^{5}+x_2^{8}x_3^{3}-4x_2^{9}x_3^{4}-4x_2^{8}x_3^{4})\\
&& \frac{1}{(1-x_1x_2^{2}x_3)} \frac{1}{(1-x_2)^4} \frac{1}{(1-x_3)^4} \\ 
&&+(-3x_2^{11}x_3^{6}-9x_2^{10}x_3^{6}-9x_2^{9}x_3^{6}-3x_2^{8}x_3^{6}+3x_2^{6}x_3^{5}+3x_2^{7}x_3^{6}-3x_2^{9}x_3^{4}-3x_2^{8}x_3^{4}\\
&& 
+6x_2^{10}x_3^{5}+12x_2^{9}x_3^{5}+6x_2^{8}x_3^{5})\\
&& \frac{1}{(1-x_1x_2^{2}x_3)} \frac{1}{(1-x_2)^3} \frac{1}{(1-x_3)^5} \\ 
\end{eqnarray*}\begin{eqnarray*}
&&+(3x_2^{12}x_3^{7}+12x_2^{11}x_3^{7}+18x_2^{10}x_3^{7}+12x_2^{9}x_3^{7}+3x_2^{8}x_3^{7}-4x_2^{6}x_3^{6}-3x_2^{7}x_3^{7}+6x_2^{10}x_3^{5}\\
&& 
+12x_2^{9}x_3^{5}+6x_2^{8}x_3^{5}-8x_2^{11}x_3^{6}-24x_2^{10}x_3^{6}-24x_2^{9}x_3^{6}-8x_2^{8}x_3^{6})\\
&& \frac{1}{(1-x_1x_2^{2}x_3)} \frac{1}{(1-x_2)^2} \frac{1}{(1-x_3)^6} \\ 
&&+(-3x_2^{13}x_3^{8}-15x_2^{12}x_3^{8}-30x_2^{11}x_3^{8}-30x_2^{10}x_3^{8}-15x_2^{9}x_3^{8}-3x_2^{8}x_3^{8}+5x_2^{6}x_3^{7}+3x_2^{7}x_3^{8}\\
&& 
-10x_2^{11}x_3^{6}-30x_2^{10}x_3^{6}-30x_2^{9}x_3^{6}-10x_2^{8}x_3^{6}+10x_2^{12}x_3^{7}+40x_2^{11}x_3^{7}+60x_2^{10}x_3^{7}+40x_2^{9}x_3^{7}\\
&& 
+10x_2^{8}x_3^{7})\\
&& \frac{1}{(1-x_1x_2^{2}x_3)} \frac{1}{(1-x_2)} \frac{1}{(1-x_3)^7} \\ 
&&+(-4x_1x_3^{4}+2x_1^{2}x_2x_3^{4}-x_1^{4}x_2^{3}x_3^{5}-2x_1^{3}x_2^{2}x_3^{5}-x_1^{2}x_2x_3^{5}-x_1x_2x_3^{5}+x_1^{3}x_2^{2}x_3^{4}\\
&& 
+4x_3^{4}+x_2x_3^{5}+2x_2^{3}x_3^{6}+3x_2^{2}x_3^{5}-4x_2x_3^{4})\\
&& \frac{1}{(1-x_1)^2} \frac{1}{(1-x_2x_3)} \frac{1}{(1-x_3)^6} \\ 
&&+(3x_1x_2^{2}x_3^{6}+3x_1^{2}x_2^{2}x_3^{5}+x_1^{3}x_2^{3}x_3^{5}-x_1^{5}x_2^{5}x_3^{7}-2x_1^{4}x_2^{4}x_3^{7}-x_1^{3}x_2^{3}x_3^{7}\\
&& 
-x_1^{4}x_2^{5}x_3^{7}-2x_1^{3}x_2^{4}x_3^{7}-x_1^{2}x_2^{3}x_3^{7}+x_1x_2^{3}x_3^{7}-5x_2x_3^{5}-x_2^{3}x_3^{7}-4x_2^{2}x_3^{6}\\
&& 
+3x_2^{5}x_3^{8}+5x_2^{4}x_3^{7}-x_2^{3}x_3^{6}-5x_2^{2}x_3^{5})\\
&& \frac{1}{(1-x_1)} \frac{1}{(1-x_2x_3)} \frac{1}{(1-x_3)^7} \\ 
&&+(-3x_1x_2^{2}x_3^{6}-3x_1x_2^{3}x_3^{7}-x_1x_2^{4}x_3^{8}+6x_2^{2}x_3^{6}+4x_2^{3}x_3^{7}+x_2^{4}x_3^{8}) \frac{1}{(1-x_1x_2x_3)} \frac{1}{(1-x_2x_3)} \frac{1}{(1-x_3)^7} \\ 
&&+(-2x_1x_2^{2}x_3^{6}-x_1x_2^{3}x_3^{7}+3x_2^{2}x_3^{6}+x_2^{3}x_3^{7}) \frac{1}{(1-x_1x_2x_3)} \frac{1}{(1-x_2x_3)^2} \frac{1}{(1-x_3)^6} \\ 
\end{eqnarray*}\begin{eqnarray*}
&&+(-x_1x_2^{2}x_3^{6}+x_2^{2}x_3^{6}) \frac{1}{(1-x_1x_2x_3)} \frac{1}{(1-x_2x_3)^3} \frac{1}{(1-x_3)^5} \\ 
&&+(x_1x_2^{5}x_3^{4}-x_2^{6}x_3^{4}-x_2^{5}x_3^{4}-x_1x_2^{4}x_3^{4}+x_2^{3}x_3^{3}+x_2^{4}x_3^{4}) \frac{1}{(1-x_1x_2x_3)} \frac{1}{(1-x_2)^4} \frac{1}{(1-x_3)^4} \\ 
&&+(-x_1x_2^{6}x_3^{5}-x_1x_2^{5}x_3^{5}+x_2^{7}x_3^{5}+2x_2^{6}x_3^{5}+x_2^{5}x_3^{5}+x_1x_2^{3}x_3^{4}+x_1x_2^{4}x_3^{5}-x_2^{2}x_3^{3}\\
&& 
-2x_2^{3}x_3^{4}-x_2^{4}x_3^{5})\\
&& \frac{1}{(1-x_1x_2x_3)} \frac{1}{(1-x_2)^3} \frac{1}{(1-x_3)^5} \\ 
&&+(x_1x_2^{7}x_3^{6}+2x_1x_2^{6}x_3^{6}+x_1x_2^{5}x_3^{6}-x_2^{8}x_3^{6}-3x_2^{7}x_3^{6}-3x_2^{6}x_3^{6}-x_2^{5}x_3^{6}-x_1x_2^{2}x_3^{4}\\
&& 
-2x_1x_2^{3}x_3^{5}-x_1x_2^{4}x_3^{6}+3x_2^{2}x_3^{4}+3x_2^{3}x_3^{5}+x_2^{4}x_3^{6})\\
&& \frac{1}{(1-x_1x_2x_3)} \frac{1}{(1-x_2)^2} \frac{1}{(1-x_3)^6} \\ 
&&+(-x_1x_2^{8}x_3^{7}-3x_1x_2^{7}x_3^{7}-3x_1x_2^{6}x_3^{7}-x_1x_2^{5}x_3^{7}+x_2^{9}x_3^{7}+4x_2^{8}x_3^{7}+6x_2^{7}x_3^{7}+4x_2^{6}x_3^{7}\\
&& 
+x_2^{5}x_3^{7}+3x_1x_2^{2}x_3^{5}+3x_1x_2^{3}x_3^{6}+x_1x_2^{4}x_3^{7}-6x_2^{2}x_3^{5}-4x_2^{3}x_3^{6}-x_2^{4}x_3^{7})\\
&& \frac{1}{(1-x_1x_2x_3)} \frac{1}{(1-x_2)} \frac{1}{(1-x_3)^7} \\ 
&&+(x_1x_2^{2}x_3^{6}+2x_1x_2x_3^{5}-x_2^{2}x_3^{6}-3x_2x_3^{5}+x_2^{4}x_3^{7}-x_2^{2}x_3^{5}) \frac{1}{(1-x_1)} \frac{1}{(1-x_2x_3)^2} \frac{1}{(1-x_3)^6} \\ 
&&+(x_1x_2x_3^{5}-x_2x_3^{5}) \frac{1}{(1-x_1)} \frac{1}{(1-x_2x_3)^3} \frac{1}{(1-x_3)^5} \\ 
&&+(-x_1^{2}x_2^{8}x_3^{4}+x_1x_2^{7}x_3^{4}-4x_1^{2}x_2^{7}x_3^{4}+4x_1x_2^{6}x_3^{4}-6x_1^{2}x_2^{6}x_3^{4}+6x_1x_2^{5}x_3^{4}-4x_1^{2}x_2^{5}x_3^{4}\\
&& 
+4x_1x_2^{4}x_3^{4}-x_1^{2}x_2^{4}x_3^{4}+x_1x_2^{3}x_3^{4})\\
&& \frac{1}{(1-x_1x_2)^2} \frac{1}{(1-x_2^{2}x_3)} \frac{1}{(1-x_3)^6} \\ 
\end{eqnarray*}\begin{eqnarray*}
&&+(x_1^{2}x_2^{7}x_3^{3}-x_1x_2^{6}x_3^{3}+3x_1^{2}x_2^{6}x_3^{3}-3x_1x_2^{5}x_3^{3}+3x_1^{2}x_2^{5}x_3^{3}-3x_1x_2^{4}x_3^{3}+x_1^{2}x_2^{4}x_3^{3}\\
&& 
-x_1x_2^{3}x_3^{3}-x_1^{2}x_2^{3}x_3^{3}+x_1x_2^{2}x_3^{3})\\
&& \frac{1}{(1-x_1x_2)^2} \frac{1}{(1-x_2)} \frac{1}{(1-x_3)^6} \\ 
&&+(-2x_1^{2}x_2^{4}+2x_1x_2^{3}+x_2^{6}x_3+x_2^{5}x_3+2x_1^{2}x_2^{3}-2x_1x_2^{2}-2x_2^{3}+x_2^{2}-x_2^{4}x_3) \frac{1}{(1-x_1x_2)} \frac{1}{(1-x_2)^5} \frac{1}{(1-x_3)^3} \\ 
&&+(2x_1^{2}x_2^{5}x_3-2x_1x_2^{4}x_3+3x_1^{2}x_2^{4}x_3-3x_1x_2^{3}x_3+x_1^{3}x_2^{5}x_3-x_1^{3}x_2^{4}x_3-3x_1^{2}x_2^{3}x_3+3x_1x_2^{2}x_3\\
&& 
-x_2^{7}x_3^{2}-2x_2^{6}x_3^{2}-x_2^{5}x_3^{2}+x_2^{4}x_3+3x_2^{3}x_3-x_2^{2}x_3+x_2^{4}x_3^{2})\\
&& \frac{1}{(1-x_1x_2)} \frac{1}{(1-x_2)^4} \frac{1}{(1-x_3)^4} \\ 
&&+(-2x_1^{2}x_2^{6}x_3^{2}+2x_1x_2^{5}x_3^{2}-5x_1^{2}x_2^{5}x_3^{2}+5x_1x_2^{4}x_3^{2}-3x_1^{2}x_2^{4}x_3^{2}+3x_1x_2^{3}x_3^{2}-x_1^{3}x_2^{6}x_3^{2}\\
&& 
-x_1^{3}x_2^{5}x_3^{2}+x_1^{3}x_2^{4}x_3^{2}+3x_1^{2}x_2^{3}x_3^{2}-3x_1x_2^{2}x_3^{2}+x_2^{8}x_3^{3}+3x_2^{7}x_3^{3}+3x_2^{6}x_3^{3}\\
&& 
+x_2^{5}x_3^{3}-x_2^{5}x_3^{2}-2x_2^{4}x_3^{2}-4x_2^{3}x_3^{2}+x_2^{2}x_3^{2}-x_2^{4}x_3^{3})\\
&& \frac{1}{(1-x_1x_2)} \frac{1}{(1-x_2)^3} \frac{1}{(1-x_3)^5} \\ 
&&+(2x_1^{2}x_2^{7}x_3^{3}-2x_1x_2^{6}x_3^{3}+7x_1^{2}x_2^{6}x_3^{3}-7x_1x_2^{5}x_3^{3}+8x_1^{2}x_2^{5}x_3^{3}-8x_1x_2^{4}x_3^{3}+3x_1^{2}x_2^{4}x_3^{3}\\
&& 
-3x_1x_2^{3}x_3^{3}+x_1^{3}x_2^{7}x_3^{3}+2x_1^{3}x_2^{6}x_3^{3}+x_1^{3}x_2^{5}x_3^{3}-x_1^{3}x_2^{4}x_3^{3}-3x_1^{2}x_2^{3}x_3^{3}\\
&& 
+3x_1x_2^{2}x_3^{3}-x_2^{9}x_3^{4}-4x_2^{8}x_3^{4}-6x_2^{7}x_3^{4}-4x_2^{6}x_3^{4}-x_2^{5}x_3^{4}+x_2^{6}x_3^{3}+3x_2^{5}x_3^{3}+3x_2^{4}x_3^{3}\\
&& 
+5x_2^{3}x_3^{3}-x_2^{2}x_3^{3}+x_2^{4}x_3^{4})\\
&& \frac{1}{(1-x_1x_2)} \frac{1}{(1-x_2)^2} \frac{1}{(1-x_3)^6} \\ 
\end{eqnarray*}\begin{eqnarray*}
&&+(-2x_1^{2}x_2^{8}x_3^{4}+2x_1x_2^{7}x_3^{4}-9x_1^{2}x_2^{7}x_3^{4}+9x_1x_2^{6}x_3^{4}-15x_1^{2}x_2^{6}x_3^{4}+15x_1x_2^{5}x_3^{4}-11x_1^{2}x_2^{5}x_3^{4}\\
&& 
+11x_1x_2^{4}x_3^{4}-3x_1^{2}x_2^{4}x_3^{4}+3x_1x_2^{3}x_3^{4}-x_1^{3}x_2^{8}x_3^{4}-3x_1^{3}x_2^{7}x_3^{4}-3x_1^{3}x_2^{6}x_3^{4}\\
&& 
-x_1^{3}x_2^{5}x_3^{4}+x_1^{3}x_2^{4}x_3^{4}+3x_1^{2}x_2^{3}x_3^{4}-3x_1x_2^{2}x_3^{4}+x_2^{10}x_3^{5}+5x_2^{9}x_3^{5}+10x_2^{8}x_3^{5}\\
&& 
+10x_2^{7}x_3^{5}+5x_2^{6}x_3^{5}+x_2^{5}x_3^{5}-x_2^{7}x_3^{4}-4x_2^{6}x_3^{4}-6x_2^{5}x_3^{4}-4x_2^{4}x_3^{4}-6x_2^{3}x_3^{4}+x_2^{2}x_3^{4}\\
&& 
-x_2^{4}x_3^{5})\\
&& \frac{1}{(1-x_1x_2)} \frac{1}{(1-x_2)} \frac{1}{(1-x_3)^7} \\ 
&&+(-x_1^{2}x_2^{6}x_3^{2}+x_1x_2^{5}x_3^{2}-2x_1^{2}x_2^{5}x_3^{2}+2x_1x_2^{4}x_3^{2}-x_1^{2}x_2^{4}x_3^{2}+x_1x_2^{3}x_3^{2}+x_1^{2}x_2^{3}x_3^{2}\\
&& 
-x_1x_2^{2}x_3^{2})\\
&& \frac{1}{(1-x_1x_2)^2} \frac{1}{(1-x_2)^2} \frac{1}{(1-x_3)^5} \\ 
&&+(x_1^{2}x_2^{5}x_3-x_1x_2^{4}x_3+x_1^{2}x_2^{4}x_3-x_1x_2^{3}x_3-x_1^{2}x_2^{3}x_3+x_1x_2^{2}x_3) \frac{1}{(1-x_1x_2)^2} \frac{1}{(1-x_2)^3} \frac{1}{(1-x_3)^4} \\ 
&&+(-x_1^{2}x_2^{4}+x_1x_2^{3}+x_1^{2}x_2^{3}-x_1x_2^{2}) \frac{1}{(1-x_1x_2)^2} \frac{1}{(1-x_2)^4} \frac{1}{(1-x_3)^3} \\ 
&&+(-x_1^{4}x_2^{6}x_3^{4}-2x_1^{3}x_2^{5}x_3^{4}-x_1^{2}x_2^{4}x_3^{4}-2x_1^{4}x_2^{5}x_3^{4}-4x_1^{3}x_2^{4}x_3^{4}-2x_1^{2}x_2^{3}x_3^{4}\\
&& 
-x_1^{4}x_2^{4}x_3^{4}-2x_1^{3}x_2^{3}x_3^{4}-x_1^{2}x_2^{2}x_3^{4}+x_2^{8}x_3^{5}+4x_2^{7}x_3^{5}+6x_2^{6}x_3^{5}+4x_2^{5}x_3^{5}\\
&& 
+x_2^{4}x_3^{5})\\
&& \frac{1}{(1-x_1)^2} \frac{1}{(1-x_2^{2}x_3)^2} \frac{1}{(1-x_3)^5} \\ 
&&+(-x_1^{4}x_2^{4}x_3^{2}-2x_1^{3}x_2^{3}x_3^{2}-x_1^{2}x_2^{2}x_3^{2}+2x_2^{6}x_3^{3}+4x_2^{5}x_3^{3}+2x_2^{4}x_3^{3}-x_1x_2^{2}x_3^{2}\\
&& 
+x_1^{2}x_2^{3}x_3-x_1x_2^{2}x_3+2x_1^{2}x_2^{2}x_3-2x_1x_2x_3-2x_1^{2}x_2x_3+2x_1x_3+x_1^{3}x_2^{3}x_3+x_1^{4}x_2^{3}x_3^{2}+2x_1^{3}x_2^{2}x_3^{2}\\
&& 
+x_1^{2}x_2x_3^{2}+x_1x_2x_3^{2}-x_1^{3}x_2^{2}x_3+x_2^{2}x_3+3x_2x_3-2x_3-2x_2^{3}x_3^{3}-4x_2^{4}x_3^{2}-x_2x_3^{2}-3x_2^{5}x_3^{2}
)\\
&& \frac{1}{(1-x_1)^2} \frac{1}{(1-x_2)^3} \frac{1}{(1-x_3)^4} \\ 
\end{eqnarray*}\begin{eqnarray*}
&&+(x_1^{4}x_2^{5}x_3^{3}+2x_1^{3}x_2^{4}x_3^{3}+x_1^{2}x_2^{3}x_3^{3}+x_1^{4}x_2^{4}x_3^{3}+2x_1^{3}x_2^{3}x_3^{3}+x_1^{2}x_2^{2}x_3^{3}\\
&& 
-2x_2^{7}x_3^{4}-6x_2^{6}x_3^{4}-6x_2^{5}x_3^{4}-2x_2^{4}x_3^{4}+x_1x_2^{3}x_3^{3}+x_1x_2^{2}x_3^{3}-3x_1x_3^{2}-x_1^{2}x_2^{4}x_3^{2}\\
&& 
+x_1x_2^{3}x_3^{2}-3x_1^{2}x_2^{3}x_3^{2}+3x_1x_2^{2}x_3^{2}-3x_1^{2}x_2^{2}x_3^{2}+2x_1x_2x_3^{2}+2x_1^{2}x_2x_3^{2}-x_1^{3}x_2^{4}x_3^{2}\\
&& 
-3x_1^{3}x_2^{3}x_3^{2}-x_1^{4}x_2^{3}x_3^{3}-2x_1^{3}x_2^{2}x_3^{3}-x_1^{2}x_2x_3^{3}-x_1x_2x_3^{3}+x_1^{3}x_2^{2}x_3^{2}-x_1^{4}x_2^{4}x_3^{2}\\
&& 
-x_2^{3}x_3^{2}-2x_2^{2}x_3^{2}-4x_2x_3^{2}+3x_3^{2}+5x_2^{4}x_3^{3}-x_2^{3}x_3^{3}+x_2^{2}x_3^{3}+x_2x_3^{3}+4x_2^{6}x_3^{3}+9x_2^{5}x_3^{3}\\
&& 
+2x_2^{3}x_3^{4})\\
&& \frac{1}{(1-x_1)^2} \frac{1}{(1-x_2)^2} \frac{1}{(1-x_3)^5} \\ 
&&+(-x_1^{4}x_2^{6}x_3^{4}-2x_1^{3}x_2^{5}x_3^{4}-x_1^{2}x_2^{4}x_3^{4}-2x_1^{4}x_2^{5}x_3^{4}-4x_1^{3}x_2^{4}x_3^{4}-2x_1^{2}x_2^{3}x_3^{4}\\
&& 
-x_1^{4}x_2^{4}x_3^{4}-2x_1^{3}x_2^{3}x_3^{4}-x_1^{2}x_2^{2}x_3^{4}+2x_2^{8}x_3^{5}+8x_2^{7}x_3^{5}+12x_2^{6}x_3^{5}+8x_2^{5}x_3^{5}\\
&& 
+2x_2^{4}x_3^{5}-x_1x_2^{4}x_3^{4}-2x_1x_2^{3}x_3^{4}-x_1x_2^{2}x_3^{4}+4x_1x_3^{3}+x_1^{2}x_2^{5}x_3^{3}-x_1x_2^{4}x_3^{3}+4x_1^{2}x_2^{4}x_3^{3}\\
&& 
-4x_1x_2^{3}x_3^{3}+7x_1^{2}x_2^{3}x_3^{3}-5x_1x_2^{2}x_3^{3}+4x_1^{2}x_2^{2}x_3^{3}-2x_1x_2x_3^{3}-2x_1^{2}x_2x_3^{3}+x_1^{3}x_2^{5}x_3^{3}\\
&& 
+6x_1^{3}x_2^{4}x_3^{3}+5x_1^{3}x_2^{3}x_3^{3}+x_1^{4}x_2^{3}x_3^{4}+2x_1^{3}x_2^{2}x_3^{4}+x_1^{2}x_2x_3^{4}+x_1x_2x_3^{4}-x_1^{3}x_2^{2}x_3^{3}\\
&& 
+2x_1^{4}x_2^{5}x_3^{3}+2x_1^{4}x_2^{4}x_3^{3}+x_2^{4}x_3^{3}+3x_2^{3}x_3^{3}+3x_2^{2}x_3^{3}+5x_2x_3^{3}-4x_3^{3}-17x_2^{5}x_3^{4}-5x_2^{4}x_3^{4}\\
&& 
+2x_2^{3}x_3^{4}-2x_2^{2}x_3^{4}-x_2x_3^{4}-5x_2^{7}x_3^{4}-16x_2^{6}x_3^{4}-2x_2^{3}x_3^{5})\\
&& \frac{1}{(1-x_1)^2} \frac{1}{(1-x_2)} \frac{1}{(1-x_3)^6} \\ 
&&+(x_1^{2}x_2^{3}x_3^{4}-x_1x_2^{2}x_3^{4}) \frac{1}{(1-x_1x_2)^2} \frac{1}{(1-x_2x_3)} \frac{1}{(1-x_3)^6} \\ 
\end{eqnarray*}\begin{eqnarray*}
&&+(2x_1^{5}x_2^{8}x_3^{5}+8x_1^{4}x_2^{7}x_3^{5}+10x_1^{3}x_2^{6}x_3^{5}+2x_1^{4}x_2^{8}x_3^{5}+4x_1^{3}x_2^{7}x_3^{5}+2x_1^{2}x_2^{6}x_3^{5}\\
&& 
+4x_1^{5}x_2^{7}x_3^{5}+10x_1^{4}x_2^{6}x_3^{5}+8x_1^{3}x_2^{5}x_3^{5}+4x_1^{2}x_2^{5}x_3^{5}+2x_1^{5}x_2^{6}x_3^{5}+4x_1^{4}x_2^{5}x_3^{5}\\
&& 
+2x_1^{3}x_2^{4}x_3^{5}+2x_1^{2}x_2^{4}x_3^{5}-x_1^{5}x_2^{9}x_3^{6}-5x_1^{4}x_2^{8}x_3^{6}-7x_1^{3}x_2^{7}x_3^{6}-x_1^{4}x_2^{9}x_3^{6}\\
&& 
-2x_1^{3}x_2^{8}x_3^{6}-x_1^{2}x_2^{7}x_3^{6}-3x_1^{5}x_2^{8}x_3^{6}-9x_1^{4}x_2^{7}x_3^{6}-9x_1^{3}x_2^{6}x_3^{6}-3x_1^{2}x_2^{6}x_3^{6}\\
&& 
-3x_1^{5}x_2^{7}x_3^{6}-7x_1^{4}x_2^{6}x_3^{6}-5x_1^{3}x_2^{5}x_3^{6}-3x_1^{2}x_2^{5}x_3^{6}-x_1^{5}x_2^{6}x_3^{6}-2x_1^{4}x_2^{5}x_3^{6}\\
&& 
-x_1^{3}x_2^{4}x_3^{6}-x_1^{2}x_2^{4}x_3^{6}+2x_2^{11}x_3^{7}+10x_2^{10}x_3^{7}+20x_2^{9}x_3^{7}+20x_2^{8}x_3^{7}+10x_2^{7}x_3^{7}+2x_2^{6}x_3^{7}\\
&& 
-4x_2^{10}x_3^{6}-16x_2^{9}x_3^{6}-24x_2^{8}x_3^{6}-16x_2^{7}x_3^{6}-4x_2^{6}x_3^{6})\\
&& \frac{1}{(1-x_1)} \frac{1}{(1-x_2^{2}x_3)^2} \frac{1}{(1-x_3)^6} \\ 
&&+(x_1^{5}x_2^{8}x_3^{5}+4x_1^{4}x_2^{7}x_3^{5}+5x_1^{3}x_2^{6}x_3^{5}+x_1^{4}x_2^{8}x_3^{5}+2x_1^{3}x_2^{7}x_3^{5}+x_1^{2}x_2^{6}x_3^{5}\\
&& 
+2x_1^{5}x_2^{7}x_3^{5}+5x_1^{4}x_2^{6}x_3^{5}+4x_1^{3}x_2^{5}x_3^{5}+2x_1^{2}x_2^{5}x_3^{5}+x_1^{5}x_2^{6}x_3^{5}+2x_1^{4}x_2^{5}x_3^{5}\\
&& 
+x_1^{3}x_2^{4}x_3^{5}+x_1^{2}x_2^{4}x_3^{5}-x_2^{10}x_3^{6}-4x_2^{9}x_3^{6}-6x_2^{8}x_3^{6}-4x_2^{7}x_3^{6}-x_2^{6}x_3^{6})\\
&& \frac{1}{(1-x_1)} \frac{1}{(1-x_2^{2}x_3)^3} \frac{1}{(1-x_3)^5} \\ 
&&+(x_1^{5}x_2^{6}x_3^{3}+x_1^{4}x_2^{5}x_3^{3}-x_1^{3}x_2^{4}x_3^{3}+x_1^{4}x_2^{6}x_3^{3}+2x_1^{3}x_2^{5}x_3^{3}+x_1^{2}x_2^{4}x_3^{3}\\
&& 
-3x_2^{8}x_3^{4}-6x_2^{7}x_3^{4}-3x_2^{6}x_3^{4}-x_1x_2^{4}x_3^{3}-2x_1^{2}x_2^{4}x_3+2x_1x_2^{3}x_3-3x_1^{2}x_2^{3}x_3+3x_1x_2^{2}x_3\\
&& 
+3x_1^{2}x_2^{2}x_3-3x_1x_2x_3-x_1^{3}x_2^{4}x_3-x_1^{5}x_2^{5}x_3^{3}-2x_1^{4}x_2^{4}x_3^{3}-x_1^{3}x_2^{3}x_3^{3}-x_1^{2}x_2^{3}x_3^{3}\\
&& 
+x_1x_2^{3}x_3^{3}+x_1^{3}x_2^{3}x_3-x_2^{3}x_3-3x_2^{2}x_3+x_2x_3+x_2^{5}x_3^{3}+3x_2^{4}x_3^{3}-x_2^{3}x_3^{3}-x_2^{2}x_3^{2}+4x_2^{7}x_3^{3}\\
&& 
+4x_2^{6}x_3^{3}+3x_2^{5}x_3^{4}+2x_2^{5}x_3^{2}+x_2^{4}x_3^{2}-x_2^{3}x_3^{2})\\
&& \frac{1}{(1-x_1)} \frac{1}{(1-x_2)^4} \frac{1}{(1-x_3)^4} \\ 
\end{eqnarray*}\begin{eqnarray*}
&&+(-x_1^{5}x_2^{7}x_3^{4}-3x_1^{4}x_2^{6}x_3^{4}-3x_1^{3}x_2^{5}x_3^{4}-x_1^{4}x_2^{7}x_3^{4}-2x_1^{3}x_2^{6}x_3^{4}-x_1^{2}x_2^{5}x_3^{4}\\
&& 
-x_1^{5}x_2^{6}x_3^{4}-x_1^{4}x_2^{5}x_3^{4}+x_1^{3}x_2^{4}x_3^{4}-x_1^{2}x_2^{4}x_3^{4}+3x_2^{9}x_3^{5}+9x_2^{8}x_3^{5}+9x_2^{7}x_3^{5}\\
&& 
+3x_2^{6}x_3^{5}+x_1x_2^{5}x_3^{4}+x_1x_2^{4}x_3^{4}-x_1x_2^{2}x_3^{3}+2x_1^{2}x_2^{5}x_3^{2}-2x_1x_2^{4}x_3^{2}+5x_1^{2}x_2^{4}x_3^{2}\\
&& 
-5x_1x_2^{3}x_3^{2}+3x_1^{2}x_2^{3}x_3^{2}-3x_1x_2^{2}x_3^{2}-3x_1^{2}x_2^{2}x_3^{2}+3x_1x_2x_3^{2}+x_1^{3}x_2^{5}x_3^{2}+x_1^{3}x_2^{4}x_3^{2}\\
&& 
+x_1^{5}x_2^{5}x_3^{4}+2x_1^{4}x_2^{4}x_3^{4}+x_1^{3}x_2^{3}x_3^{4}+x_1^{2}x_2^{3}x_3^{4}-x_1x_2^{3}x_3^{4}-x_1^{3}x_2^{3}x_3^{2}+x_1^{5}x_2^{6}x_3^{3}\\
&& 
+2x_1^{4}x_2^{5}x_3^{3}+x_1^{3}x_2^{4}x_3^{3}+x_1^{4}x_2^{6}x_3^{3}+2x_1^{3}x_2^{5}x_3^{3}+x_1^{2}x_2^{4}x_3^{3}+x_2^{4}x_3^{2}+2x_2^{3}x_3^{2}\\
&& 
+4x_2^{2}x_3^{2}+2x_2^{2}x_3^{3}-7x_2^{6}x_3^{4}-2x_2^{5}x_3^{4}-4x_2^{4}x_3^{4}+x_2^{3}x_3^{4}-6x_2^{8}x_3^{4}-12x_2^{7}x_3^{4}+2x_2^{7}x_3^{3}\\
&& 
-3x_2^{5}x_3^{5}-3x_2^{5}x_3^{3}-x_2^{4}x_3^{3}+x_2^{3}x_3^{3})\\
&& \frac{1}{(1-x_1)} \frac{1}{(1-x_2)^3} \frac{1}{(1-x_3)^5} \\ 
&&+(x_1^{5}x_2^{8}x_3^{5}+4x_1^{4}x_2^{7}x_3^{5}+5x_1^{3}x_2^{6}x_3^{5}+x_1^{4}x_2^{8}x_3^{5}+2x_1^{3}x_2^{7}x_3^{5}+x_1^{2}x_2^{6}x_3^{5}\\
&& 
+2x_1^{5}x_2^{7}x_3^{5}+5x_1^{4}x_2^{6}x_3^{5}+4x_1^{3}x_2^{5}x_3^{5}+2x_1^{2}x_2^{5}x_3^{5}+x_1^{5}x_2^{6}x_3^{5}+x_1^{4}x_2^{5}x_3^{5}\\
&& 
-x_1^{3}x_2^{4}x_3^{5}+x_1^{2}x_2^{4}x_3^{5}-3x_2^{10}x_3^{6}-12x_2^{9}x_3^{6}-18x_2^{8}x_3^{6}-12x_2^{7}x_3^{6}-3x_2^{6}x_3^{6}\\
&& 
-x_1x_2^{6}x_3^{5}-2x_1x_2^{5}x_3^{5}-x_1x_2^{4}x_3^{5}+2x_1x_2^{2}x_3^{4}-2x_1^{2}x_2^{6}x_3^{3}+2x_1x_2^{5}x_3^{3}-7x_1^{2}x_2^{5}x_3^{3}\\
&& 
+7x_1x_2^{4}x_3^{3}-7x_1^{2}x_2^{4}x_3^{3}+8x_1x_2^{3}x_3^{3}-3x_1^{2}x_2^{3}x_3^{3}+3x_1x_2^{2}x_3^{3}+3x_1^{2}x_2^{2}x_3^{3}-2x_1x_2x_3^{3}\\
&& 
-x_1^{3}x_2^{6}x_3^{3}+x_1^{4}x_2^{6}x_3^{3}+2x_1^{4}x_2^{5}x_3^{3}-x_1^{5}x_2^{5}x_3^{5}-2x_1^{4}x_2^{4}x_3^{5}-x_1^{3}x_2^{3}x_3^{5}\\
&& 
-x_1^{2}x_2^{3}x_3^{5}+x_1x_2^{3}x_3^{5}+x_1^{3}x_2^{3}x_3^{3}+x_1^{5}x_2^{6}x_3^{3}-2x_1^{5}x_2^{7}x_3^{4}-6x_1^{4}x_2^{6}x_3^{4}\\
&& 
-6x_1^{3}x_2^{5}x_3^{4}-2x_1^{4}x_2^{7}x_3^{4}-4x_1^{3}x_2^{6}x_3^{4}-2x_1^{2}x_2^{5}x_3^{4}-2x_1^{5}x_2^{6}x_3^{4}-4x_1^{4}x_2^{5}x_3^{4}\\
&& 
-2x_1^{3}x_2^{4}x_3^{4}-2x_1^{2}x_2^{4}x_3^{4}-x_2^{5}x_3^{3}-3x_2^{4}x_3^{3}-3x_2^{3}x_3^{3}-5x_2^{2}x_3^{3}-2x_2x_3^{3}+25x_2^{7}x_3^{5}\\
&& 
+11x_2^{6}x_3^{5}+3x_2^{5}x_3^{5}+5x_2^{4}x_3^{5}-x_2^{3}x_3^{5}-3x_2^{2}x_3^{4}+8x_2^{9}x_3^{5}+24x_2^{8}x_3^{5}-5x_2^{8}x_3^{4}-8x_2^{7}x_3^{4}\\
&& 
+3x_2^{5}x_3^{6}+4x_2^{5}x_3^{4}+x_2^{4}x_3^{4}-x_2^{3}x_3^{4})\\
&& \frac{1}{(1-x_1)} \frac{1}{(1-x_2)^2} \frac{1}{(1-x_3)^6} \\ 
\end{eqnarray*}\begin{eqnarray*}
&&+(-x_1^{5}x_2^{9}x_3^{6}-5x_1^{4}x_2^{8}x_3^{6}-7x_1^{3}x_2^{7}x_3^{6}-x_1^{4}x_2^{9}x_3^{6}-2x_1^{3}x_2^{8}x_3^{6}-x_1^{2}x_2^{7}x_3^{6}\\
&& 
-3x_1^{5}x_2^{8}x_3^{6}-9x_1^{4}x_2^{7}x_3^{6}-9x_1^{3}x_2^{6}x_3^{6}-3x_1^{2}x_2^{6}x_3^{6}-3x_1^{5}x_2^{7}x_3^{6}-7x_1^{4}x_2^{6}x_3^{6}\\
&& 
-5x_1^{3}x_2^{5}x_3^{6}-3x_1^{2}x_2^{5}x_3^{6}-x_1^{5}x_2^{6}x_3^{6}-x_1^{4}x_2^{5}x_3^{6}+x_1^{3}x_2^{4}x_3^{6}-x_1^{2}x_2^{4}x_3^{6}\\
&& 
+3x_2^{11}x_3^{7}+15x_2^{10}x_3^{7}+30x_2^{9}x_3^{7}+30x_2^{8}x_3^{7}+15x_2^{7}x_3^{7}+3x_2^{6}x_3^{7}+x_1x_2^{7}x_3^{6}+3x_1x_2^{6}x_3^{6}\\
&& 
+3x_1x_2^{5}x_3^{6}+x_1x_2^{4}x_3^{6}-3x_1x_2^{2}x_3^{5}+2x_1^{2}x_2^{7}x_3^{4}-2x_1x_2^{6}x_3^{4}+9x_1^{2}x_2^{6}x_3^{4}-9x_1x_2^{5}x_3^{4}\\
&& 
+12x_1^{2}x_2^{5}x_3^{4}-15x_1x_2^{4}x_3^{4}+8x_1^{2}x_2^{4}x_3^{4}-11x_1x_2^{3}x_3^{4}+3x_1^{2}x_2^{3}x_3^{4}-3x_1x_2^{2}x_3^{4}\\
&& 
-3x_1^{2}x_2^{2}x_3^{4}-x_1^{3}x_2^{3}x_3^{4}+x_1^{3}x_2^{7}x_3^{4}-3x_1^{3}x_2^{6}x_3^{4}-6x_1^{3}x_2^{5}x_3^{4}-2x_1^{3}x_2^{4}x_3^{4}\\
&& 
+x_1^{5}x_2^{5}x_3^{6}+2x_1^{4}x_2^{4}x_3^{6}+x_1^{3}x_2^{3}x_3^{6}+x_1^{2}x_2^{3}x_3^{6}-x_1x_2^{3}x_3^{6}-3x_1^{5}x_2^{7}x_3^{4}\\
&& 
-9x_1^{4}x_2^{6}x_3^{4}-3x_1^{4}x_2^{7}x_3^{4}-3x_1^{5}x_2^{6}x_3^{4}-6x_1^{4}x_2^{5}x_3^{4}+3x_1^{5}x_2^{8}x_3^{5}+12x_1^{4}x_2^{7}x_3^{5}\\
&& 
+15x_1^{3}x_2^{6}x_3^{5}+3x_1^{4}x_2^{8}x_3^{5}+6x_1^{3}x_2^{7}x_3^{5}+3x_1^{2}x_2^{6}x_3^{5}+6x_1^{5}x_2^{7}x_3^{5}+15x_1^{4}x_2^{6}x_3^{5}\\
&& 
+12x_1^{3}x_2^{5}x_3^{5}+6x_1^{2}x_2^{5}x_3^{5}+3x_1^{5}x_2^{6}x_3^{5}+6x_1^{4}x_2^{5}x_3^{5}+3x_1^{3}x_2^{4}x_3^{5}+3x_1^{2}x_2^{4}x_3^{5}\\
&& 
+x_2^{6}x_3^{4}+4x_2^{5}x_3^{4}+6x_2^{4}x_3^{4}+4x_2^{3}x_3^{4}+6x_2^{2}x_3^{4}+5x_2x_3^{4}-61x_2^{8}x_3^{6}-44x_2^{7}x_3^{6}-16x_2^{6}x_3^{6}\\
&& 
-4x_2^{5}x_3^{6}-6x_2^{4}x_3^{6}+x_2^{3}x_3^{6}+4x_2^{2}x_3^{5}-10x_2^{10}x_3^{6}-40x_2^{9}x_3^{6}+9x_2^{9}x_3^{5}+25x_2^{8}x_3^{5}\\
&& 
+20x_2^{7}x_3^{5}-3x_2^{5}x_3^{7}-5x_2^{5}x_3^{5}-x_2^{4}x_3^{5}+x_2^{3}x_3^{5})\\
&& \frac{1}{(1-x_1)} \frac{1}{(1-x_2)} \frac{1}{(1-x_3)^7} \\ 
&&+(-x_1x_3^{4}+x_3^{4}+x_2^{2}x_3^{5}-x_2x_3^{4}) \frac{1}{(1-x_1)^2} \frac{1}{(1-x_2x_3)^2} \frac{1}{(1-x_3)^5} \\ 
\end{eqnarray*}\begin{eqnarray*}
&&+(-3x_2^{6}x_3^{2}+3x_2^{5}x_3^{2}+x_2^{4}-x_2^{3}) \frac{1}{(1-x_1)} \frac{1}{(1-x_2)^6} \frac{1}{(1-x_3)^2} \\ 
&&+(2x_2^{4}x_3-2x_2^{3}x_3-x_2^{3}+x_2^{2}) \frac{1}{(1-x_1)^2} \frac{1}{(1-x_2)^5} \frac{1}{(1-x_3)^2} \\ 
&&+x_2^{3}x_3^{5} \frac{1}{(1-x_1x_2)} \frac{1}{(1-x_2x_3)^2} \frac{1}{(1-x_3)^6} \\ 
\end{eqnarray*}
\end{document}